\documentclass[a4paper,10pt,francais]{smfart}
\usepackage[english,francais]{babel}
\usepackage{amsmath,amssymb,amsthm}
\usepackage{amsfonts}
\usepackage[T1]{fontenc}
\usepackage[all]{xy}
\usepackage{smfthm}
\usepackage{smfenum}
\usepackage{bull}
\usepackage{vmargin}
\usepackage{cite}
\usepackage{pifont}
\usepackage{euscript}
\usepackage{verbatim}
\RequirePackage{mathrsfs}
\usepackage{frcursive}

\DeclareMathOperator{\Hom}{Hom}
\DeclareMathOperator{\GL}{GL}

\DeclareMathOperator{\spec}{Spec}
\DeclareMathOperator{\card}{card} 
\DeclareMathOperator{\codim}{codim} 
\DeclareMathOperator{\jet}{jet}

\DeclareMathOperator{\moy}{moy}
\DeclareMathOperator{\Gal}{Gal}
\DeclareMathOperator{\dist}{dist}
\DeclareMathOperator{\tr}{tr}

\theoremstyle{plain}

\newtheorem*{lemmesiegel}{Lemme de Siegel}

\newtheorem*{lemmeBV}{Lemme de Bombieri \& Vaaler}

\newtheorem*{lemmeSA}{Lemme de Siegel absolu}

\newtheorem*{lemmepv}{Lemme de Siegel approch{\'e}}

\newtheorem*{lemmepva}{Lemme de Siegel approch{\'e} absolu}

\newtheorem*{lemmeperturbation}{Lemme de perturbation}

\author{{\'E}ric Gaudron}
\address{Universit{\'e} Grenoble I, Institut Fourier\\ UMR $5582$, BP $74$, $38402$ Saint-Martin-d'H{\`e}res Cedex, France}
\email{Eric.Gaudron@ujf-grenoble.fr}
\urladdr{http://www-fourier.ujf-grenoble.fr/\~{}gaudron}
\title[Minorations simultan\'ees de formes lin\'eaires de logarithmes]{Minorations simultan\'ees de formes lin\'eaires de logarithmes de nombres alg\'ebriques}
\alttitle{Simultaneous lower bounds for linear forms in logarithms of algebraic numbers}
\date{{\rm\texttt{Version r{\'e}vis{\'e}e}}~: Lundi 21 mai 2012}
\begin{document}
\frontmatter

\begin{abstract}
Soit $p$ un nombre premier ou $p=\infty$ et $k$ un corps de nombres plong{\'e} dans $\mathbf{C}_{p}$. Soit $n\in\mathbf{N}\setminus\{0\}$ et $u_{1},\ldots,u_{n}\in\mathbf{C}_{p}$ tels que $e^{u_{j}}\in k$ pour tout $j\in\{1,\ldots,n\}$. Soit $(\beta_{i,j})$, $1\le i\le t$, $1\le j\le n$, une matrice $t\times n$ {\`a} coefficients dans $k$. Soit $(\beta_{1,0},\ldots,\beta_{t,0})\in k^{t}$. Posons $\Lambda_{i}:=\beta_{i,0}+\sum_{j=1}^{n}{\beta_{i,j}u_{j}}\in\mathbf{C}_{p}$ pour tout $i\in\{1,\ldots,t\}$. Nous obtenons des minorations de $\max{\{\vert\Lambda_{i}\vert_{p}\,;\ 1\le i\le t\}}$ explicites en tous les param{\`e}tres, lorsque ce maximum n'est pas nul. Ces minorations englobent de nombreux r{\'e}sultats ant{\'e}rieurs. La d{\'e}monstration repose sur la m{\'e}thode de Baker-Philippon-Waldschmidt, la r{\'e}duction d'Hirata-Kohno, le proc{\'e}d{\'e} de changement de variables de Chudnovsky, repens{\'e}s avec les outils modernes de la th{\'e}orie des pentes ad{\'e}liques (\emph{m{\'e}thode de la section auxiliaire}). Au passage, nous montrons comment {\'e}tendre le lemme de Siegel absolu de Zhang au cadre des fibr{\'e}s ad{\'e}liques hermitiens et nous {\'e}tablissons un nouveau lemme de Siegel approch{\'e}.\end{abstract}
\begin{altabstract}This work falls within the theory of linear forms in logarithms over a commutative \emph{linear} algebraic group defined over a number field. We give lower bounds for simultaneous linear forms in logarithms of algebraic numbers, treating both the archimedean and $p$-adic cases. The proof includes Baker's method, Hirata's reduction, Chudnovsky's process of variable change. The novelty is that we integrate into the proof the modern tools of adelic slope theory (\emph{auxiliary section method}), using a new small values Siegel's lemma. \end{altabstract}


\altkeywords{Linear forms in logarithms, simultaneous approximations, Baker's method, Hirata's reduction, Chudnovsky's process of variable change, adelic slope, hermitian vector bundle, auxiliary section method, absolute Siegel's lemma, small values Siegel's lemma.}
\maketitle

\tableofcontents

\mainmatter

\section{Introduction}
Ce texte pr{\'e}sente une minoration simultan{\'e}e de formes lin{\'e}aires de logarithmes de nombres alg{\'e}briques, valide pour un plongement quelconque (complexe ou ultram{\'e}trique) du corps de nombres ambiant. En d'autres termes il s'agit de minorer $\max{\{\vert\Lambda_{i}\vert_{p_{0}}\,;\ 1\le i\le t\}}$ o{\`u} $\Lambda_{i}$ est de la forme $\beta_{i,0}+\sum_{j=1}^{n}{\beta_{i,j}u_{j}}$ avec $\beta_{i,j}$ et $e^{u_{j}}$ dans un corps de nombres $k$ plong{\'e} dans $\mathbf{C}_{p_{0}}$ et $\vert\cdot\vert_{p_{0}}$ est une valeur absolue sur le compl{\'e}t{\'e} $\mathbf{C}_{p_{0}}$ d'une cl{\^o}ture alg{\'e}brique de $\mathbf{Q}_{p_{0}}$ ($p_{0}$ est un nombre premier ou $p_{0}=\infty$ avec la notation $\mathbf{C}_{\infty}:=\mathbf{C}$) .\par Depuis les travaux de Baker entrepris au milieu des ann{\'e}es soixante, de nombreux articles ont {\'e}t{\'e} consacr{\'e}s {\`a} l'{\'e}tude de cas particuliers de ce probl{\`e}me. En tr{\`e}s grande majorit{\'e} les auteurs s'int{\'e}ressaient {\`a} \emph{une} seule forme lin{\'e}aire ($t=1$) et, souvent, seulement au cas d'un plongement archim{\'e}dien. Si le plongement {\'e}tait ultram{\'e}trique alors ils supposaient $\beta_{i,j}\in\mathbf{Q}$. Plus pr{\'e}cis{\'e}ment, d{\`e}s que $t\ge 2$, et si l'on {\'e}carte les travaux de Philippon \& Waldschmidt~\cite{pphmiw2} et d'Hirata-Kohno~\cite{hirata3} qui traitent le cas tr{\`e}s g{\'e}n{\'e}ral d'un groupe alg{\'e}brique commutatif quelconque, seuls Ramachandra~\cite{ramachandra} et Loxton~\cite{loxton} ont {\'e}tudi{\'e} la question qui nous int{\'e}resse ici (cas archim{\'e}dien). Lorsque le plongement est ultram{\'e}trique, ne semble exister que l'{\'e}tude de Dong~\cite{dong}, qui traite un cas particulier de la question ($\beta_{i,j}\in\mathbf{Z}$, $\beta_{i,0}=0$).  Ainsi, en d{\'e}pit de la richesse de la litt{\'e}rature sur le th{\`e}me des formes lin{\'e}aires de logarithmes\footnote{Nous renvoyons le lecteur au livre tr{\`e}s complet de Waldschmidt~\cite{miw4} pour en saisir l'{\'e}tendue.}, l'on ne trouve gu{\`e}re de r{\'e}sultats g{\'e}n{\'e}raux qui prennent en compte plusieurs formes lin{\'e}aires et aucun qui ne consid{\`e}re {\`a} la fois le cas archim{\'e}dien et le cas $p$-adique. Un des objectifs de cet article est de combler cette lacune, tout en donnant une minoration qui soit en phase avec celles pr{\'e}sent{\'e}es dans~\cite{miw4}.\par Si le r{\'e}sultat le plus g{\'e}n{\'e}ral que nous avons ne sera donn{\'e} qu'au paragraphe suivant, nous en proposons ici un cas particulier. {\'E}tant donn{\'e} un nombre alg{\'e}brique $\alpha$, on d{\'e}signe par $h(\alpha)$ la hauteur de Weil logarithmique absolue de $\alpha$. Si $x$ est un nombre r{\'e}el alors $[x]$ est la partie enti{\`e}re de $x$.
\begin{theo}\label{theointro} Soit $n\in\mathbf{N}\setminus\{0\}$. Soit $k$ un sous-corps de nombres de $\mathbf{C}$, de degr{\'e} $D$ sur $\mathbf{Q}$. Soit $t\in\{1,\ldots,n\}$ et \begin{equation*}(\beta_{i,j})_{\genfrac{}{}{0pt}{}{1\le i\le t}{1\le j\le n}}\in\mathrm{M}_{t,n}(k)\end{equation*}une matrice de rang \emph{maximal} $t$. Soit $(\beta_{i,0})_{1\le i\le t}\in k^{t}$. Pour tout $j\in\{1,\ldots,n\}$, soit $u_{j}\in\mathbf{C}$ tel que $\alpha_{j}:=e^{u_{j}}\in k$. Soit $b,a,\mathfrak{e}$ des nombres r{\'e}els positifs tels que $\mathfrak{e}\ge e$, \begin{equation*}\log a\ge\max_{1\le j\le n}{\left\{h(\alpha_{j}),\frac{\mathfrak{e}\vert u_{j}\vert}{D}\right\}},\qquad\log b\ge D\max_{\genfrac{}{}{0pt}{}{1\le i\le t}{0\le j\le n}}{\{1,h(\beta_{i,j})\}}.\end{equation*}Soit $\mathfrak{a}$ l'entier d{\'e}fini par \begin{equation*}\mathfrak{a}:=\left[\frac{D}{\log\mathfrak{e}}\log\left(e+\frac{D}{\log\mathfrak{e}}+\log a\right)\right]+1.\end{equation*}Si $\{u_{1},\ldots,u_{n}\}$ est une famille libre sur $\mathbf{Q}$ alors les formes lin{\'e}aires de logarithmes \begin{equation*}\Lambda_{i}:=\beta_{i,0}+\beta_{i,1}u_{1}+\cdots+\beta_{i,n}u_{n}\quad (1\le i\le t)\end{equation*} ne sont pas toutes nulles et elles v{\'e}rifient la minoration \begin{equation}\label{ineq:minorationstandard}\log\max_{1\le i\le t}{\vert\Lambda_{i}\vert}\ge -(4n)^{91n^{2}}\mathfrak{a}^{1/t}(\log b+\mathfrak{a}\log\mathfrak{e})\left(1+\frac{D\log a}{\log\mathfrak{e}}\right)^{n/t}.\end{equation}
\end{theo}Un {\'e}nonc{\'e} similaire sera aussi d{\'e}montr{\'e} dans le cas ultram{\'e}trique lorsque $k$ est un sous-corps de $\mathbf{C}_{p_{0}}$. Dans le minorant de~\eqref{ineq:minorationstandard}, on notera la d{\'e}pendance lin{\'e}aire (et donc optimale) en $\log b$ (hauteur des formes lin{\'e}aires) et la d{\'e}pendance usuelle $(\log a)^{n/t}(\log\log a)^{1+1/t}$ en $\log a$ (hauteur du point $(\alpha_{1},\ldots,\alpha_{n})$). Hormis l'aspect g{\'e}n{\'e}ral de cette minoration d{\'e}j{\`a} {\'e}voqu{\'e} ($t$ et $k\hookrightarrow\mathbf{C}_{p_{0}}$ quelconques), il n'y a plus le discriminant du corps de nombres, qui intervenait auparavant dans le cas simultan{\'e} (comme dans l'article~\cite{hirata3} d'Hirata-Kohno).
\par Cet article de synth{\`e}se n'apporte pas d'id{\'e}e originale qui am{\'e}liorerait de mani{\`e}re significative les r{\'e}sultats connus pour \emph{une} forme lin{\'e}aire. En revanche, la d{\'e}monstration que nous proposons est, elle, plus novatrice dans sa forme car elle combine la trame de la m{\'e}thode des fonctions auxiliaires classique avec les outils de la th{\'e}orie des pentes ad{\'e}liques (sans \og m{\'e}thode des pentes\fg{} proprement dite). Il en r{\'e}sulte un cumul des avantages propres {\`a} ces techniques~: souplesse d'utilisation et simplicit{\'e} de la d{\'e}marche pour l'une, conservation de l'aspect intrins{\`e}que des donn{\'e}es et obtention ais{\'e}e des constantes num{\'e}riques pour l'autre. Nous expliquerons plus en d{\'e}tail notre approche au \S~\ref{sec:canevasdelademonstration}. Si nous avons pris le parti d'{\'e}crire ce texte dans le cas d'un groupe lin{\'e}aire, il est cependant possible de g{\'e}n{\'e}raliser {\`a} un groupe alg{\'e}brique commutatif quelconque et m{\^e}me d'obtenir des constantes num{\'e}riques explicites pour une vari{\'e}t{\'e} ab{\'e}lienne, comme dans~\cite{artepredeux}, \emph{sans hypoth{\`e}se} de non-torsion du point rationnel consid{\'e}r{\'e} (ici $(\alpha_{1},\ldots,\alpha_{n})\in\mathbf{G}_{\mathrm{m}}^{n}(k)$).  
\subsection*{Remerciements}Je remercie Daniel Bertrand et Michel Waldschmidt de leurs commentaires sur une premi{\`e}re version de ce texte. Je remercie {\'e}galement Ga{\"e}l R{\'e}mond pour les discussions {\'e}clairantes que nous avons eues sur les parties~\ref{sec:elements} et~\ref{lemmedepetitesvaleurs}. Par ailleurs, la lecture attentive et critique du rapporteur a permis d'am{\'e}liorer nombre de d{\'e}tails. Je lui en suis extr{\^{e}}mement reconnaissant. Mentionnons enfin que ce texte a b{\'e}n{\'e}fici{\'e} du soutien de l'ANR Diophante, 06-JCJC-0028, pilot{\'e} par Lucia Di Vizio.
\section{{\'E}nonc{\'e} principal}
\subsection{}
\label{subsection:donnees}
On note $\EuScript{P}$ l'ensemble des nombres premiers auquel est adjoint un symbole suppl{\'e}mentaire not{\'e} $\infty$. Pour $p\ne+\infty$, on note $(\mathbf{C}_{p},\vert\cdot\vert_{p})$ le corps valu{\'e} ultram{\'e}trique complet et alg{\'e}briquement clos usuel, avec la normalisation $\vert p\vert_{p}=p^{-1}$ de la valeur absolue. Si $p=\infty$, le corps valu{\'e} $(\mathbf{C}_{\infty},\vert\cdot\vert_{\infty})$ est le corps des nombres complexes $\mathbf{C}$ muni de la valeur absolue usuelle. Soit $k$ un corps de nombres de degr{\'e} $D=[k:\mathbf{Q}]$. Pour chaque $p\in\EuScript{P}$ il y a $D$ morphismes de corps (appel{\'e}s plongements) $\sigma:k\to\mathbf{C}_{p}$ qui correspondent aux racines dans $\mathbf{C}_{p}$ du polyn{\^{o}}me minimal sur $\mathbf{Q}$ d'un {\'e}l{\'e}ment primitif de $k$. Avec ces conventions, la formule du produit s'{\'e}crit $$\forall x\in k\setminus\{0\},\quad\prod_{p\in\EuScript{P}}{\prod_{\sigma:k\hookrightarrow\mathbf{C}_{p}}{\vert\sigma(x)\vert_{p}}}=1$$et la \emph{hauteur de Weil} (logarithmique absolue) d'un {\'e}l{\'e}ment $x$ de $k$ est la somme (finie) $$h(x):=\frac{1}{D}\sum_{p\in\EuScript{P}}{\sum_{\sigma:k\hookrightarrow\mathbf{C}_{p}}{\log\max{\{1,\vert\sigma(x)\vert_{p}\}}}}.$$Soit $r_{p}:=p^{-1/(p-1)}$ si $p\ne\infty$ et $r_{\infty}:=+\infty$ sinon. D{\'e}signons par $\mathcal{T}_{p}$ le disque ouvert de $\mathbf{C}_{p}$, centr{\'e} en $0$ et de rayon $r_{p}$ ($\mathcal{T}_{\infty}=\mathbf{C}$). L'exponentielle du groupe de Lie $\mathbf{G}_{\mathrm{m}}(\mathbf{C}_{p})$ est d{\'e}finie sur $\mathcal{T}_{p}$ par la s{\'e}rie convergente usuelle $$\forall\,z\in\mathcal{T}_{p},\quad e^{z}=\sum_{i=0}^{\infty}{\frac{z^{i}}{i!}}\cdotp$$Dans la suite on fixe $p_{0}\in\EuScript{P}$ et $\sigma_{0}:k\hookrightarrow\mathbf{C}_{p_{0}}$ un plongement avec lequel nous identifierons $k$ {\`a} un sous-corps de $\mathbf{C}_{p_{0}}$. Soit $n\in\mathbf{N}\setminus\{0\}$ et $u_{1},\ldots,u_{n}$ des {\'e}l{\'e}ments de $\mathcal{T}_{p_{0}}$, d'exponentielles respectives $\alpha_{1},\ldots,\alpha_{n}$. Nous supposons que tous les $\alpha_{i}$ appartiennent {\`a} $k$.\par Soit $t\in\{1,\ldots,n\}$ et $\ell_{1},\ldots,\ell_{t}:k^{n}\to k$ des formes lin{\'e}aires que l'on {\'e}crit dans la base canonique sous la forme $$\forall\,z=(z_{1},\ldots,z_{n})\in k^{n},\quad \ell_{i}(z)=\beta_{i,1}z_{1}+\cdots+\beta_{i,n}z_{n}.$$Soit $W_{0}$ le sous-espace vectoriel de $k^{n}$ intersection des noyaux des formes lin{\'e}aires $\ell_{i}$, $i\in\{1,\ldots,t\}$. Soit $u_{0}:=-(\beta_{1,0},\ldots,\beta_{t,0})\in k^{t}$ et $u:=(u_{1},\ldots,u_{n})\in\mathbf{C}_{p_{0}}^{n}$. Pour tout $i\in\{1,\ldots,t\}$, posons $\Lambda_{i}:=\beta_{i,0}+\ell_{i}(u)$. Soit $\epsilon_{0}:=0$ si $p_{0}$ est premier et $\epsilon_{0}:=1$ si $p_{0}=\infty$.
\begin{theo}
\label{theoremeprincipal}Avec les donn{\'e}es ci-dessus, notons $\mathsf{T}_{u}$ le sous-espace vectoriel de dimension minimale de $\mathbf{Q}^{n}$ tel que $u\in\mathsf{T}_{u}\otimes_{\mathbf{Q}}\mathbf{C}_{p_{0}}$. Consid{\'e}rons une partie $I$ de $\{1,\ldots,n\}$ telle que $(u_{i})_{i\in I}$ soit une famille libre sur $\mathbf{Q}$ et maximale pour cette propri{\'e}t{\'e}. Nous supposons que $\vert u_{i}\vert_{p_{0}}\le r_{p_{0}}^{2}e^{-1}$ pour tout $i\in I$. Soit $a_{1},\ldots,a_{n},\mathfrak{e},\mathfrak{r}$ des nombres r{\'e}els v{\'e}rifiant les conditions suivantes~: \begin{equation*}\mathfrak{e}\in\begin{cases}[e,\min{\{r_{p_{0}}^{2}/\vert u_{i}\vert_{p_{0}};\ i\in I\}}] & \text{si $p_{0}\ne\infty$},\\ [e,+\infty[ & \text{si $p_{0}=\infty$},\end{cases}\end{equation*}
 \begin{equation*}\forall\,j\in\{1,\ldots,n\},\quad \log a_{j}\ge\max{\left\{h(\alpha_{j}),\frac{\epsilon_{0}\mathfrak{e}\vert u_{j}\vert_{p_{0}}}{D}\right\}}\cdotp\end{equation*}Soit $\mathfrak{a}$ l'entier d{\'e}fini par la formule \begin{equation*}\mathfrak{a}:=\left[\frac{D}{\log\mathfrak{e}}\log\left(e+\frac{D}{\log\mathfrak{e}}+\log\prod_{j=1}^{n}{a_{j}}\right)\right]+1.\end{equation*}Soit $b$ un nombre r{\'e}el v{\'e}rifiant \begin{equation*}\log b\ge D\max_{\genfrac{}{}{0pt}{}{1\le j\le t}{0\le\ell\le n}}{\{1,h(\beta_{j,\ell})\}}\end{equation*}et $s:=\dim\mathsf{T}_{u}-\dim(W_{0}\cap\mathsf{T}_{u})$. Alors, s'il existe $j\in\{1,\ldots, t\}$ tel que $\Lambda_{j}\ne 0$, on a \begin{equation*}\log\max_{1\le j\le t}{\vert\Lambda_{j}\vert_{p_{0}}}\ge-(4n)^{91n^{2}}\mathfrak{a}^{1/s}(\log b+\mathfrak{a}\log\mathfrak{e}+D\log\log\mathfrak{e})\prod_{i\in I}{\left(1+\frac{D\log a_{i}}{\log\mathfrak{e}}\right)^{1/s}}.\end{equation*}
\end{theo}
Le th{\'e}or{\`e}me de W{\"u}stholz~\cite{wustholz} affirme que $\mathsf{T}_{u}$ est l'espace tangent d'un sous-groupe alg{\'e}brique $G_{u}$ de $\mathbf{G}_{\mathrm{m}}^{n}$. Le param{\`e}tre $s$ est donc identique {\`a} celui de~\cite{davidhiratacrelles}, article dans lequel il appara{\^\i}t pour la premi{\`e}re fois. On notera au passage que $\dim\mathsf{T}_{u}=\card I$. Par ailleurs, nous verrons au \S~\ref{paragraphereductions} que si $s=1$, la quantit{\'e} $\log b+\mathfrak{a}\log\mathfrak{e}+D\log\log\mathfrak{e}$ peut {\^{e}}tre remplac{\'e}e par $\log(b\mathfrak{e})+D\log\mathfrak{a}+D\log\log\mathfrak{e}$.
\subsection{Commentaires}
\begin{enumerate}\item[1)] La diff{\'e}rence avec le r{\'e}sultat de Loxton~\cite{loxton} est assez modeste. Elle r{\'e}side surtout dans la d{\'e}pendance en le degr{\'e} $D$ (ici de l'ordre de $D^{1+\varepsilon+(n+1)/s}$ au lieu de $D^{200n}$) et dans le caract{\`e}re plus g{\'e}n{\'e}ral du th{\'e}or{\`e}me~\ref{theoremeprincipal} (plongement $\sigma_{0}$ quelconque, pas d'hypoth{\`e}se d'ind{\'e}pendance multiplicative des $\alpha_{j}$, $1\le j\le n$). En revanche  Loxton a une constante {\'e}gale {\`a} $(16n)^{200n}$. 
\item[2)] Lorsque $p_{0}$ est premier, je n'ai pas trouv{\'e} dans la litt{\'e}rature de minoration analogue {\`a} celle du th{\'e}or{\`e}me~\ref{theoremeprincipal}, bien que la m{\'e}thode de Baker originelle puisse donner un tel r{\'e}sultat. Il semble que l'obtention de minorations g{\'e}n{\'e}rales et explicites ait {\'e}t{\'e} d{\'e}laiss{\'e}e au profit de l'{\'e}tude de cas particuliers comme le cas rationnel ($\beta_{i,j}\in\mathbf{Z}$, $\beta_{i,0}=0$), \textit{a priori} plus utiles pour les applications (voir le survol~\cite{yu} et~\cite[p.~$547$]{miw4}). 
\item[3)] Toujours dans le cas $p$-adique, on observe que la minoration de $\log\max_{i}{\vert\Lambda_{i}\vert_{p_{0}}}$ ne d{\'e}pend pas de $p_{0}$. En r{\'e}alit{\'e} cette d{\'e}pendance est cach{\'e}e dans la condition $\vert u_{i}\vert_{p_{0}}\le r_{p_{0}}^{2}e^{-1}$ qui est assez forte. Nous aurions pu affaiblir la condition en $\vert u_{i}\vert_{p_{0}}<r_{p_{0}}^{2}$ mais au prix d'une complexit{\'e} accrue du minorant de $\log\max_{i}{\vert\Lambda_{i}\vert_{p_{0}}}$, {\`a} l'instar du th{\'e}or{\`e}me de R{\'e}mond \& Urfels~\cite{gu} (qui est dans un contexte un peu diff{\'e}rent). Nous avons donc opt{\'e} pour un compromis entre une hypoth{\`e}se un peu plus forte et un r{\'e}sultat un peu plus lisible et maniable. Dans le cas rationnel, Yu parvient {\`a} des hypoth{\`e}ses tr{\`e}s faibles du type $\vert u_{i}\vert_{p_{0}}\le 1$.
\item[4)] La constante $(4n)^{91n^{2}}$ qui appara{\^{\i}}t dans le th{\'e}or{\`e}me~\ref{theoremeprincipal} a une moins bonne d{\'e}pendance en $n$ que celle de Loxton car nous n'avons pas mis en {\oe}uvre la m{\'e}thode de Kummer. Rappelons que dans le cas rationnel, avec $t=1$ et $p_{0}=\infty$, Matveev a montr{\'e} que l'on pouvait avoir une constante du type $c^{n}$ avec $c$ absolue (voir~\cite{matveev}). Il serait int{\'e}ressant d'incorporer sa m{\'e}thode dans le cadre de la th{\'e}orie des pentes que nous utilisons ici. Quant {\`a} l'aspect num{\'e}rique de la constante $(4n)^{91n^{2}}$, une grand partie de sa taille vient du fait qu'il n'a pas {\'e}t{\'e} possible de supposer $n\ge 2$ dans la d{\'e}monstration ({\`a} cause du cas $p$-adique, voir plus haut). Une meilleure constante serait facilement accessible si l'on traitait s{\'e}par{\'e}ment les cas $n=1$ et $n\ge 2$ par exemple. 
 \end{enumerate}
 
 \par La d{\'e}monstration du th{\'e}or{\`e}me~\ref{theoremeprincipal} commence au \S~\ref{demonstrationtheoprincipal}, avec une r{\'e}duction au cas $\mathsf{T}_{u}=\mathbf{Q}^{n}$ (la famille $\{u_{1},\ldots,u_{n}\}$ est libre sur $\mathbf{Q}$). Nous expliquerons son d{\'e}roulement {\`a} cette occasion (\S~\ref{sec:canevasdelademonstration}) mais, auparavant, nous rappelons quelques notions de th{\'e}orie des pentes ad{\'e}liques et nous {\'e}non{\c{c}}ons un nouveau lemme de Siegel approch{\'e} qui jouera un r{\^{o}}le clef dans la d{\'e}monstration.
\section{{\'E}l{\'e}ments de th{\'e}orie des pentes  ad{\'e}liques}
\label{sec:elements}
La th{\'e}orie des pentes des fibr{\'e}s vectoriels ad{\'e}liques sur un corps de nombres $k$ a fait l'objet d'une pr{\'e}sentation syst{\'e}matique dans l'article~\cite{rendiconti}. Nous adoptons ici une pr{\'e}sentation un peu d{\'e}cal{\'e}e o{\`u} les places du corps de nombres sont remplac{\'e}es par les plongements dans les  corps $\mathbf{C}_{p}$ pour $p$ un nombre premier ou $p=\infty$ (dans ce dernier cas $\mathbf{C}_{\infty}$ d{\'e}signe le corps des nombres complexes $\mathbf{C}$). Rappelons que $\EuScript{P}$ d{\'e}signe l'ensemble des nombres premiers et de l'infini, $\EuScript{P}:=\{2,3,5,\ldots,\infty\}$.
\subsection{D{\'e}finitions} Un \emph{fibr{\'e} vectoriel pr{\'e}-ad{\'e}lique} $\overline{E}=(E,(\Vert\cdot\Vert_{\overline{E},\sigma})_{\sigma})$ est la donn{\'e}e d'un $k$-espace vectoriel $E$, de dimension finie $\nu\ge 0$, et, pour tout $p\in\EuScript{P}$ et tout plongement $\sigma:k\hookrightarrow\mathbf{C}_{p}$, d'une norme $\Vert\cdot\Vert_{\overline{E},\sigma}$ sur $E\otimes_{\sigma}\mathbf{C}_{p}$, qui satisfont aux contraintes suivantes~:\begin{enumerate}\item[1)] il existe une $k$-base $(e_{1},\ldots,e_{\nu})$ de $E$ telle que, pour tout $p\in\EuScript{P}$ sauf un nombre fini, pour tout plongement $\sigma:k\hookrightarrow\mathbf{C}_{p}$, on a $\Vert\sum_{i=1}^{\nu}{e_{i}\otimes_{\sigma}x_{i}}\Vert_{\overline{E},\sigma}=\max{\{\vert x_{1}\vert_{p},\ldots,\vert x_{\nu}\vert_{p}\}}$ pour tout $(x_{1},\ldots,x_{\nu})\in\mathbf{C}_{p}^{\nu}$,\item[2)] si $p\ne\infty$ alors $\Vert\cdot\Vert_{\overline{E},\sigma}$ est une ultranorme.
\end{enumerate}
Dans cette d{\'e}finition, la notation $E\otimes_{\sigma}\mathbf{C}_{p}$ d{\'e}signe le produit tensoriel $E\otimes_{k}\mathbf{C}_{p}$ o{\`u} $\mathbf{C}_{p}$ est muni de la structure de $k$-espace vectoriel donn{\'e}e par $\lambda.x:=\sigma(\lambda)x$ pour $\lambda\in k$ et $x\in\mathbf{C}_{p}$. Dans la suite, l'espace $E\otimes_{\sigma}\mathbf{C}_{p}$ est vu comme espace vectoriel sur $\mathbf{C}_{p}$ avec la loi externe $\mu_{1}.(x\otimes\mu_{2}):=x\otimes(\mu_{1}\mu_{2})$ pour $x\in E$ et $\mu_{1},\mu_{2}\in\mathbf{C}_{p}$. En particulier, pour $\lambda\in k$ et $x\in E\otimes_{\sigma}\mathbf{C}_{p}$, on a $\Vert\lambda.x\Vert_{\overline{E},\sigma}=\vert\sigma(\lambda)\vert_{p}\Vert x\Vert_{\overline{E},\sigma}$. On rappelle que deux plongements $\sigma,\sigma'$ de $k$ dans $\mathbf{C}_{p}$ sont dits \emph{conjugu{\'e}s} s'il existe un automorphisme continu $\iota:\mathbf{C}_{p}\to\mathbf{C}_{p}$ tel que $\sigma'=\iota\circ\sigma$. Ceci d{\'e}finit une relation d'{\'e}quivalence et une \emph{place} de $k$ est une classe d'{\'e}quivalence de plongements.
\begin{defi}\label{defiinvariance} Soit $p\in\EuScript{P}$ et $v$ une place de $k$ au-dessus de $p$. On dit que la collection de normes $(\Vert\cdot\Vert_{\overline{E},\sigma})_{\sigma\in v}$ est \emph{invariante par Galois} si, pour \emph{tout} automorphisme continu $\iota:\mathbf{C}_{p}\to\mathbf{C}_{p}$, l'application canonique $$\mathrm{id}\otimes\iota:(E\otimes_{\sigma}\mathbf{C}_{p},\Vert\cdot\Vert_{\overline{E},\sigma})\to(E\otimes_{\iota\circ\sigma}\mathbf{C}_{p},\Vert\cdot\Vert_{\overline{E},\iota\circ\sigma})$$est une isom{\'e}trie.
\end{defi}
 
 Quand la collection de normes $(\Vert\cdot\Vert_{\overline{E},\sigma})_{\sigma\in v}$ est invariante par Galois, elle ne d{\'e}pend que de la place $v$ de $k$ induite par $\sigma$ et il est naturel alors de la noter $\Vert\cdot\Vert_{\overline{E},v}$. Par exemple, si $E=k^{\nu}$ et si $\sigma:k\hookrightarrow\mathbf{C}$ est un plongement complexe, la paire $\{\Vert\cdot\Vert_{\overline{E},\sigma},\Vert\cdot\Vert_{\overline{E},\overline{\sigma}}\}$ (ou, par abus de langage, la norme $\Vert\cdot\Vert_{\overline{E},\sigma}$) est invariante par Galois si $\Vert(x_{1},\ldots,x_{\nu})\Vert_{\overline{E},\sigma}=\Vert(\overline{x_{1}},\ldots,\overline{x_{\nu}})\Vert_{\overline{E},\overline{\sigma}}$ pour tout $(x_{1},\ldots,x_{\nu})\in\mathbf{C}^{\nu}$, car il n'y a que deux automorphismes continus de $\mathbf{C}$ qui sont l'identit{\'e} et la conjugaison complexe. Si, de plus, $\sigma$ est r{\'e}el ($\sigma(k)\subseteq\mathbf{R}$) et si $\Vert\cdot\Vert_{\overline{E},\sigma}$ est donn{\'e} par une matrice hermitienne $A\in\mathrm{M}_{\nu}(\mathbf{C})$ d{\'e}finie positive (c'est-{\`a}-dire $\Vert x\Vert_{\overline{E},\sigma}=({}^{\mathrm{t}}\overline{x}Ax)^{1/2}$ pour tout $x\in\mathbf{C}^{\nu}$) alors n{\'e}cessairement $A$ est une matrice sym{\'e}trique \emph{r{\'e}elle}.
 
 \begin{defi}
 Un \emph{fibr{\'e} vectoriel ad{\'e}lique} sur $k$ est un fibr{\'e} vectoriel pr{\'e}-ad{\'e}lique $\overline{E}=(E,(\Vert\cdot\Vert_{\overline{E},\sigma})_{\sigma})$ sur $k$ tel que, pour toute place $v$ de $k$, la collection de normes $(\Vert\cdot\Vert_{\overline{E},\sigma})_{\sigma\in v}$ est invariante par Galois. \end{defi}On le note alors habituellement $\overline{E}=(E,(\Vert\cdot\Vert_{\overline{E},v})_{v})$ o{\`u} $v$ parcourt les places de $k$. C'est le point de vue adopt{\'e} dans~\cite{rendiconti}.\par Soit $\overline{E}$ un fibr{\'e} vectoriel pr{\'e}-ad{\'e}lique. Si $F$ est un sous-espace vectoriel de $E$, il induit un sous-fibr{\'e} vectoriel pr{\'e}-ad{\'e}lique $\overline{F}$, dont les normes sont les restrictions de celles de $\overline{E}$ {\`a} $F$. Le fibr{\'e} vectoriel pr{\'e}-ad{\'e}lique est dit \emph{hermitien} lorsque toutes les normes correspondant aux plongements complexes $\sigma:k\hookrightarrow\mathbf{C}$ sont hermitiennes. Il est dit \emph{pur}\footnote{La notion de puret{\'e} d'une norme est un cas particulier de la propri{\'e}t{\'e} (N) de Serre~\cite[p.~$69$]{serre}. La terminologie \emph{solid norm} est {\'e}galement parfois employ{\'e}e dans d'autres contextes~\cite[p.~$19$]{pgs}.} lorsque, pour tout plongement ultram{\'e}trique $\sigma:k\hookrightarrow\mathbf{C}_{p}$ et tout $x\in E$, on a $\Vert x\Vert_{\overline{E},\sigma}\in\vert k_{\sigma}\vert_{p}$ o{\`u} $k_{\sigma}$ est l'adh{\'e}rence topologique de $\sigma(k)$ dans $\mathbf{C}_{p}$ (voir~\cite{aneuf}). Contrairement {\`a} la d{\'e}finition~$2.3$ de~\cite{aneuf}, il n'est pas suppos{\'e} ici qu'un fibr{\'e} pr{\'e}-ad{\'e}lique hermitien est n{\'e}cessairement pur. Les fibr{\'e}s ad{\'e}liques hermitiens purs sont exactement les \emph{fibr{\'e}s vectoriels hermitiens sur $\spec\mathcal{O}_{k}$} de la g{\'e}om{\'e}trie d'Arakelov classique (voir proposition~\ref{proppure}). L'exemple le plus simple d'un tel fibr{\'e} est celui du fibr{\'e} standard $(k^{\nu},\vert\cdot\vert_{2})$ que l'on obtient de la mani{\`e}re suivante. \'Etant donn{\'e} $p\in\EuScript{P}$ et $x_{1},\ldots,x_{\nu}$ des nombres r{\'e}els, notons $$\moy_{p}(x_{1},\ldots,x_{\nu}):=\begin{cases}(\sum_{i=1}^{\nu}{x_{i}^{2}})^{1/2} & \text{si $p=\infty$},\\ \max{\{x_{1},\ldots,x_{\nu}\}} & \text{si $p\ne\infty$}.\end{cases}$$Le fibr{\'e} vectoriel ad{\'e}lique $(k^{\nu},\vert\cdot\vert_{2})$, hermitien et pur, est par d{\'e}finition le fibr{\'e} ad{\'e}lique d'espace sous-jacent $k^{\nu}$ et de normes $\vert x\vert_{2,\sigma}:=\moy_{p}(\vert x_{1}\vert_{p},\ldots,\vert x_{\nu}\vert_{p})$ pour tout $x=(x_{1},\ldots,x_{\nu})\in\mathbf{C}_{p}^{\nu}$. Apr{\`e}s le choix d'une $k$-base (quelconque) $(e_{1},\ldots,e_{\nu})$ de $E$, les fibr{\'e}s vectoriels hermitiens sur $\spec\,\mathcal{O}_{k}$ s'obtiennent {\`a} partir du fibr{\'e} standard $(k^{\nu},\vert\cdot\vert_{2})$ au moyen d'une collection de matrices $\{A_{\sigma}\in\GL_{\nu}(k_{\sigma})\,;\ \sigma:k\hookrightarrow\mathbf{C}_{p}\}_{p\in\EuScript{P}}$ telles que
\begin{enumerate}
\item[1)] pour tout $\sigma$, sauf un nombre fini, $A_{\sigma}$ est une isom{\'e}trie de $(\mathbf{C}_{p}^{\nu},\vert\cdot\vert_{2,\sigma})$,
\item[2)] pour tout automorphisme continu $\iota:\mathbf{C}_{p}\to\mathbf{C}_{p}$, on a $A_{\iota\circ\sigma}=\iota(A_{\sigma})$ (l'action de $\iota$ est sur les coefficients de la matrice),\item[3)] pour tous $x_{1},\ldots,x_{\nu}\in\mathbf{C}_{p}$ on a  $$\Vert e_{1}\otimes_{\sigma}x_{1}+\cdots+e_{\nu}\otimes_{\sigma}x_{\nu}\Vert_{\overline{E},\sigma}=\Big\vert A_{\sigma}\begin{pmatrix}x_{1}\\ \vdots\\ x_{\nu}\end{pmatrix}\Big\vert_{2,\sigma}.$$
\end{enumerate}

 \subsection{Extension des scalaires}\label{sectionextension} Soit $\overline{E}$ un fibr{\'e} vectoriel pr{\'e}-ad{\'e}lique sur $k$ et $K$ une extension finie de $k$ (on identifie $k$ {\`a} un sous-corps de $K$). Le $K$-espace vectoriel $E_{K}:=E\otimes_{k}K$ est muni naturellement d'une structure de fibr{\'e} vectoriel pr{\'e}-ad{\'e}lique $\overline{E_{K}}$ de la mani{\`e}re suivante. Soit $\tau:K\hookrightarrow\mathbf{C}_{p}$ un plongement qui {\'e}tend $\sigma:k\hookrightarrow\mathbf{C}_{p}$. Pour des {\'e}l{\'e}ments $e_{i}\in E$, $\lambda_{i}\in K$, $x_{i}\in\mathbf{C}_{p}$ posons $$\left\Vert\sum_{i}{(e_{i}\otimes_{k}\lambda_{i})\otimes_{\tau}x_{i}}\right\Vert_{\overline{E_{K}},\tau}:=\left\Vert\sum_{i}{e_{i}\otimes_{\sigma}(\tau(\lambda_{i})x_{i})}\right\Vert_{\overline{E},\sigma}.$$En choisissant des bases de $E$ sur $k$ et de $K$ sur $k$, on v{\'e}rifie que cette d{\'e}finition a bien un sens et que la norme d'un {\'e}l{\'e}ment $x\in(E\otimes_{k}K)\otimes_{\tau}\mathbf{C}_{p}$ ne d{\'e}pend pas du choix de sa repr{\'e}sentation en somme de tenseurs {\'e}l{\'e}mentaires. Le couple $\overline{E_{K}}=(E\otimes_{k}K,(\Vert\cdot\Vert_{\overline{E_{K}},\tau})_{\tau})$ forme alors un fibr{\'e} vectoriel pr{\'e}-ad{\'e}lique sur $K$. Il est ad{\'e}lique si $\overline{E}$ l'est car, si $\iota:\mathbf{C}_{p}\to\mathbf{C}_{p}$ est un automorphisme continu de corps, on a \begin{equation*}\begin{split}\left\Vert\sum_{i}{(e_{i}\otimes_{k}\lambda_{i})\otimes_{\iota\circ\tau}\iota(x_{i})}\right\Vert_{\overline{E_{K}},\iota\circ\tau}&=\left\Vert\sum_{i}{e_{i}\otimes_{\iota\circ\sigma}\left((\iota\circ\tau)(\lambda_{i})\iota(x_{i})\right)}\right\Vert_{\overline{E},\iota\circ\sigma}\quad\text{par d{\'e}finition}\\ & =\left\Vert\sum_{i}{e_{i}\otimes_{\iota\circ\sigma}\iota(\tau(\lambda_{i})x_{i})}\right\Vert_{\overline{E},\iota\circ\sigma}\\ &=\left\Vert\sum_{i}{e_{i}\otimes_{\sigma}(\tau(\lambda_{i})x_{i})}\right\Vert_{\overline{E},\sigma}\quad\text{par isom{\'e}trie}\\ & =\left\Vert\sum_{i}{(e_{i}\otimes_{k}\lambda_{i})\otimes_{\tau}x_{i}}\right\Vert_{\overline{E_{K}},\tau}\quad\text{par d{\'e}finition}. \end{split}\end{equation*}De m{\^{e}}me, si la norme $\Vert\cdot\Vert_{\overline{E},\sigma}$ est hermitienne alors $\Vert\cdot\Vert_{\overline{E_{K}},\tau}$ l'est aussi pour tout plongement $\tau$ qui {\'e}tend $\sigma$.

 \subsection{Produit tensoriel de fibr{\'e}s hermitiens} Soit $\overline{E}$ et $\overline{F}$ des fibr{\'e}s pr{\'e}-ad{\'e}liques hermitiens sur $k$. On munit $E\otimes_{k}F$ d'une structure de fibr{\'e} pr{\'e}-ad{\'e}lique hermitien de la mani{\`e}re suivante. Si $\sigma:k\hookrightarrow\mathbf{C}$ est un plongement complexe, on consid{\`e}re des bases orthonorm{\'e}es $\mathsf{e}_{1},\ldots,\mathsf{e}_{\nu}$ de $(E\otimes_{\sigma}\mathbf{C},\Vert\cdot\Vert_{\overline{E},\sigma})$ et $\mathsf{f}_{1},\ldots,\mathsf{f}_{\mu}$ de $(F\otimes_{\sigma}\mathbf{C},\Vert\cdot\Vert_{\overline{F},\sigma})$ et on choisit sur $(E\otimes_{k}F)\otimes_{\sigma}\mathbf{C}\simeq(E\otimes_{\sigma}\mathbf{C})\otimes_{\mathbf{C}}(F\otimes_{\sigma}\mathbf{C})$ l'unique produit hermitien pour lequel $(\mathsf{e}_{i}\otimes\mathsf{f}_{j})_{i,j}$ forme une base orthonorm{\'e}e (ind{\'e}pendant du choix des bases). Si $\sigma:k\hookrightarrow\mathbf{C}_{p}$ est un plongement ultram{\'e}trique, on identifie $E\otimes_{k}F$ {\`a} $\Hom_{k}(E^{\mathsf{v}},F)$ o{\`u} $E^{\mathsf{v}}=\Hom_{k}(E,k)$ d{\'e}signe le dual de $E$ et l'on munit $(E\otimes_{k}F)\otimes_{\sigma}\mathbf{C}_{p}$ de la norme d'op{\'e}rateur sur $\Hom_{k}(E^{\mathsf{v}},F)\otimes_{\sigma}\mathbf{C}_{p}$. Autrement dit, si $x=\sum_{i=1}^{N}{x_{i}\otimes y_{i}}$ avec $x_{i}\in E\otimes_{\sigma}\mathbf{C}_{p}$ et $y_{i}\in F\otimes_{\sigma}\mathbf{C}_{p}$, on a $$\Vert x\Vert_{\overline{E\otimes F},\sigma}:=\sup{\left\{\frac{\Vert\sum_{i=1}^{N}{\varphi(x_{i})y_{i}}\Vert_{\overline{F},\sigma}}{\Vert\varphi\Vert_{\overline{E^{\mathsf{v}}},\sigma}}\,;\ \varphi\in(E^{\mathsf{v}}\otimes_{\sigma}\mathbf{C}_{p})\setminus\{0\}\right\}}.$$On v{\'e}rifie que le couple $\overline{E}\otimes\overline{F}=(E\otimes_{k}F,(\Vert\cdot\Vert_{\overline{E\otimes F},\sigma})_{\sigma})$ forme un fibr{\'e} pr{\'e}-ad{\'e}lique hermitien. De plus si $\overline{E}$ et $\overline{F}$ sont des fibr{\'e}s ad{\'e}liques hermitiens (\emph{resp}. purs), il en est de m{\^{e}}me pour $\overline{E}\otimes\overline{F}$. Lorsque $\overline{E}$ est un fibr{\'e} pr{\'e}-ad{\'e}lique hermitien, de dimension $\nu\ge 1$, et $\mu\in\{1,\ldots,\nu\}$, on peut alors munir la puissance ext{\'e}rieure $\wedge^{\mu}E$ de $E$ de normes lui conf{\'e}rant une structure de fibr{\'e} pr{\'e}-ad{\'e}lique $\overline{\wedge^{\mu}E}$ sur $k$ de la mani{\`e}re suivante. Si $\sigma$ est un plongement complexe, on consid{\`e}re une base orthonorm{\'e}e $\mathsf{e}_{1},\ldots,\mathsf{e}_{\nu}$ de $(E\otimes_{\sigma}\mathbf{C},\Vert\cdot\Vert_{\overline{E},\sigma})$ et on d{\'e}cr{\`e}te que la base $\{\mathsf{e}_{i_{1}}\wedge\cdots\wedge\mathsf{e}_{i_{\mu}}\,;\ 1\le i_{1}<\cdots<i_{\mu}\le\nu\}$ de $\wedge^{\mu}E\otimes_{\sigma}\mathbf{C}$ est orthonorm{\'e}e. Si $\sigma$ est un plongement ultram{\'e}trique, on voit $\wedge^{\mu}E\otimes_{\sigma}\mathbf{C}_{p}$ comme un quotient de $(E\otimes_{\sigma}\mathbf{C}_{p})^{\otimes\mu}$ que l'on munit de la norme quotient induite par $\Vert\cdot\Vert_{\overline{E}^{\otimes\mu},\sigma}$. 
 \subsection{Pentes des fibr{\'e}s pr{\'e}-ad{\'e}liques hermitiens}On note $\det E$ la puissance ext{\'e}rieure maximale $\wedge^{\nu}E$ de $E$.
\begin{defi}
Le \emph{degr{\'e} d'Arakelov (normalis{\'e}) du fibr{\'e} pr{\'e}-ad{\'e}lique hermitien $\overline{E}$} est le nom\-bre r{\'e}el $$\widehat{\deg}_{\mathrm{n}}\overline{E}:=-\frac{1}{D}\sum_{p\in\EuScript{P}}{\sum_{\sigma:k\hookrightarrow\mathbf{C}_{p}}{\log\Vert e_{1}\wedge\cdots\wedge e_{\nu}\Vert_{\overline{\det E},\sigma}}}$$pour tout choix d'une $k$-base $(e_{1},\ldots,e_{\nu})$ de $E$. Lorsque $E=\{0\}$ on pose $\widehat{\deg}_{\mathrm{n}}\overline{E}:=0$. \end{defi}L'ind{\'e}pendance du degr{\'e} par rapport au choix de cette base est une cons{\'e}quence de la formule du produit. La \emph{hauteur (logarithmique absolue)} de $\overline{E}$ est $h(\overline{E}):=-\widehat{\deg}_{\mathrm{n}}\overline{E}$, sa \emph{pente d'Arakelov normalis{\'e}e} (si $E\ne\{0\}$) est \begin{equation*}\widehat{\mu}(\overline{E}):=\frac{\widehat{\deg}_{\mathrm{n}}\overline{E}}{\dim E},\end{equation*}sa \emph{pente maximale} (si $E\ne\{0\}$) est \begin{equation*}\widehat{\mu}_{\mathrm{max}}(\overline{E}):=\max{\left\{\widehat{\mu}(\overline{F})\,;\ \{0\}\ne F\subseteq E\right\}}.\end{equation*}On a {\'e}galement une notion de \emph{hauteur pour les {\'e}l{\'e}ments $x$ de $E$}~: $$\forall\,x\in E\setminus\{0\},\quad h_{\overline{E}}(x):=\frac{1}{D}\sum_{p\in\EuScript{P}}{\sum_{\sigma:k\hookrightarrow\mathbf{C}_{p}}{\log\Vert x\Vert_{\overline{E},\sigma}}}\qquad\text{($h_{\overline{E}}(0):=-\infty$)}.$$
\begin{lemm}Soit $\overline{E}$ un fibr{\'e} pr{\'e}-ad{\'e}lique hermitien sur $k$. 
\begin{enumerate}
\item[1)] Si $\overline{F}$ est un sous-fibr{\'e} de $\overline{E}$ alors on a $\widehat{\deg}_{\mathrm{n}}\overline{E/F}=\widehat{\deg}_{\mathrm{n}}\overline{E}-\widehat{\deg}_{\mathrm{n}}\overline{F}$.
\item[2)] Si $\overline{E_{1}},\overline{E_{2}}$ sont des sous-fibr{\'e}s de $\overline{E}$ alors on a 
 $$\widehat{\deg}_{\mathrm{n}}\overline{E_{1}}+\widehat{\deg}_{\mathrm{n}}\overline{E_{2}}\le\widehat{\deg}_{\mathrm{n}}\overline{E_{1}+E_{2}}+\widehat{\deg}_{\mathrm{n}}\overline{E_{1}\cap E_{2}}.$$
\end{enumerate}
\end{lemm}
\begin{proof} La d{\'e}monstration est identique {\`a} celle des fibr{\'e}s vectoriels hermitiens sur $\spec\mathcal{O}_{k}$ (voir~\cite[propositions~$4.22$ et $4.23$]{rendiconti}). 
\end{proof}
Toutes ces notions sont invariantes par extension des scalaires~:
\begin{prop}\label{propmumax}Soit $\overline{E}$ un fibr{\'e} pr{\'e}-ad{\'e}lique hermitien sur $k$ et $K/k$ une extension finie. Alors on a $\widehat{\deg}_{\mathrm{n}}\overline{E_{K}}=\widehat{\deg}_{\mathrm{n}}\overline{E}$ et $\widehat{\mu}_{\mathrm{max}}(\overline{E_{K}})=\widehat{\mu}_{\mathrm{max}}(\overline{E})$.
\end{prop}
\begin{proof}Le nombre de plongements $\tau:K\hookrightarrow\mathbf{C}_{p}$ qui prolongent $\sigma:k\hookrightarrow\mathbf{C}_{p}$ est {\'e}gal {\`a} $[K:k]$, ce qui donne la premi{\`e}re propri{\'e}t{\'e}. Pour la seconde, on montre comme dans le cas classique l'existence et l'unicit{\'e} du fibr{\'e} d{\'e}stabilisant (en utilisant le lemme pr{\'e}c{\'e}dent), c'est-{\`a}-dire qu'il existe un unique sous-fibr{\'e} $\overline{D_{k}}$ de $\overline{E}$, de dimension maximale, tel que $\widehat{\mu}_{\mathrm{max}}(\overline{E})=\widehat{\mu}(\overline{D_{k}})$. Comme $D_{k}\otimes_{k}K$ est un sous-espace de $E_{K}$, on a $$\widehat{\mu}_{\mathrm{max}}(\overline{E})=\widehat{\mu}(\overline{D_{k}})=\widehat{\mu}((\overline{D_{k}})_{K})\le\widehat{\mu}_{\mathrm{max}}(\overline{E_{K}}).$$Pour d{\'e}montrer l'{\'e}galit{\'e}, on peut donc supposer que $K/k$ est galoisienne. Consid{\'e}rons le morphisme de groupes $\rho:\Gal(K/k)\to\GL(E_{K})$ qui {\`a} $u\in\Gal(K/k)$ associe l'isomorphisme $\rho(u)=\mathrm{id}_{E}\otimes u$ de $E_{K}=E\otimes_{k}K$. Pour tout $m\in\{1,\ldots,\dim E\}$, pour tous $e_{1},\ldots,e_{m}\in E_{K}$, pour tout plongement $\tau:K\hookrightarrow\mathbf{C}_{p}$, on a $$\Vert\rho(u)(e_{1})\wedge\cdots\wedge\rho(u)(e_{m})\Vert_{\overline{\wedge^{m}E},\tau}=\Vert e_{1}\wedge\cdots\wedge e_{m}\Vert_{\overline{\wedge^{m}E},\tau\circ u}.$$En appliquant cette {\'e}galit{\'e} {\`a} une $K$-base de $D_{K}$, on obtient $\widehat{\mu}(\overline{\rho(u)(D_{K})})=\widehat{\mu}(\overline{D_{K}})$. Par unicit{\'e} de $\overline{D_{K}}$ on a $\rho(u)(D_{K})=D_{K}$ pour tout $u\in\Gal(K/k)$. Il existe alors un sous-espace vectoriel $D_{0}$ de $E$ tel que $D_{K}=D_{0}\otimes_{k}K$. En effet, si, pour $x\in E_{K}$, on pose $\tr(x):=\sum_{u\in\Gal(K/k)}{\rho(u)(x)}\in E$, alors $D_{0}=\{\tr(x)\,;\ x\in D_{K}\}$ convient. On conclut avec $\widehat{\mu}_{\mathrm{max}}(\overline{E_{K}})=\widehat{\mu}(\overline{D_{K}})=\widehat{\mu}(\overline{D_{0}})\le\widehat{\mu}_{\mathrm{max}}(\overline{E})$.
\end{proof}
L'argument sur le nombre de plongements $\tau:K\hookrightarrow\mathbf{C}_{p}$ qui prolongent $\sigma:k\hookrightarrow\mathbf{C}_{p}$ permet de voir {\'e}galement que le nombre r{\'e}el $h_{\overline{E_{K}}}(x)$ ne d{\'e}pend pas du choix de l'extension $K/k$ pour laquelle $x\in E\otimes_{k}K$ (lorsque $\overline{E}$ est un fibr{\'e} vectoriel pr{\'e}-ad{\'e}lique).\par L'{\'e}nonc{\'e} suivant remplace l'in{\'e}galit{\'e} de Liouville usuelle.
\begin{lemm}\label{liouville}Soit $\overline{E}$ un fibr{\'e} pr{\'e}-ad{\'e}lique hermitien non nul. Alors \begin{equation*}\forall\,x\in E\setminus\{0\},\quad h_{\overline{E}}(x)\ge-\widehat{\mu}_{\mathrm{max}}(\overline{E}).\end{equation*}
\end{lemm}La d{\'e}monstration est imm{\'e}diate {\`a} partir des d{\'e}finitions car $h_{\overline{E}}(x)=h((k.x,(\Vert\cdot\Vert_{\overline{E},\sigma})_{\sigma}))$.

\subsection{Compl{\'e}ments sur les fibr{\'e}s ad{\'e}liques}\label{paragraphecomplements} Soit $(\mathsf{K},\vert\cdot\vert)$ un corps valu{\'e} ultram{\'e}trique complet. Soit $\mathsf{E}$ un $\mathsf{K}$-espace vectoriel de dimension $\nu\ge 1$ et $\Vert\cdot\Vert$ une ultranorme sur $\mathsf{E}$. \begin{defi}Soit $a\in\,]0,1]$. On dit qu'une famille $\{\mathsf{e}_{1},\ldots,\mathsf{e}_{\nu}\}$ de vecteurs de $\mathsf{E}$ est \emph{$a$-orthogonale} si, pour tous $x_{1},\ldots,x_{\nu}\in\mathsf{K}$, on a $$a\max_{1\le i\le\nu}{\vert x_{i}\vert\Vert\mathsf{e}_{i}\Vert}\le\Vert x_{1}\mathsf{e}_{1}+\cdots+x_{\nu}\mathsf{e}_{\nu}\Vert.$$Si $a=1$ on dit que la base est \emph{orthogonale} et si, de plus, $\Vert\mathsf{e}_{1}\Vert=\cdots=\Vert\mathsf{e}_{\nu}\Vert=1$ alors la base est dite \emph{orthonorm{\'e}e}. 
\end{defi}Lorsque $a<1$, de telles bases existent par le crit{\`e}re suivant d{\^{u}} {\`a} Van der Put (voir~\cite[th{\'e}or{\`e}me~$2.2.16$]{pgs}). Si $\mathsf{e}\in\mathsf{E}$ et $\mathsf{F}\subseteq\mathsf{E}$ est un sous-espace vectoriel, on note $\dist(\mathsf{e},\mathsf{F})$ la borne inf{\'e}rieure des $\Vert\mathsf{e}-\mathsf{f}\Vert$ pour $\mathsf{f}\in\mathsf{F}$.

\begin{lemm}\label{lemmeortho} Soit $t_{2},\ldots,t_{\nu}\in\,]0,1]$ et $(\mathsf{e}_{1},\ldots,\mathsf{e}_{\nu})$ une base de $\mathsf{E}$ tels que, pour tout $i\in\{2,\ldots,\nu\}$, on a  $$\dist(\mathsf{e}_{i},\mathsf{K}.\mathsf{e}_{1}+\cdots+\mathsf{K}.\mathsf{e}_{i-1})\ge t_{i}\Vert\mathsf{e}_{i}\Vert.$$Alors $(\mathsf{e}_{1},\ldots,\mathsf{e}_{\nu})$ est une base $t_{2}\cdots t_{\nu}$-orthogonale de $\mathsf{E}$. En particulier, pour tout $a\in\,]0,1[$, il existe une base $a$-orthogonale de $\mathsf{E}$.\end{lemm} Dans le cas d'une base orthogonale, la r{\'e}ciproque du lemme~\ref{lemmeortho} est vraie~: si $(\mathsf{e}_{1},\ldots,\mathsf{e}_{\nu})$ est une base orthogonale de $(\mathsf{E},\Vert\cdot\Vert)$ alors, pour tout $i\in\{2,\ldots,\nu\}$, on a $\dist(\mathsf{e}_{i},\mathsf{K}.\mathsf{e}_{1}+\cdots+\mathsf{K}.\mathsf{e}_{i-1})=\Vert\mathsf{e}_{i}\Vert.$

 \begin{defi}L'espace norm{\'e} $(\mathsf{E},\Vert\cdot\Vert)$ est dit \emph{sph{\'e}riquement complet} si toute suite de boules embo{\^{\i}}t{\'e}es de $\mathsf{E}$ a une intersection non vide. \end{defi}Un corps local\footnote{Corps complet muni d'une valeur absolue discr{\`e}te (non triviale).}, par exemple $\mathbf{Q}_{p}$, est sph{\'e}riquement complet. On montre que $\mathsf{E}$ (de dimension finie) est sph{\'e}riquement complet si et seulement si $(\mathsf{K},\vert\cdot\vert)$ l'est. Dans ce cas, la distance $\dist(\mathsf{e},\mathsf{F})$ est un minimum et, en particulier, $(\mathsf{E},\Vert\cdot\Vert)$ poss{\`e}de une base orthogonale (cette propri{\'e}t{\'e} est en g{\'e}n{\'e}ral fausse si $\mathsf{K}$ n'est pas sph{\'e}riquement complet, comme par exemple $\mathsf{K}=\mathbf{C}_{p}$, voir~\cite[exemple $2.3.26$]{pgs}). 
 
\begin{lemmeperturbation}Soit $\mathsf{E}$ un espace vectoriel de dimension finie sur un corps valu{\'e} ultram{\'e}trique $(\mathsf{K},\vert\cdot\vert)$ et $\Vert\cdot\Vert$ une ultranorme sur $\mathsf{E}$. Soit $a\in\,]0,1]$ et $(\mathsf{e}_{1},\ldots,\mathsf{e}_{\nu})$ une base $a$-orthogonale de $\mathsf{E}$. Soit $\mathsf{f}_{1},\ldots,\mathsf{f}_{\nu}\in\mathsf{E}$ tels que $\Vert\mathsf{e}_{i}-\mathsf{f}_{i}\Vert<a\Vert\mathsf{e}_{i}\Vert$ pour tout $i\in\{1,\ldots,\nu\}$. Alors $(\mathsf{f}_{1},\ldots,\mathsf{f}_{\nu})$ est une base $a$-orthogonale de $\mathsf{E}$.
\end{lemmeperturbation}
\begin{proof}Par in{\'e}galit{\'e} ultram{\'e}trique on a $\Vert\mathsf{f}_{i}\Vert=\Vert\mathsf{e}_{i}\Vert$. On en d{\'e}duit, pour tous $x_{1},\ldots,x_{\nu}\in\mathsf{K}$, \begin{equation*}\begin{split} a\max_{1\le i\le\nu}{\vert x_{i}\vert\Vert\mathsf{f}_{i}\Vert}&=a\max_{1\le i\le\nu}{\vert x_{i}\vert\Vert\mathsf{e}_{i}\Vert}\le\Vert\sum_{i=1}^{\nu}{x_{i}\mathsf{e}_{i}}\Vert\\ &\le\max{\left(\Vert\sum_{i=1}^{\nu}{x_{i}\mathsf{f}_{i}}\Vert,\Vert\sum_{i=1}^{\nu}{x_{i}(\mathsf{f}_{i}-\mathsf{e}_{i})}\Vert\right)}.\end{split}\end{equation*}Si au moins un des $x_{i}$ n'est pas nul on a $\Vert\sum_{i=1}^{\nu}{x_{i}(\mathsf{f}_{i}-\mathsf{e}_{i})}\Vert<a\max_{1\le i\le\nu}{\vert x_{i}\vert\Vert\mathsf{e}_{i}\Vert}$, ce qui montre alors que $(\mathsf{f}_{1},\ldots,\mathsf{f}_{\nu})$ est une famille $a$-orthogonale. Elle est donc libre et elle forme une base de $\mathsf{E}$ (ceci est un cas particulier du th{\'e}or{\`e}me~$2.3.16$ de~\cite{pgs}).
\end{proof}

L'{\'e}nonc{\'e} suivant est d{\^{u}} {\`a} Ga{\"e}l R{\'e}mond.

\begin{prop}\label{proppure}Les fibr{\'e}s ad{\'e}liques hermitiens purs sur un corps de nombres $k$ sont exactement les fibr{\'e}s vectoriels hermitiens sur $\spec\mathcal{O}_{k}$.
\end{prop}
Un fibr{\'e} vectoriel hermitien sur $\spec\mathcal{O}_{k}$ est toujours pur et l'int{\'e}r{\^{e}}t de cette proposition r{\'e}side dans le fait que la r{\'e}ciproque est vraie. La d{\'e}monstration repose sur le lemme suivant.

\begin{lemm}\label{lemmepurfibre} Soit $\mathsf{K}/\mathsf{k}$ une extension galoisienne de corps valu{\'e}s ultram{\'e}triques. Notons $\vert\cdot\vert$ la valeur absolue sur $\mathsf{K}$ et supposons que $(\mathsf{K},\vert\cdot\vert)$ est sph{\'e}riquement complet. Supposons {\'e}galement que l'extension des corps r{\'e}siduels $\kappa(\mathsf{K})/\kappa(\mathsf{k})$ est s{\'e}parable. Soit $\mathsf{E}$ un $\mathsf{k}$-espace vectoriel de dimension finie et $\Vert\cdot\Vert$ une ultranorme sur $\mathsf{E}\otimes_{\mathsf{k}}\mathsf{K}$, invariante sous l'action du groupe de Galois $\Gal(\mathsf{K}/\mathsf{k})$. Alors toute base orthonorm{\'e}e de $(\mathsf{E},\Vert\cdot\Vert)$ reste une base orthonorm{\'e}e de $(\mathsf{E}\otimes_{\mathsf{k}}\mathsf{K},\Vert\cdot\Vert)$. 
\end{lemm}
L'hypoth{\`e}se d'invariance par Galois de la norme signifie que, pour toute $\mathsf{k}$-base $(\mathsf{e}_{1},\ldots,\mathsf{e}_{\nu})$ de $\mathsf{E}$, pour tout $u\in\Gal(\mathsf{K}/\mathsf{k})$, pour tout $x_{1},\ldots,x_{\nu}\in\mathsf{K}$, on a $$\Vert\sum_{i=1}^{\nu}{\mathsf{e}_{i}\otimes u(x_{i})}\Vert=\Vert\sum_{i=1}^{\nu}{\mathsf{e}_{i}\otimes x_{i}}\Vert.$$L'hypoth{\`e}se de s{\'e}parabilit{\'e} des corps r{\'e}siduels est v{\'e}rifi{\'e}e lorsque $\kappa(\mathsf{k})$ est parfait, par exemple si $\mathsf{k}$ est une extension finie de $\mathbf{Q}_{p}$.

\begin{proof}On raisonne par r{\'e}currence sur la dimension $\nu\ge 1$. Le cas $\nu=1$ {\'e}tant clair, supposons le lemme vrai pour tout espace de dimension $\le\nu-1$. Soit $\mathsf{E}$ de dimension $\nu\ge 2$ et $(\mathsf{e}_{1},\ldots,\mathsf{e}_{\nu})$ une base orthonorm{\'e}e de $(\mathsf{E},\Vert\cdot\Vert)$. Il s'agit de montrer que $$\Vert\sum_{i=1}^{\nu}{\mathsf{e}_{i}\otimes x_{i}}\Vert=\max{\{\vert x_{1}\vert,\ldots,\vert x_{\nu}\vert\}}$$pour tout $x_{1},\ldots,x_{\nu}\in\mathsf{K}$. Soit $\mathsf{F}$ le sous-espace de $\mathsf{E}$ engendr{\'e} par $\mathsf{e}_{1},\ldots,\mathsf{e}_{\nu-1}$. D'apr{\`e}s l'hypoth{\`e}se de r{\'e}currence appliqu{\'e}e {\`a} $(\mathsf{F},\Vert\cdot\Vert)$, la base orthonorm{\'e}e $(\mathsf{e}_{1},\ldots,\mathsf{e}_{\nu-1})$ de $\mathsf{F}$ reste une base $\mathsf{K}$-orthonorm{\'e}e de $(\mathsf{F}\otimes_{\mathsf{k}}\mathsf{K},\Vert\cdot\Vert)$. Comme $\mathsf{K}$ est sph{\'e}riquement complet, il existe $\mathsf{f}\in\mathsf{E}\otimes_{\mathsf{k}}\mathsf{K}\setminus\{0\}$ tel que $\dist(\mathsf{f},\mathsf{F}\otimes_{\mathsf{k}}\mathsf{K})=\Vert\mathsf{f}\Vert$. Quitte {\`a} diviser $\mathsf{f}$ par son coefficient devant $\mathsf{e}_{\nu}$ (n{\'e}cessairement non nul), on peut supposer que $\mathsf{f}$ est de la forme $\mathsf{e}_{1}\otimes\lambda_{1}+\cdots+\mathsf{e}_{\nu-1}\otimes\lambda_{\nu-1}+\mathsf{e}_{\nu}$ avec $\lambda_{1},\ldots,\lambda_{\nu-1}\in\mathsf{K}$. D'apr{\`e}s le lemme~\ref{lemmeortho}, la base $(\mathsf{e}_{1},\ldots,\mathsf{e}_{\nu-1},\mathsf{f})$ de $\mathsf{E}\otimes_{\mathsf{k}}\mathsf{K}$ est orthogonale. Ainsi on a $$1=\Vert\mathsf{e}_{\nu}\Vert=\max{\{\Vert\mathsf{f}\Vert,\vert\lambda_{1}\vert,\ldots,\vert\lambda_{\nu-1}\vert\}}.$$Supposons $\Vert\mathsf{f}\Vert<1$. Pour un {\'e}l{\'e}ment $u$ du groupe de Galois de $\mathsf{K}/\mathsf{k}$, notons $$\mathsf{f}_{u}:=\mathsf{e}_{\nu}+\sum_{i=1}^{\nu-1}{\mathsf{e}_{i}\otimes u(\lambda_{i})}.$$Comme la norme est invariante par Galois on a $\Vert\mathsf{f}_{u}\Vert=\Vert\mathsf{f}\Vert$ donc $$\max_{1\le i\le\nu-1}{\vert u(\lambda_{i})-\lambda_{i}\vert}=\Vert\mathsf{f}_{u}-\mathsf{f}\Vert\le\Vert\mathsf{f}\Vert<1.$$Comme l'extension $\kappa(\mathsf{K})/\kappa(\mathsf{k})$ est s{\'e}parable, elle est galoisienne et l'application canonique $\Gal(\mathsf{K}/\mathsf{k})\to\Gal(\kappa(\mathsf{K})/\kappa(\mathsf{k}))$ qui {\`a} $u$ associe $$\overline{u}: x\bmod\mathsf{M}\mapsto u(x)\bmod\mathsf{M}\quad \text{(o{\`u} $\mathsf{M}:=\{x\in\mathsf{K}\,;\ \vert x\vert<1\}$)}$$est surjective (voir par exemple~\cite[II.9.9,~p.~$172$]{neukirch}). Ainsi, pour tout $i\in\{1,\ldots,\nu-1\}$, l'in{\'e}galit{\'e} $\vert u(\lambda_{i})-\lambda_{i}\vert<1$ donne $\overline{u}(\lambda_{i}\bmod\mathsf{M})=\lambda_{i}\bmod\mathsf{M}$. En faisant varier $u$, l'{\'e}l{\'e}ment $\lambda_{i}\bmod\mathsf{M}$ est invariant par $\Gal(\kappa(\mathsf{K})/\kappa(\mathsf{k}))$ et donc $\lambda_{i}\bmod\mathsf{M}\in\kappa(\mathsf{k})$, ce qui signifie qu'il existe $a_{i}\in\mathsf{k}$ tel que $\vert\lambda_{i}-a_{i}\vert<1$ (en particulier $\vert a_{i}\vert\le 1$). En posant $\mathsf{g}:=\mathsf{e}_{\nu}+\sum_{i=1}^{\nu-1}{\mathsf{e}_{i}\otimes a_{i}}$, on a $\Vert\mathsf{f}-\mathsf{g}\Vert=\max_{i}{\vert\lambda_{i}-a_{i}\vert}<1$ et $\Vert\mathsf{g}\Vert=1$ car $(\mathsf{e}_{1},\ldots,\mathsf{e}_{\nu})$ est une base orthonorm{\'e}e de $\mathsf{E}$. Ceci contredit l'in{\'e}galit{\'e} ultram{\'e}trique $\Vert\mathsf{g}\Vert\le\max{(\Vert\mathsf{g}-\mathsf{f}\Vert,\Vert\mathsf{f}\Vert)}$. Le cas $\Vert\mathsf{f}\Vert<1$ ne peut donc pas  se produire. Ainsi on a $\Vert\mathsf{f}\Vert=1$ puis $$\dist(\mathsf{e}_{\nu},\mathsf{F}\otimes_{\mathsf{k}}\mathsf{K})=\dist(\mathsf{f},\mathsf{F}\otimes_{\mathsf{k}}\mathsf{K})=\Vert\mathsf{f}\Vert=1=\Vert\mathsf{e}_{\nu}\Vert.$$Par cons{\'e}quent $(\mathsf{e}_{1},\ldots,\mathsf{e}_{\nu})$ est une base orthonorm{\'e}e de $(\mathsf{E}\otimes_{\mathsf{k}}\mathsf{K},\Vert\cdot\Vert)$ par la remarque qui suit le lemme~\ref{lemmeortho}. 
\end{proof}

\begin{proof}[D{\'e}monstration de la proposition~\ref{proppure}] Soit $\overline{E}=(E,(\Vert\cdot\Vert_{\overline{E},\sigma})_{\sigma})$ un fibr{\'e} ad{\'e}lique hermitien pur sur $k$. Fixons un plongement ultram{\'e}trique $\sigma:k\hookrightarrow\mathbf{C}_{p}$ et une base orthonorm{\'e}e $(\mathsf{e}_{1},\ldots,\mathsf{e}_{\nu})$ de $(E\otimes_{\sigma}k_{\sigma},\Vert\cdot\Vert_{\overline{E},\sigma})$. Soit $x_{1},\ldots,x_{\nu}\in\mathbf{C}_{p}$ des {\'e}l{\'e}ments alg{\'e}briques sur $\mathbf{Q}_{p}$ et $K/k$ une extension galoisienne finie qui contient ces nombres. Soit $\tau:K\hookrightarrow\mathbf{C}_{p}$ un plongement qui prolonge $\sigma$. On note $w$ la place de $K$ induite par $\tau$. L'extension $K_{\tau}/k_{\sigma}$ est galoisienne et, pour tout $u\in\Gal(K_{\tau}/k_{\sigma})$ il existe un automorphisme continu $\iota:\mathbf{C}_{p}\to\mathbf{C}_{p}$ tel que $\tau\circ u=\iota\circ\tau$ (ces plongements correspondent tous {\`a} $w$). Ainsi la norme $\Vert\cdot\Vert_{\overline{E},\sigma}$ sur $(E\otimes_{\sigma}k_{\sigma})\otimes_{k_{\sigma}}K_{\tau}$ est invariante par le groupe de Galois de $K_{\tau}/k_{\sigma}$ (au sens du lemme~\ref{lemmepurfibre}) car la collection de normes $(\Vert\cdot\Vert_{\overline{E_{K}},\tau})_{\tau\in w}$ est invariante par Galois (au sens de la d{\'e}finition~\ref{defiinvariance}). Le lemme~\ref{lemmepurfibre} s'applique donc {\`a} $(\mathsf{E},\Vert\cdot\Vert)=(E\otimes_{\sigma}k_{\sigma},\Vert\cdot\Vert_{\overline{E},\sigma})$ et $\mathsf{K}=K_{\tau}$ (qui est une extension finie de $\mathbf{Q}_{p}$ donc sph{\'e}riquement complet). On a ainsi $$\Vert\sum_{i=1}^{\nu}{\mathsf{e}_{i}\otimes_{\sigma}x_{i}}\Vert_{\overline{E},\sigma}=\max_{1\le i\le\nu}{\vert x_{i}\vert_{p}}.$$Cette {\'e}galit{\'e} est vraie pour tout $(x_{1},\ldots,x_{\nu})\in\overline{\mathbf{Q}_{p}}^{\nu}$, puis, par continuit{\'e}, sur $\mathbf{C}_{p}^{\nu}$. La matrice de passage de $(\mathsf{e}_{1},\ldots,\mathsf{e}_{\nu})$ {\`a} une $k$-base $(e_{1},\ldots,e_{\nu})$ de $E$ fournit une matrice $A_{\sigma}\in\GL_{\nu}(k_{\sigma})$ telle que $$\Vert\sum_{i=1}^{\nu}{e_{i}\otimes_{\sigma}x_{i}}\Vert=\vert A_{\sigma}\begin{pmatrix}x_{1}\\ \vdots\\ x_{\nu}\end{pmatrix}\vert_{2,\sigma}.$$De plus, l'invariance par Galois de $(\Vert\cdot\Vert_{\overline{E},\sigma})_{\sigma\in v}$ ($v$ est la place de $k$ induite par $\sigma$) permet de choisir ces matrices de sorte que $A_{\iota\circ\sigma}=\iota(A_{\sigma})$ pour tout automorphisme continu $\iota$ de $\mathbf{C}_{p}$. Le fibr{\'e} $\overline{E}$ est donc un fibr{\'e} vectoriel hermitien sur $\spec\mathcal{O}_{k}$.
\end{proof}

Nous allons montrer maintenant que l'on peut approcher un fibr{\'e} ad{\'e}lique hermitien par un fibr{\'e} ad{\'e}lique hermitien \emph{pur}, pourvu que l'on autorise une extension du corps de base.

\begin{lemm}\label{lemmetroisquatre}
Soit $\sigma:k\hookrightarrow\mathbf{C}_{p}$ un plongement ultram{\'e}trique d'un corps de nombres $k$ et $E$ un $k$-espace vectoriel de dimension $\nu\ge 1$. Soit $\Vert\cdot\Vert$ une ultranorme sur $E\otimes_{\sigma}\mathbf{C}_{p}$. Pour tout $\varepsilon>0$, il existe une extension finie $K$ de $k$, il existe une base $f_{1},\ldots,f_{\nu}$ de $E\otimes_{k}K$ et un plongement $\tau:K\hookrightarrow\mathbf{C}_{p}$ qui prolonge $\sigma$ tels que, pour tous $x_{1},\ldots,x_{\nu}\in\mathbf{C}_{p}$, on a $$\max_{1\le i\le\nu}{\vert x_{i}\vert_{p}}\le\Vert\sum_{i=1}^{\nu}{f_{i}\otimes_{\tau}x_{i}}\Vert\le(1+\varepsilon)\max_{1\le i\le\nu}{\vert x_{i}\vert_{p}}.$$
\end{lemm}
\begin{proof}Par le lemme~\ref{lemmeortho}, il existe une base $\mathsf{e}_{1},\ldots,\mathsf{e}_{\nu}$ de $E\otimes_{\sigma}\mathbf{C}_{p}$ telle que, pour tous $x_{1},\ldots,x_{\nu}\in\mathbf{C}_{p}$, on a $$
(1+\varepsilon)^{-1/2}\max_{1\le i\le\nu}{\vert x_{i}\vert_{p}\Vert\mathsf{e}_{i}\Vert}\le\Vert\sum_{i=1}^{\nu}{x_{i}\mathsf{e}_{i}}\Vert\le\max_{1\le i\le\nu}{\vert x_{i}\vert_{p}\Vert\mathsf{e}_{i}\Vert}.$$Comme les nombres alg{\'e}briques sur $\mathbf{Q}$ sont denses dans $\mathbf{C}_{p}$, le lemme de perturbation permet de supposer qu'il existe une extension finie $K/k$ et un plongement $\tau:K\hookrightarrow\mathbf{C}_{p}$ qui prolonge $\sigma$ tels que $\mathsf{e}_{i}\in E\otimes_{\sigma}\tau(K)$. Consid{\'e}rons un entier $M\ge 1$ tel que $p^{1/M}\le(1+\varepsilon)^{1/2}$. Quitte {\`a} prendre une extension finie de $K$ on peut supposer que $p^{1/M}\in K$. Par ailleurs, il existe $\alpha_{i}\in\mathbf{R}$ tel que $(1+\varepsilon)^{-1/2}\Vert\mathsf{e}_{i}\Vert=p^{\alpha_{i}}$. Soit $f_{i}\in E\otimes_{k}K$ tel que $f_{i}\otimes_{\tau}1=p^{[\alpha_{i}M]/M}\mathsf{e}_{i}$. La famille $\{f_{1},\ldots,f_{\nu}\}$ forme une base de $E\otimes_{k}K$ et la norme de $\sum_{i=1}^{\nu}{f_{i}\otimes_{\tau}x_{i}}$ v{\'e}rifie l'encadrement voulu.
\end{proof}On notera que toute extension finie du corps $K$ donn{\'e} par ce lemme convient encore.
\begin{coro}\label{coroimpur}Soit $\overline{E}$ un fibr{\'e} ad{\'e}lique hermitien sur $k$. Alors, pour tout nombre r{\'e}el $\varepsilon>0$, il existe une extension finie $K$ de $k$ et un fibr{\'e} ad{\'e}lique  hermitien pur $\overline{E}_{\varepsilon}$ sur $K$, d'espace sous-jacent $E\otimes_{k}K$, tels que, pour tout plongement $\tau:K\hookrightarrow\mathbf{C}_{p}$ qui prolonge $\sigma:k\hookrightarrow\mathbf{C}_{p}$, pour tout $x\in (E\otimes_{k}K)\otimes_{\tau}\mathbf{C}_{p}$, on a $$\Vert x\Vert_{\overline{E}_{\varepsilon},\tau}\le\Vert x\Vert_{\overline{E_{K}},\tau}\le(1+\varepsilon)\Vert x\Vert_{\overline{E}_{\varepsilon},\tau}.$$De plus, l'in{\'e}galit{\'e} $\Vert\cdot\Vert_{\overline{E}_{\varepsilon},\tau}\le\Vert\cdot\Vert_{\overline{E_{K}},\tau}$ n'est stricte que pour un nombre fini de plongements $\tau$.
\end{coro}
\begin{proof}Fixons une $k$-base $(e_{1},\ldots,e_{\nu})$ de $E$ et un sous-ensemble fini $S$ de $\EuScript{P}$ tel que, pour tout $p\not\in S$, pour tout plongement $\sigma:k\hookrightarrow\mathbf{C}_{p}$, pour tout $(x_{1},\ldots,x_{\nu})\in\mathbf{C}_{p}^{\nu}$, on a$$\Vert e_{1}\otimes_{\sigma}x_{1}+\cdots+e_{\nu}\otimes_{\sigma}x_{\nu}\Vert_{\overline{E},\sigma}=\max{(\vert x_{1}\vert_{p},\ldots,\vert x_{\nu}\vert_{p})}$$(premi{\`e}re condition de la d{\'e}finition de fibr{\'e} pr{\'e}-ad{\'e}lique). On applique le lemme~\ref{lemmetroisquatre} {\`a} chacune des normes $\Vert\cdot\Vert_{\overline{E},\sigma}$ pour $\sigma:k\hookrightarrow\mathbf{C}_{p}$ avec $p\in S\setminus\{\infty\}$. L'on peut supposer que les extensions finies $K_{\sigma}/k$ que l'on obtient sont toutes les m{\^{e}}mes, {\'e}gales {\`a} une extension galoisienne finie $K$ de $k$. Ainsi, pour chacun de ces $\sigma$, il existe un plongement $\tau_{\sigma}:K\hookrightarrow\mathbf{C}_{p}$ prolongeant $\sigma$ et une base $f_{1,1},\ldots,f_{\nu,1}$ de $E\otimes_{k}K$ (qui d{\'e}pend de $\sigma$) tels que, pour tous $x_{1},\ldots,x_{\nu}\in\mathbf{C}_{p}$, on a \begin{equation}\label{ineqessentielle}\max_{1\le i\le\nu}{\vert x_{i}\vert_{p}}\le\Vert\sum_{i=1}^{\nu}{f_{i,1}\otimes_{\tau_{\sigma}}x_{i}}\Vert_{\overline{E_{K}},\tau_{\sigma}}\le(1+\varepsilon)\max_{1\le i\le\nu}{\vert x_{i}\vert_{p}}.\end{equation}Soit $a\in K$ tel que $K=k(a)$ et $\pi_{a}$ le polyn{\^{o}}me minimal de $a$ sur $k$. Comme $K/k$ est galoisienne, l'on peut choisir des {\'e}l{\'e}ments $u_{1}=\mathrm{id},u_{2},\ldots,u_{g}$ du groupe de Galois de $K/k$ tels que les nombres $\tau_{\sigma}\circ u_{\ell}(a)\in\mathbf{C}_{p}$ pour $\ell\in\{1,\ldots,g\}$ soient chacune racine d'un  des $g$ facteurs irr{\'e}ductibles de $\sigma(\pi_{a})$ dans $k_{\sigma}[X]$. Les plongements $ \tau:K\hookrightarrow\mathbf{C}_{p}$ au-dessus de $\sigma$ sont alors de la forme $\iota\circ\tau_{\sigma}\circ u_{\ell}$ avec $\iota$ automorphisme continu de $\mathbf{C}_{p}$ qui laisse fixe les {\'e}l{\'e}ments de $k_{\sigma}$ (on a rang{\'e} les plongements selon les places de $K$ au-dessus de $\sigma$). Consid{\'e}rons des {\'e}l{\'e}ments $\lambda_{i,j}\in K$, $1\le i,j\le\nu$, tels que $f_{i,1}=\sum_{j=1}^{\nu}{e_{j}\otimes_{k}\lambda_{i,j}}$ pour tout $i\in\{1,\ldots,\nu\}$. Posons alors $$f_{i,\ell}:=\sum_{j=1}^{\nu}{e_{j}\otimes_{k}u_{\ell}^{-1}(\lambda_{i,j})}\in E\otimes_{k}K.$$Pour tous $x_{1},\ldots,x_{\nu}\in\mathbf{C}_{p}$ et $\tau=\iota\circ\tau_{\sigma}\circ u_{\ell}$ et par d{\'e}finition des normes sur $\overline{E_{K}}$, on a alors $$\Vert\sum_{i=1}^{\nu}{f_{i,\ell}\otimes_{\tau}x_{i}}\Vert_{\overline{E_{K}},\tau}=\Vert\sum_{j=1}^{\nu}{e_{j}\otimes_{\sigma}\left(\sum_{i=1}^{\nu}{\iota\circ\tau_{\sigma}(\lambda_{i,j})x_{i}}\right)}\Vert_{\overline{E},\sigma}.$$Comme la norme $\Vert\cdot\Vert_{\overline{E},\sigma}$ est invariante par Galois et comme $\iota\circ\sigma=\sigma$, on a \begin{equation}\label{ineqtroises}\Vert\sum_{i=1}^{\nu}{f_{i,\ell}\otimes_{\tau}x_{i}}\Vert_{\overline{E_{K}},\tau}=\Vert\sum_{i=1}^{\nu}{f_{i,1}\otimes_{\tau_{\sigma}}x_{i}}\Vert_{\overline{E_{K}},\tau_{\sigma}}.\end{equation}Sur $(E\otimes_{k}K)\otimes_{\tau}\mathbf{C}_{p}$ (pour $\tau=\iota\circ\tau_{\sigma}\circ u_{\ell}$), d{\'e}finissons la norme $\Vert\cdot\Vert_{\overline{E}_{\varepsilon},\tau}$ par la formule $$\Vert\sum_{i=1}^{\nu}{f_{i,\ell}\otimes_{\tau}x_{i}}\Vert_{\overline{E}_{\varepsilon},\tau}:=\max_{1\le i\le\nu}{\vert x_{i}\vert_{p}}.$$Le couple $\overline{E}_{\varepsilon}=(E\otimes_{k}K,(\Vert\cdot\Vert_{\overline{E}_{\varepsilon},\tau})_{\tau})$ forme un fibr{\'e} hermitien pur sur $K$. De plus, l'encadrement $\Vert\cdot\Vert_{\overline{E}_{\varepsilon},\tau}\le\Vert\cdot\Vert_{\overline{E_{K}},\tau}\le(1+\varepsilon)\Vert\cdot\Vert_{\overline{E}_{\varepsilon},\tau}$ r{\'e}sulte de~\eqref{ineqessentielle} et~\eqref{ineqtroises}.
\end{proof}
Ce corollaire implique que, pour tout $\varepsilon>0$, pour tout fibr{\'e} ad{\'e}lique hermitien $\overline{E}$ sur $k$, il existe une extension finie $K$ de $k$ et un fibr{\'e} hermitien pur $\overline{E'}$  sur $K$, avec $E'=E\otimes_{k}K$, tels que, pour tout $x\in E\otimes_{k}\overline{k}$, on a \begin{equation}\label{inegalitehauteurs}h_{\overline{E'}}(x)\le h_{\overline{E}}(x)\le h_{\overline{E'}}(x)+\varepsilon.\end{equation}De plus, on peut choisir $\overline{E'}$ de sorte que $$h(\overline{E'})\le h(\overline{E})\le h(\overline{E'})+\varepsilon\quad\text{et}\quad\widehat{\mu}_{\mathrm{max}}(\overline{E'})-\varepsilon\le\widehat{\mu}_{\mathrm{max}}(\overline{E})\le\widehat{\mu}_{\mathrm{max}}(\overline{E'})$$car les in{\'e}galit{\'e}s de normes dans le corollaire~\ref{coroimpur} induisent des in{\'e}galit{\'e}s similaires pour les sous-fibr{\'e}s et les d{\'e}terminants d'iceux.

\subsection{Normes d'op{\'e}rateurs}\label{noperateurs}Soit $(\mathsf{K},\vert\cdot\vert)$ un corps valu{\'e} complet. Soit $(\mathsf{E},\Vert\cdot\Vert_{\mathsf{E}})$ (\emph{resp}. $(\mathsf{F},\Vert\cdot\Vert_{\mathsf{F}})$) un $\mathsf{K}$-espace vectoriel norm{\'e} de dimension $\nu\ge 1$ (\emph{resp}. $\mu\ge 1$). Lorsque $\mathsf{K}$ est ultram{\'e}trique, on suppose que les normes sont des ultranormes. Soit $\mathsf{a}:\mathsf{E}\to\mathsf{F}$ une application $\mathsf{K}$-lin{\'e}aire. L'{\'e}quivalence des normes en dimension finie (qui d{\'e}coule par exemple de l'existence de bases $1/2$-orthogonales) permet de poser la
\begin{defi}La \emph{norme d'op{\'e}rateur} de $\mathsf{a}$ est le nombre r{\'e}el $$\Vert\mathsf{a}\Vert:=\sup{\left\{\frac{\Vert\mathsf{a}(x)\Vert_{\mathsf{F}}}{\Vert x\Vert_{\mathsf{E}}}\,;\ x\in\mathsf{E}\setminus\{0\}\right\}}.$$ 
\end{defi}

\begin{prop}\label{propnop} Supposons que $\mathsf{K}$ est ultram{\'e}trique et que $(\mathsf{E},\Vert\cdot\Vert_{\mathsf{E}})$ et $(\mathsf{F},\Vert\cdot\Vert_{\mathsf{F}})$ poss{\`e}dent des bases orthogonales, not{\'e}es respectivement $\mathsf{e}=(\mathsf{e}_{1},\ldots,\mathsf{e}_{\nu})$ et $\mathsf{f}=(\mathsf{f}_{1},\ldots,\mathsf{f}_{\mu})$. Soit $\mathsf{A}=(\mathsf{a}_{i,j})\in\mathrm{M}_{\mu,\nu}(\mathsf{K})$ la matrice qui repr{\'e}sente l'application lin{\'e}aire $\mathsf{a}$ dans les bases $\mathsf{e}$ et $\mathsf{f}$. Alors on a $$\Vert\mathsf{a}\Vert=\max{\left\{\vert\mathsf{a}_{i,j}\vert\Vert\mathsf{f}_{i}\Vert_{\mathsf{F}}/\Vert\mathsf{e}_{j}\Vert_{\mathsf{E}}\,;\ 1\le i\le\mu,\ 1\le j\le\nu\right\}}$$
\end{prop}
\begin{proof}Notons $\alpha$ le maximum du membre de droite. Soit $x=\sum_{j=1}^{\nu}{x_{j}\mathsf{e}_{j}}\in\mathsf{E}$. On a $$\mathsf{a}(x)=\sum_{i=1}^{\mu}{\left(\sum_{j=1}^{\nu}{\mathsf{a}_{i,j}x_{j}}\right)\mathsf{f}_{i}}$$ donc $$\Vert\mathsf{a}(x)\Vert_{\mathsf{F}}\le\max_{i,j}{\vert\mathsf{a}_{i,j}x_{j}\vert\Vert\mathsf{f}_{i}\Vert_{\mathsf{F}}}\le\alpha\max_{1\le j\le\nu}{\vert x_{j}\vert\Vert\mathsf{e}_{j}\Vert_{\mathsf{E}}}=\alpha\Vert x\Vert_{\mathsf{E}}$$puis $\Vert\mathsf{a}\Vert\le\alpha$. R{\'e}ciproquement, consid{\'e}rons des indices $i_{0},j_{0}$ tels que $\alpha=\vert\mathsf{a}_{i_{0},j_{0}}\vert\Vert\mathsf{f}_{i_{0}}\Vert_{\mathsf{F}}/\Vert\mathsf{e}_{j_{0}}\Vert_{\mathsf{E}}$ et prenons $x=\mathsf{e}_{j_{0}}$. On a $\mathsf{a}(x)=\sum_{i=1}^{\mu}{\mathsf{a}_{i,j_{0}}\mathsf{f}_{i}}$ et, par orthogonalit{\'e} de $\mathsf{f}$, on en d{\'e}duit  $$\Vert\mathsf{a}(x)\Vert_{\mathsf{F}}=\max_{1\le i\le\mu}{\vert\mathsf{a}_{i,j_{0}}\vert\Vert\mathsf{f}_{i}\Vert_{\mathsf{F}}}\ge\vert\mathsf{a}_{i_{0},j_{0}}\vert\Vert\mathsf{f}_{i_{0}}\Vert_{\mathsf{F}}=\alpha\Vert x\Vert_{\mathsf{E}}$$puis $\Vert\mathsf{a}\Vert\ge\alpha$.
\end{proof}

\subsection{Puissances sym{\'e}triques des fibr{\'e}s ad{\'e}liques hermitiens}\label{paragraphepuissance} \'Etant donn{\'e} des multiplets $m=(m_{1},\ldots,m_{\nu})$ et $n=(n_{1},\ldots,n_{\nu})$, on note $\vert m\vert$ la \emph{longueur} $\sum_{i=1}^{\nu}{m_{i}}$ de $m$ et on d{\'e}signe par $m!$ (\emph{resp}. $m^{n}$) le produit $m_{1}!\cdots m_{\nu}!$ (\emph{resp}. $m^{n}:=m_{1}^{n_{1}}\cdots m_{\nu}^{n_{\nu}}$). Soit $\ell\ge 1$ un entier et $\overline{E}$ un fibr{\'e} pr{\'e}-ad{\'e}lique hermitien sur $k$. La puissance sym{\'e}trique $S^{\ell}(E)$ est un quotient de $E^{\otimes\ell}$ (l'alg{\`e}bre sym{\'e}trique $\mathbf{S}(E)$ est le quotient de l'alg{\`e}bre tensorielle $\mathbf{T}(E)$ par l'id{\'e}al engendr{\'e} par les {\'e}l{\'e}ments $x\otimes y-y\otimes x$). Cette structure quotient conf{\`e}re {\`a} $S^{\ell}(E)$ une structure de fibr{\'e} pr{\'e}-ad{\'e}lique hermitien $\overline{S^{\ell}(E)}$ induite par $\overline{E}^{\otimes\ell}$. 
\begin{prop}\label{propestipentesym} Pour tout fibr{\'e} ad{\'e}lique hermitien $\overline{E}$ sur $k$ de dimension $\nu\ge 1$, pour tout entier $\ell\ge 1$, on a \begin{equation*}\widehat{\mu}_{\mathrm{max}}(\overline{S^{\ell}(E)})\le\ell(\widehat{\mu}_{\mathrm{max}}(\overline{E})+2\log\nu).\end{equation*}
\end{prop}
\begin{proof}Si $\overline{E}$ est un fibr{\'e} vectoriel hermitien sur $\spec\mathcal{O}_{k}$, ce r{\'e}sultat est compris dans la proposition~$8.4$ de~\cite{ga}. Le passage aux fibr{\'e}s hermitiens quelconques s'effectue au moyen du corollaire~\ref{coroimpur}, en utilisant l'invariance par extension des scalaires de la pente maximale de $\overline{E}$ (proposition~\ref{propmumax}).
\end{proof}

Dans la suite, la puissance sym{\'e}trique sera utilis{\'e}e {\`a} la fois pour d{\'e}crire l'ensemble des polyn{\^{o}}mes auxiliaires et pour d{\'e}crire l'espace des d{\'e}rivations. Si $E$ est un $k$-espace vectoriel, on rappelle que $E^{\mathsf{v}}$ d{\'e}signe l'espace dual $\Hom_{k}(E,k)$.
\begin{defi}
Un \emph{polyn{\^{o}}me homog{\`e}ne} $s$ sur $E$ de degr{\'e} $\ell$ est un {\'e}l{\'e}ment de $S^{\ell}(E^{\mathsf{v}})$ et l'\emph{application polynomiale} associ{\'e}e est l'application de $E$ dans $k$ qui envoie $x\in E$ sur $s(x)$. Plus g{\'e}n{\'e}ralement, si $n\in\mathbf{N}$ et si $E_{0},\ldots,E_{n}$ sont des $k$-espaces vectoriels, un \emph{polyn{\^{o}}me multihomog{\`e}ne} $s$ de degr{\'e} $(\ell_{0},\ldots,\ell_{n})\in\mathbf{N}^{n+1}$ est un {\'e}l{\'e}ment de $\mathfrak{E}:=S^{\ell_{0}}(E_{0}^{\mathsf{v}})\otimes\cdots\otimes S^{\ell_{n}}(E_{n}^{\mathsf{v}})$ et l'application polynomiale associ{\'e}e est l'application de $\prod_{j=0}^{n}{E_{j}}$ dans $k$ qui envoie $(x_{0},\ldots,x_{n})$ sur $s(x_{0},\ldots,x_{n})$. \end{defi}Lorsque, pour tout $j\in\{0,\ldots,n\}$, l'espace $E_{j}$ poss{\`e}de une structure de fibr{\'e} ad{\'e}lique hermitien $\overline{E_{j}}$ (il en est alors de m{\^{e}}me pour $E_{j}^{\mathsf{v}}$ avec les m{\'e}triques duales), la norme de $s$ en une place $v$ de $k$ est calcul{\'e}e dans $\bigotimes_{j=0}^{n}{\overline{S^{\ell_{j}}(E_{j}^{\mathsf{v}})}}$. Dans la pratique, si $\overline{E}$ est un fibr{\'e} ad{\'e}lique hermitien pur sur $k$, si l'on fixe une base orthonorm{\'e}e $\mathsf{e}_{j}=(\mathsf{e}_{j,1},\ldots,\mathsf{e}_{j,\nu_{j}})$ de $E_{j}\otimes_{\sigma}\mathbf{C}_{p}$ (ici $\nu_{j}=\dim E_{j}$), de base duale $\mathsf{e}_{j}^{\mathsf{v}}=(\mathsf{e}_{j,1}^{\mathsf{v}},\ldots,\mathsf{e}_{j,\nu_{j}}^{\mathsf{v}})$, et si l'on {\'e}crit l'{\'e}l{\'e}ment $s$ comme une somme \begin{equation}\label{ecrituredes}s:=\sum_{i}{p_{i}\prod_{j=0}^{n}{(\mathsf{e}_{j}^{\mathsf{v}})^{i_{j}}}}\end{equation}avec $i=(i_{0},\ldots,i_{n})\in\mathbf{N}^{\nu_{0}}\times\cdots\times\mathbf{N}^{\nu_{n}}$, $i_{j}$ de longueur $\ell_{j}$, et $p_{i}\in\mathbf{C}_{p}$, alors la norme $\Vert s\Vert_{\overline{\mathfrak{E}},\sigma}$ de $s$ est \begin{equation}\label{puissancedeux}\Vert s\Vert_{\overline{\mathfrak{E}},\sigma}=\begin{cases}\left(\sum_{i}{\vert\sigma(p_{i})\vert^{2}\frac{i!}{\ell_{0}!\cdots\ell_{n}!}}\right)^{1/2} & \text{si $p=\infty$}\\ \max_{i}{\vert p_{i}\vert_{p}} & \text{si $p\ne\infty$}.\end{cases}\end{equation}Remarquons {\'e}galement que gr{\^{a}}ce {\`a} cette formule~\eqref{puissancedeux} et {\`a} l'in{\'e}galit{\'e} de Cauchy-Schwarz, la somme $\sum_{i}{\vert\sigma(p_{i})\vert}$ est plus petite que $\Vert s\Vert_{\overline{\mathfrak{E}},\sigma}\times\prod_{j=0}^{n}{\nu_{j}^{\ell_{j}/2}}$ (si $\sigma$ est complexe).

\section{Lemmes de Siegel approch{\'e}s}\label{lemmedepetitesvaleurs}
\subsection{} La plupart des d{\'e}monstrations de transcendance requiert l'utilisation d'une fonction auxiliaire qui doit satisfaire {\`a} un nombre fini de conditions lin{\'e}aires. En g{\'e}n{\'e}ral il s'agit de conditions d'annulations en des points particuliers avec des ordres de multiplicit{\'e}s dans certaines directions (on parle parfois de \og points {\'e}paissis\fg{}). Aux pr{\'e}mices de la th{\'e}orie (travaux d'Hermite et Lindemann par exemple), on exhibait hardiment une fonction auxiliaire \emph{explicite}. Toutefois sont rapidement apparues les difficult{\'e}s et les limitations inh{\'e}rentes {\`a} cette approche presque impudique\footnote{Cependant, ce proc{\'e}d{\'e} reste encore tr{\`e}s actuel au travers par exemple des approximants de Pad{\'e} et des fonctions hyperg{\'e}om{\'e}triques car, comme l'avait not{\'e} Chudnovsky, il conduit souvent {\`a} de meilleurs r{\'e}sultats.}. En filigrane dans les articles de Thue~\cite{thuep,thue}, l'on doit {\`a} Siegel d'avoir conceptualis{\'e} en $1929$ l'id{\'e}e qu'il suffisait de conna{\^\i}tre une estimation de la \og taille\fg{} de la fonction auxiliaire $F$. Demander l'annulation de $F$ en un nombre fini de points {\'e}paissis {\'e}quivaut {\`a} r{\'e}clamer que les coefficients de $F$ satisfassent {\`a} un syst{\`e}me lin{\'e}aire\begin{equation}\label{syssiegel} \forall\,i\in\{1,\ldots,\mu\},\quad\sum_{j=1}^{\nu}{a_{i,j}x_{j}}=0\end{equation}d'inconnues $x_{1},\ldots,x_{\nu}$. Lors de l'{\'e}tude de la transcendance des valeurs des fonctions de Bessel~\cite{siegel}, Siegel formula l'{\'e}nonc{\'e} pr{\'e}cis suivant.
\begin{lemmesiegel}Supposons que $\mu<\nu$ et que, pour tous $i,j$, on a $a_{i,j}\in\mathbf{Z}$. Soit $A=\max_{i,j}{\vert a_{i,j}\vert}$. Alors il existe une solution $(x_{1},\ldots,x_{\nu})\in\mathbf{Z}^{\nu}\setminus\{0\}$ au syst{\`e}me~\eqref{syssiegel} telle que \begin{equation*}\max{\{\vert x_{1}\vert,\ldots,\vert x_{\nu}\vert\}}\le1+(\nu A)^{\frac{\mu}{\nu-\mu}}.\end{equation*}\end{lemmesiegel}Ce r{\'e}sultat d{\'e}coule simplement du \emph{principe des tiroirs de Dirichlet}. Comme l'a remarqu{\'e} Mignotte, il est possible de l'{\'e}tendre {\`a} un syst{\`e}me {\`a} coefficients alg{\'e}briques (voir lemme~$1.3.1$ de~\cite{miw402}). Cependant l'on s'est aper{\c c}u qu'un tel lemme n'{\'e}tait rien d'autre qu'une variante du \emph{premier th{\'e}or{\`e}me de Minkowski} sur les corps convexes. Ce point de vue s'est r{\'e}v{\'e}l{\'e} f{\'e}cond en d{\'e}bouchant sur une version ad{\'e}lique du lemme de Siegel, d{\'e}montr{\'e}e par Bombieri \& Vaaler~\cite{BombieriVaaler} (voir~\cite[\S~$3.2$]{evitement} pour la constante $(1/2)\log\nu$).
\begin{lemmeBV}\label{lemmebombierivaaler}
Soit $k$ un corps de nombres de discriminant absolu $D_{k}$. Notons $\mathrm{rd}_{k}:=\vert D_{k}\vert^{1/[k:\mathbf{Q}]}$ son discriminant racine. Soit $\overline{E}$ un fibr{\'e} ad{\'e}lique hermitien  pur sur $k$, de dimension $\nu\ge 1$. Alors il existe $x\in E\setminus\{0\}$ tel que \begin{equation*}h_{\overline{E}}(x)\le-\widehat{\mu}(\overline{E})+\frac{1}{2}(\log\nu+\log\mathrm{rd}_{k}).\end{equation*}
\end{lemmeBV}
Plusieurs d{\'e}veloppements de ce lemme existent.  On peut demander que l'{\'e}l{\'e}ment $x\in E\setminus\{0\}$ {\'e}vite un nombre fini de sous-espaces vectoriels de $E$. Ceci constitue le fil directeur des articles~\cite{aneuf,evitement}. Mais la variante la plus importante est celle qui permet de s'affranchir de la d{\'e}pendance en le discriminant du corps de nombres. On recherche une solution $x$ non nulle non pas dans $E$ mais dans $E\otimes\overline{\mathbf{Q}}$. Un tel lemme de Siegel est dit \emph{absolu}. Roy \& Thunder ont obtenu un  {\'e}nonc{\'e} de ce type dans~\cite{roythunder}. Nous pr{\'e}sentons ici un raffinement de leur r{\'e}sultat, signal{\'e} par David \& Philippon~\cite{david-pph}, qui se d{\'e}duit d'une in{\'e}galit{\'e} de Zhang relative aux minima successifs d'une vari{\'e}t{\'e} arithm{\'e}tique~\cite{Zhangdeux}. 
\begin{lemmeSA}\label{lemmesiegelabsolu}
Soit $\overline{E}$ un fibr{\'e} ad{\'e}lique hermitien sur $k$, de dimension $\nu\ge 1$. Alors il existe $x\in(E\otimes\overline{\mathbf{Q}})\setminus\{0\}$ tel que \begin{equation*}h_{\overline{E}}(x)\le-\widehat{\mu}(\overline{E})+\frac{1}{2}\log\nu.\end{equation*}
\end{lemmeSA}L'{\'e}nonc{\'e} de Roy \& Thunder donne $(\nu-1)/4+\varepsilon$ au lieu de $\frac{1}{2}\log\nu$ lorsque $\overline{E}$ est pur. Avec la m{\^{e}}me hypoth{\`e}se, Zhang montre que, pour tout $\varepsilon>0$, il existe $x\in E\otimes_{k}\overline{\mathbf{Q}}$ de hauteur plus petite que $-\widehat{\mu}(\overline{E})+(H_{\nu}-1)/2+\varepsilon$, o{\`u} $H_{\nu}=1+1/2+\cdots+1/\nu$ est le nombre harmonique (la constante $(1/2)\log\nu$ {\'e}tant un majorant strict de $(H_{\nu}-1)/2$ lorsque $\nu\ge 2$). \textit{A priori} valide pour les fibr{\'e}s ad{\'e}liques hermitiens purs sur $k$, cet {\'e}nonc{\'e} reste vrai pour les fibr{\'e}s ad{\'e}liques hermitiens gr{\^{a}}ce au corollaire~\ref{coroimpur} (et aux cons{\'e}quences qui suivent).
\par Si, autrefois, dans les d{\'e}monstrations de mesures d'ind{\'e}pendance lin{\'e}aire de logarithmes, le discriminant du corps de nombres {\'e}tait absorb{\'e} par des termes plus gros, aujourd'hui les mesures sont devenues assez fines pour que ce discriminant devienne un facteur limitant. C'est la pr{\'e}sence de ce discriminant qui explique pourquoi l'{\'e}nonc{\'e} principal de~\cite{gaudron1} s'exprime au moyen d'un maximum sur deux quantit{\'e}s. Comme il supprime cette imperfection, le lemme de Siegel absolu a pris une grande importance ces derni{\`e}res ann{\'e}es, ainsi que le prouve son utilisation dans plusieurs articles r{\'e}cents de la th{\'e}orie des formes lin{\'e}aires logarithmes~\cite{davidhiratacrelles,AblyGaudron,casrationnel}. Le fait que l'on ne ma{\^\i}trise pas le (degr{\'e} du) corps de nombres dans lequel vit la solution $x$ est sans cons{\'e}quence car, en d{\'e}finitive, l'on est amen{\'e} {\`a} effectuer une somme sur les diff{\'e}rents plongements de ce corps, somme de laquelle {\'e}merge directement $h_{\overline{E}}(x)$. N{\'e}anmoins quelques soucis techniques peuvent surgir.\par Tout d'abord, la hauteur de $\overline{E}$ peut s'av{\'e}rer difficile {\`a} {\'e}valuer avec pr{\'e}cision. Si, comme dans le lemme de Siegel original, l'espace vectoriel $E$ est d{\'e}fini par le syst{\`e}me lin{\'e}aire~\eqref{syssiegel}, la hauteur de $\overline{E}$ se calcule au moyen des mineurs maximaux de la matrice $\mathsf{A}:=(a_{i,j})_{i,j}$. En g{\'e}n{\'e}ral, on estime la taille de ces mineurs avec l'in{\'e}galit{\'e} d'Hadamard (le d{\'e}terminant d'une famille de vecteurs est plus petit que le produit des normes hermitiennes de ces vecteurs), qui, \textit{in fine}, fait ressortir la hauteur des $a_{i,j}$. Il arrive que les $a_{i,j}$ soient petits aux places archim{\'e}diennes (ce qui est tr{\`e}s bien) mais avec un d{\'e}nominateur trop grand aux places ultram{\'e}triques. L'astuce consiste alors {\`a} trouver un syst{\`e}me \begin{equation}\label{syssiegeldeux}\forall\,i\in\{1,\ldots,\mu\},\quad\sum_{j=1}^{\nu}{b_{i,j}x_{j}}=0\ ,\end{equation}{\'e}quivalent {\`a}~\eqref{syssiegel} et d{\'e}finissant ainsi le m{\^e}me espace vectoriel $E$, mais, o{\`u}, cette fois-ci, les nombres alg{\'e}briques $b_{i,j}$ sont petits aux places ultram{\'e}triques et de tailles quelconques aux autres places. Dans ce cas, la partie archim{\'e}dienne de la hauteur de $\overline{E}$ est {\'e}valu{\'e}e avec~\eqref{syssiegel} et la partie ultram{\'e}trique avec~\eqref{syssiegeldeux}. Cette technique fonctionne assez bien m{\^{e}}me si elle est parfois d{\'e}licate {\`a} mettre en {\oe}uvre (voir p.~ex. la d{\'e}monstration de la proposition~$4.15$ de~\cite{casrationnel}). \par Une autre difficult{\'e} est que les coefficients $a_{i,j}$ qui se pr{\'e}sentent naturellement peuvent ne pas {\^{e}}tre alg{\'e}briques, m{\^e}me apr{\`e}s renormalisation. Le syst{\`e}me~\eqref{syssiegel} ne poss{\`e}de alors en g{\'e}n{\'e}ral aucune solution alg{\'e}brique hormis $(0,\ldots,0)$. C'est pourquoi il est plus raisonnable de demander au lemme de fournir une solution alg{\'e}brique $(x_{1},\ldots,x_{\nu})$ non nulle au syst{\`e}me d'\emph{in{\'e}quations}\begin{equation}\label{sysbis}\forall\,i\in\{1,\ldots,\mu\},\quad\left\vert\sum_{j=1}^{\nu}{a_{i,j}x_{j}}\right\vert\le\varepsilon\end{equation}(ici $\varepsilon$ est un nombre r{\'e}el strictement positif). Un {\'e}nonc{\'e} qui garantit l'existence d'une solution alg{\'e}brique non nulle {\`a} ce syst{\`e}me d'in{\'e}quations est appel{\'e} \emph{lemme de Siegel approch{\'e}}. En voici un exemple\footnote{Cet exemple est d{\^{u}} en partie {\`a} Mignotte~\cite{mignotte1975}, qui s'est int{\'e}ress{\'e} {\`a} la question dans les ann{\'e}es 70. On pourra aussi consulter le \S~$1.3$ du chapitre $1$ de~\cite{miw402}.}, extrait de l'article de Philippon \& Waldschmidt~\cite{pph-miw}, qui a servi {\`a} la construction de la fonction auxiliaire de plusieurs articles marquants de la th{\'e}orie des formes lin{\'e}aires de logarithmes~\cite{pph-miw,pphmiw2,hirata2,hirata3} (aussi utilis{\'e} dans~\cite{gaudron1})~:
\begin{lemmepv}   
Soit $(a_{i,j})$, $1\le i\le\mu$, $1\le j\le\nu$, une matrice de nombres complexes de rang $\rho$ et $A$ un nombre r{\'e}el tel que $\max_{1\le i\le\mu}{\sum_{j=1}^{\nu}{\vert a_{i,j}\vert}}\le A$. Soit $H\in\mathbf{N}\setminus\{0\}$ et $\varepsilon\in\,]0,+\infty[$ tels que $$\left(\frac{2\mu HA}{\varepsilon}+1\right)^{2\rho}<(H+1)^{\nu}.$$Alors il existe $(x_{1},\ldots,x_{\nu})\in\mathbf{Z}^{\nu}\setminus\{0\}$ tel que \begin{equation*}\max_{1\le j\le\nu}{\vert x_{j}\vert}\le H\quad\text{et}\quad\max_{1\le i\le\mu}{\left\vert\sum_{j=1}^{\nu}{a_{i,j}x_{j}}\right\vert}\le\varepsilon.\end{equation*} 
\end{lemmepv}
Tout comme pour le lemme de Siegel original, la d{\'e}monstration repose sur le principe des tiroirs. S'il est facile de g{\'e}n{\'e}raliser {\`a} un corps de nombres au moyen d'une $\mathbf{Q}$-base $(\xi_{1},\ldots,\xi_{D})$ de ce corps, ceci fait intervenir la hauteur de cette base. Le corps consid{\'e}r{\'e} est souvent le corps de nombres dans lequel vivent tous les nombres alg{\'e}briques de la d{\'e}monstration. En particulier, pour les formes lin{\'e}aires de logarithmes, la hauteur de $(\xi_{1},\ldots,\xi_{D})$ est li{\'e}e aux param{\`e}tres $\log a$ et $\log b$. Pour faire dispara{\^\i}tre la d{\'e}pendance en le corps de nombres ambiant, l'on peut imaginer {\'e}crire un lemme de Siegel approch{\'e} et \emph{absolu}. Il s'av{\`e}re qu'un lemme de Siegel (classique ou absolu) donne automatiquement un lemme de Siegel approch{\'e} par \emph{d{\'e}formation des normes}\footnote{\textit{A posteriori}, nous nous sommes rendu compte que cette id{\'e}e {\'e}tait d{\'e}j{\`a} en substance dans l'article~\cite{mignotte1975} de Mignotte et qu'elle serait due {\`a} Mahler (voir~\cite{mahlerbook}).}, comme nous allons l'expliquer maintenant. \par Soit $\overline{E}=(E,(\Vert\cdot\Vert_{\overline{E},\sigma})_{\sigma})$ et $\overline{F}=(F,(\Vert\cdot\Vert_{\overline{F},\sigma})_{\sigma})$ des fibr{\'e}s ad{\'e}liques hermitiens sur un corps de nombres $k$. Soit $S$ un ensemble fini (non vide) de plongements de $k$. Pour tout $\sigma\in S$, soit $\mathsf{a}_{\sigma}:E\otimes_{\sigma}\mathbf{C}_{p}\to F\otimes_{\sigma}\mathbf{C}_{p}$ une application $\mathbf{C}_{p}$-lin{\'e}aire et posons $\mathsf{a}:=(\mathsf{a}_{\sigma})_{\sigma\in S}$. On consid{\`e}re sur $E\otimes_{\sigma}\mathbf{C}_{p}$ la norme tordue par $\mathsf{a}_{\sigma}$~:\begin{equation*}\forall\,x\in E\otimes_{\sigma}\mathbf{C}_{p},\quad \Vert x\Vert_{\overline{E}_{\mathsf{a}},\sigma}:=\moy_{\sigma}(\Vert x\Vert_{\overline{E},\sigma},\Vert\mathsf{a}_{\sigma}(x)\Vert_{\overline{F},\sigma}).\end{equation*}Si $\sigma\not\in S$, on pose $\Vert\cdot\Vert_{\overline{E}_{\mathsf{a}},\sigma}:=\Vert\cdot\Vert_{\overline{E},\sigma}$. Le couple $\overline{E}_{\mathsf{a}}=(E,(\Vert\cdot\Vert_{\overline{E}_{\mathsf{a}},\sigma}))$ forme un fibr{\'e} pr{\'e}-ad{\'e}lique hermitien sur $k$. Une condition suffisante pour qu'il soit ad{\'e}lique est que, pour tout $\sigma:k\hookrightarrow\mathbf{C}_{p}$ dans $S$, pour tout automorphisme continu $\iota:\mathbf{C}_{p}\to\mathbf{C}_{p}$, on a $\iota\circ\sigma\in S$ et les deux {\'e}l{\'e}ments de $F\otimes_{\iota\circ\sigma}\mathbf{C}_{p}$ que l'on obtient par les deux chemins du diagramme
\begin{equation}\label{diagramme}\xymatrix{E\otimes_{\sigma}\mathbf{C}_{p}\ar[r]^{\mathsf{a}_{\sigma}}\ar[d]_{\mathrm{id}\otimes\iota} & F\otimes_{\sigma}\mathbf{C}_{p}\ar[d]^{\mathrm{id}\otimes\iota}\\ E\otimes_{\iota\circ\sigma}\mathbf{C}_{p}\ar[r]_{\mathsf{a}_{\iota\circ\sigma}} & F\otimes_{\iota\circ\sigma}\mathbf{C}_{p}}
\end{equation} sont de m{\^{e}}me norme. De plus, si $x\in E$, on a $h_{\overline{E}}(x)\le h_{\overline{E}_{\mathsf{a}}}(x)$ et $\Vert\mathsf{a}_{\sigma}(x)\Vert_{\overline{F},\sigma}\le\Vert x\Vert_{\overline{E}_{\mathsf{a}},\sigma}$ pour $\sigma\in S$. D{\`e}s lors, savoir majorer finement $h_{\overline{E}_{\mathsf{a}}}(x)$, pour un vecteur $x$ particulier, am{\`e}ne {\`a} {\'e}tablir un lemme de Siegel approch{\'e}. On note $\Vert\mathsf{a}_{\sigma}\Vert$ la norme d'op{\'e}rateur de $\mathsf{a}_{\sigma}$ (voir paragraphe~\ref{noperateurs})~: $$\Vert\mathsf{a}_{\sigma}\Vert:=\sup{\{\Vert\mathsf{a}_{\sigma}(x)\Vert_{\overline{F},\sigma}/\Vert x\Vert_{\overline{E},\sigma}\,;\ x\in E\otimes_{\sigma}\mathbf{C}_{p}\setminus\{0\}\}}.$$

\begin{lemmepva}
Soit $k$ un corps de nombres de degr{\'e} $D$ et $S$ un ensemble fini de plongements de $k$. Soit $\overline{E}$ et $\overline{F}$ deux fibr{\'e}s ad{\'e}liques hermitiens sur $k$. Pour tout $\sigma\in S$, soit $\mathsf{a}_{\sigma}:E\otimes_{\sigma}\mathbf{C}_{p}\to F\otimes_{\sigma}\mathbf{C}_{p}$ une application $\mathbf{C}_{p}$-lin{\'e}aire, de rang $\rho_{\sigma}$, et de norme d'op{\'e}rateur $\Vert\mathsf{a}_{\sigma}\Vert$. On suppose que $\overline{E}_{\mathsf{a}}$ est un fibr{\'e} ad{\'e}lique hermitien. Alors il existe $x\in (E\otimes\overline{\mathbf{Q}})\setminus\{0\}$ tel que \begin{equation*}h_{\overline{E}_{\mathsf{a}}}(x)\le\frac{1}{\nu D}\sum_{\sigma\in S}{\rho_{\sigma}\log\moy_{\sigma}{(1,\Vert\mathsf{a}_{\sigma}\Vert)}}+\frac{1}{2}\log\nu-\widehat{\mu}(\overline{E})\end{equation*}(ici $\moy_{\sigma}=\moy_{p}$ o{\`u} $p\in\EuScript{P}$ est celui correspondant {\`a} $\sigma:k\hookrightarrow\mathbf{C}_{p}$).
\end{lemmepva}
Une caract{\'e}ristique importante du majorant est la pr{\'e}sence des quotients $\rho_{\sigma}/\nu$, qui seront petits dans le contexte des formes lin{\'e}aires de logarithmes (voir proposition~\ref{proprang}). Dans un souci d'efficacit{\'e}, nous n'avons pas s{\'e}par{\'e} les contributions des normes de $x$ en les diff{\'e}rents plongements de $\overline{\mathbf{Q}}$, en exprimant simplement le r{\'e}sultat en termes de hauteur (globale) de $x$. Par ailleurs une mani{\`e}re de proc{\'e}der pour faire en sorte que $\overline{E}_{\mathsf{a}}$ soit un fibr{\'e} ad{\'e}lique hermitien est de consid{\'e}rer un ensemble $S_{0}$ de plongements $\sigma :k\hookrightarrow\mathbf{C}_{p}$ qui n'induisent pas la m{\^{e}}me valuation. Pour un plongement conjugu{\'e} {\`a} $\sigma$, c'est-{\`a}-dire pour $\sigma'=\iota\circ\sigma$ avec $\iota$ automorphisme continu de $\mathbf{C}_{p}$, on choisit $\mathsf{a}_{\sigma'}$ de sorte que le diagramme~\eqref{diagramme}  v{\'e}rifie la condition d'isom{\'e}trie et l'on prend $S=\{\sigma'\,;\ \sigma\in S_{0}\}$ (le cardinal de $S$ est $\sum_{\sigma\in S_{0}}{[k_{\sigma}:\mathbf{Q}_{p}]}$).

\subsection{D{\'e}monstration du lemme de Siegel approch{\'e} absolu}En vertu du lemme de Siegel absolu de Zhang appliqu{\'e} {\`a} $\overline{E}_{\mathsf{a}}$, il suffit de montrer que  \begin{equation*}-\widehat{\mu}(\overline{E}_{\mathsf{a}})\le\frac{1}{\nu D}\sum_{\sigma\in S}{\rho_{\sigma}\log\moy_{\sigma}{(1,\Vert\mathsf{a}_{\sigma}\Vert)}}-\widehat{\mu}(\overline{E}).\end{equation*}Par d{\'e}finition du degr{\'e} ad{\'e}lique, pour toute $k$-base $e_{1},\ldots,e_{\nu}$ de $E$, on a  \begin{equation*}\widehat{\mu}(\overline{E}_{\mathsf{a}})-\widehat{\mu}(\overline{E})=-\frac{1}{\nu D}\sum_{p\in\EuScript{P}}{\sum_{\sigma:k\hookrightarrow\mathbf{C}_{p}}{\log\frac{\Vert e_{1}\wedge\cdots\wedge e_{\nu}\Vert_{\overline{\det E_{\mathsf{a}}},\sigma}}{\Vert e_{1}\wedge\cdots\wedge e_{\nu}\Vert_{\overline{\det E},\sigma}}}}.\end{equation*}Dans cette somme, seuls les $\sigma\in S$ interviennent, les autres donnant des termes nuls. De plus, {\`a} $\sigma$ fix{\'e}, chacun des quotients $\Vert e_{1}\wedge\cdots\wedge e_{\nu}\Vert_{\overline{\det E_{\mathsf{a}}},\sigma}/\Vert e_{1}\wedge\cdots\wedge e_{\nu}\Vert_{\overline{\det E},\sigma}$ est ind{\'e}pendant du choix de la base $e_{1},\ldots,e_{\nu}$ de $E\otimes_{\sigma}\mathbf{C}_{p}$. Pour calculer ce quotient, supposons dans un premier temps que $\sigma$ est ultram{\'e}trique et consid{\'e}rons un nombre r{\'e}el $t\in\,]0,1]$ et une base $t$-orthogonale $\mathsf{e}_{1},\ldots,\mathsf{e}_{\nu}$ de $E\otimes_{\sigma}\mathbf{C}_{p}$ telle que $\mathsf{e}_{\rho_{\sigma}+1},\ldots,\mathsf{e}_{\nu}$ soit une base de $\ker\mathsf{a}_{\sigma}$~\cite[corollaire~$2.3.21$]{pgs}. Alors, d'une part, comme la norme sur le d{\'e}terminant est une norme quotient, on a \begin{equation*}\Vert\mathsf{e}_{1}\wedge\cdots\wedge\mathsf{e}_{\nu}\Vert_{\overline{\det E_{\mathsf{a}}},\sigma}\le\prod_{i=1}^{\nu}{\Vert\mathsf{e}_{i}\Vert_{\overline{E}_{\mathsf{a}},\sigma}}\le\left(\prod_{i=1}^{\nu}{\Vert\mathsf{e}_{i}\Vert_{\overline{E},\sigma}}\right)\moy_{\sigma}{(1,\Vert\mathsf{a}_{\sigma}\Vert)}^{\rho_{\sigma}}.\end{equation*}D'autre part, la base $\{\mathsf{e}_{i_{1}}\otimes\cdots\otimes\mathsf{e}_{i_{\nu}},\ 1\le i_{1},\ldots,i_{\nu}\le\nu\}$ de $(E\otimes_{\sigma}\mathbf{C}_{p})^{\otimes\nu}$ est $t^{\nu}$-orthogonale~\cite[corollaire~$10.2.10$, (vi)]{pgs}. Si un vecteur $\mathsf{x}\in(E\otimes_{\sigma}\mathbf{C}_{p})^{\otimes\nu}$ a pour coordonn{\'e}es $(\mathsf{x}_{i_{1},\ldots,i_{\nu}})$ dans cette base alors son image dans $\det E\otimes_{\sigma}\mathbf{C}_{p}$ est le vecteur $\left(\sum_{s\in\mathfrak{S}_{\nu}}{\varepsilon(s)\mathsf{x}_{s(1),\ldots,s(\nu)}}\right)\mathsf{e}_{1}\wedge\cdots\wedge\mathsf{e}_{\nu}$. On en d{\'e}duit que, pour tout $\mathsf{x}\in(E\otimes_{\sigma}\mathbf{C}_{p})^{\otimes\nu}$, on a $$\Vert\mathsf{e}_{1}\otimes\cdots\otimes\mathsf{e}_{\nu}+\mathsf{x}\Vert_{\overline{\det E},\sigma}\ge t^{\nu}\left(\prod_{i=1}^{\nu}{\Vert\mathsf{e}_{i}\Vert_{\overline{E},\sigma}}\right)\max_{s\in\mathfrak{S}_{\nu}\setminus\{\mathrm{id}\}}{\left\{\vert 1+\mathsf{x}_{1,\ldots,\nu}\vert_{p},\vert\mathsf{x}_{s(1),\ldots,s(\nu)}\vert_{p}\right\}}.
$$Si, de plus, l'image de $\mathsf{x}$ dans $\det E\otimes_{\sigma}\mathbf{C}_{p}$ est nulle alors $\sum_{s\in\mathfrak{S}_{\nu}}{\varepsilon(s)\mathsf{x}_{s(1),\ldots,s(\nu)}}=0$. De cette relation et par in{\'e}galit{\'e} ultram{\'e}trique, l'on d{\'e}duit que le maximum ci-dessus est sup{\'e}rieur {\`a} $1$ puis que $$\Vert\mathsf{e}_{1}\wedge\cdots\wedge\mathsf{e}_{\nu}\Vert_{\overline{\det E},\sigma}\ge t^{\nu}\prod_{i=1}^{\nu}{\Vert\mathsf{e}_{i}\Vert_{\overline{E},\sigma}}.$$On a donc $$\frac{\Vert\mathsf{e}_{1}\wedge\cdots\wedge\mathsf{e}_{\nu}\Vert_{\overline{\det E_{\mathsf{a}}},\sigma}}{\Vert\mathsf{e}_{1}\wedge\cdots\wedge\mathsf{e}_{\nu}\Vert_{\overline{\det E},\sigma}}\le t^{-\nu}\moy_{\sigma}{(1,\Vert\mathsf{a}_{\sigma}\Vert)}^{\rho_{\sigma}}.$$Lorsque $\sigma$ est archim{\'e}dien, la m{\^{e}}me in{\'e}galit{\'e} reste valide en choisissant la base $(\mathsf{e}_{1},\ldots,\mathsf{e}_{\nu})$ orthonorm{\'e}e (la preuve se simplifie car $t=1$ convient). Ceci conduit {\`a} l'estimation $$-\widehat{\mu}(\overline{E}_{\mathsf{a}})+\widehat{\mu}(\overline{E})\le\frac{1}{\nu D}\sum_{\sigma\in S}{\rho_{\sigma}\log\moy_{\sigma}{(1,\Vert\mathsf{a}_{\sigma}\Vert)}}-\frac{\card S}{D}\log t.$$On conclut en faisant tendre $t$ vers $1$.
\section{D{\'e}monstration du th{\'e}or{\`e}me~\ref{theoremeprincipal}}\label{demonstrationtheoprincipal}
\subsection{Canevas de la d{\'e}monstration}\label{sec:canevasdelademonstration} Nous utilisons la m{\'e}thode de Baker non pas avec des fonctions auxiliaires, ni avec des d{\'e}terminants d'interpolation, ni avec la m{\'e}thode des pentes. En r{\'e}alit{\'e}, nous faisons un m{\'e}lange entre la premi{\`e}re et la derni{\`e}re de ces m{\'e}thodes, c'est-{\`a}-dire que nous passons par la construction d'une section auxiliaire d'un certain fibr{\'e} ad{\'e}lique hermitien. {\`A} chaque {\'e}tape de la d{\'e}monstration nous conservons l'aspect intrins{\`e}que des donn{\'e}es. De la sorte nous b{\'e}n{\'e}ficions de la souplesse de la m{\'e}thode des fonctions auxiliaires et de la possibilit{\'e} d'acc{\'e}der naturellement aux constantes num{\'e}riques de la m{\'e}thode des pentes. Concr{\`e}tement, apr{\`e}s avoir modifi{\'e} les donn{\'e}es initiales de mani{\`e}re {\`a} {\^e}tre en mesure d'utiliser la r{\'e}duction d'Hirata-Kohno et le proc{\'e}d{\'e} de changement de variables de Chudnovsky, nous construisons un fibr{\'e} ad{\'e}lique hermitien $\overline{E}$. L'espace vectoriel $E$ sous-jacent est un espace de polyn{\^o}mes en plusieurs variables, auquel sont adjointes des normes en tous les plongements de $k$. Ce fibr{\'e} ad{\'e}lique $\overline{E}$ a la particularit{\'e} d'avoir une $\sigma_{0}$-norme \og tordue\fg{}, ce qui constitue une des nouveaut{\'e}s de ce texte. Nous construisons alors un {\'e}l{\'e}ment $s\ne 0$ de $E$ de petite hauteur au moyen du lemme de Siegel absolu. Le choix des param{\`e}tres et le lemme de multiplicit{\'e}s de Philippon assurent qu'il existe un jet de $s$ le long de $W$ (sous-espace construit {\`a} partir de $W_{0}$) en un multiple de $\mathsf{p}=(u_{0},\alpha_{1},\ldots,\alpha_{n})$ qui est non nul. Ensuite nous {\'e}valuons la hauteur de ce jet, en distinguant les normes relatives aux plongements ultram{\'e}triques de celles relatives aux plongements archim{\'e}diens. Le comportement du jet en $\sigma_{0}$ et tous ses conjugu{\'e}s est {\'e}tudi{\'e} {\`a} part. La majoration de sa norme repose sur une extrapolation sur les d{\'e}rivations (cas p{\'e}riodique) ou sur les multiples de $\mathsf{p}$ (cas non p{\'e}riodique), qui a {\'e}t{\'e} pr{\'e}par{\'e}e par la construction de $s$. C'est {\`a} cet endroit qu'apparaissent les valeurs absolues des formes lin{\'e}aires que nous cherchons {\`a} {\'e}valuer. Pour conclure, nous utilisons une variante de l'in{\'e}galit{\'e} de Liouville (lemme~\ref{liouville}) pour minorer la hauteur de ce jet, qui fait intervenir la pente maximale arakelovienne de l'espace naturel dans lequel vit le jet consid{\'e}r{\'e}.  
\subsection{R{\'e}ductions}\label{paragraphereductions}
Pour d{\'e}montrer le th{\'e}or{\`e}me~\ref{theoremeprincipal}, l'on peut supposer que \begin{enumerate}\item[(i)] $\{u_{1},\ldots,u_{n}\}$ est une famille libre sur $\mathbf{Q}$,
\item[(ii)] $\{\ell_{1},\ldots,\ell_{t}\}$ est une famille libre du dual $(k^{n})^{\mathsf{v}}$,
\item[(iii)] $\vert(\beta_{1,0},\ldots,\beta_{t,0})\vert_{2,\sigma_{0}}\le 1$ lorsque $\sigma_{0}$ est archim{\'e}dien\footnote{Cette hypoth{\`e}se sera utile uniquement pour le cas p{\'e}riodique, \S~\ref{sectioncasperiodique}.},
\item[(iv)] si $n=1$ (et donc $t=1$) alors $\beta_{1,0}\ne 0$.
\end{enumerate}En effet, consid{\'e}rons un ensemble $I$ tel que $(u_{i})_{i\in I}$ soit une base du $\mathbf{Q}$-espace vectoriel engendr{\'e} par $u_{1},\ldots,u_{n}$. Posons $J=\{1,\ldots,n\}\setminus I$. Pour $j\in J$, la famille $\{u_{i}\}_{i\in I}\cup\{u_{j}\}$ est li{\'e}e sur $\mathbf{Q}$.  Le lemme 7.19 de~\cite[p.~$222$]{miw4} assure l'existence de nombres rationnels $\theta_{j,i}$ tels que $u_{j}=\sum_{i\in I}{\theta_{j,i}u_{i}}$ et \begin{equation*}h(\theta_{j,i})\le (n-1)\log(11(n-1)D^{3})+\log\prod_{l=1}^{n}{\log a_{l}}\end{equation*}pour tous $j,i$. En utilisant la d{\'e}finition de $\mathfrak{a}$, on a $$h(\theta_{j,i})\le(n-1)\log(11(n-1))+3(n-1)\left(\log\frac{D}{\log\mathfrak{e}}+\log\log\mathfrak{e}\right)+n\frac{\mathfrak{a}\log\mathfrak{e}}{D}$$puis, comme $\frac{\mathfrak{a}\log\mathfrak{e}}{D}\ge\max{(1,\log\frac{D}{\log\mathfrak{e}})}$, on trouve \begin{align*}h(\theta_{j,i})&\le\left((n-1)\log(11(n-1))+4n-3\right)\left(\frac{\mathfrak{a}\log\mathfrak{e}}{D}+\log\log\mathfrak{e}\right)\\ &\le 2n^{2}\left(\frac{\mathfrak{a}\log\mathfrak{e}}{D}+\log\log\mathfrak{e}\right).\end{align*} Par ailleurs, chaque $\ell_{\mathsf{j}}(u)$, $1\le\mathsf{j}\le t$, est la valeur d'une forme lin{\'e}aire $L_{\mathsf{j}}$ en les $u_{i}$, $i\in I$, avec des coefficients de la forme $$\beta_{\mathsf{j},i}':=\beta_{\mathsf{j},i}+\sum_{m\not\in I}{\beta_{\mathsf{j},m}\theta_{m,i}}.$$Pour un tel coefficient on a $$\vert\beta_{\mathsf{j},i}'\vert_{\sigma}\le n^{\epsilon_{\sigma}}\max{(1,\vert\beta_{\mathsf{j},i}\vert_{\sigma})}\prod_{m\not\in I}{\max{(1,\vert\beta_{\mathsf{j},m}\vert_{\sigma})}\max{(1,\vert\theta_{m,i}\vert_{\sigma})}}$$ce qui conduit {\`a} la borne de la hauteur \begin{align*}h(\beta_{\mathsf{j},i}')&\le\log n+n\max_{m,\ell}{\{h(\beta_{m,\ell})\}}+n\max_{m,i}{\{h(\theta_{m,i})\}}\\ &\le\frac{(n+\log n)\log b}{D}+2n^{3}\left(\frac{\mathfrak{a}\log\mathfrak{e}}{D}+\log\log\mathfrak{e}\right)\\ &\le\frac{2n^{3}}{D}\left(\log b+\mathfrak{a}\log\mathfrak{e}+D\log\log\mathfrak{e}\right).\end{align*}De cette famille $\{L_{1},\ldots,L_{t}\}$ de formes lin{\'e}aires sur $k^{I}$,  l'on peut extraire une famille libre maximale, qui comporte exactement $s=\dim\mathsf{T}_{u}-\dim(W_{0}\cap\mathsf{T}_{u})$ {\'e}l{\'e}ments. Notons $\mathsf{S}$ le sous-ensemble de $\{1,\ldots,t\}$ qui indexe la famille libre choisie. Pour tout $j\in\{1,\ldots,t\}$, on a $\Lambda_{j}=\beta_{j,0}+L_{j}((u_{i})_{i\in I})$ et le maximum des $\vert\Lambda_{j}\vert_{p_{0}}$ que l'on cherche {\`a} minorer est plus grand que $\max_{\mathsf{s}\in\mathsf{S}}{\{\vert\beta_{\mathsf{s},0}+L_{\mathsf{s}}((u_{i})_{i\in I})\vert_{p_{0}}\}}$. De la sorte l'on s'est donc bien ramen{\'e} aux conditions (i) et (ii) ci-dessus. Pour la condition (iii), il suffit de diviser tous les nombres $\beta_{i,j}$ par le premier entier sup{\'e}rieur {\`a} $\vert(\beta_{1,0},\ldots,\beta_{t,0})\vert_{2,\sigma_{0}}$. Cet entier est plus petit que $nb$. Enfin pour la condition~(iv), lorsque $n=1$ et $\beta_{1,0}=0$, il n'y a qu'une seule forme lin{\'e}aire $\Lambda_{1}=\beta_{1,1}u_{1}$. L'in{\'e}galit{\'e} de Liouville donne $\log\vert\beta_{1,1}\vert_{\sigma_{0}}\ge-Dh(\beta_{1,1})$. De plus, si $\sigma_{0}$ est ultram{\'e}trique, on a $\vert u_{1}\vert_{\sigma_{0}}=\vert\alpha_{1}-1\vert_{\sigma_{0}}$ et, si $\sigma_{0}$ est archim{\'e}dien, on a $$\vert\alpha_{1}-1\vert_{\sigma_{0}}=\left\vert u_{1}\int_{0}^{1}{e^{tu_{1}}\,\mathrm{d}t}\right\vert_{\sigma_{0}}\le\vert u_{1}\vert_{\sigma_{0}}\int_{0}^{1}{e^{t}\,\mathrm{d}t}\le 2\vert u_{1}\vert_{\sigma_{0}}$$lorsque $\vert u_{1}\vert_{\sigma_{0}}\le 1$. De plus, dans tous les cas, $\log\vert\alpha_{1}-1\vert_{\sigma_{0}}\ge-Dh(\alpha_{1}-1)\ge-D(h(\alpha_{1})+\log 2)$. On obtient donc \begin{align*}\log\vert\Lambda_{1}\vert_{\sigma_{0}}&\ge-\log b-D(\log a_{1}+2\log2)\\ &\le-2(\log b+\log\mathfrak{e})\left(1+\frac{D\log a_{1}}{\log\mathfrak{e}}\right)\end{align*}et le th{\'e}or{\`e}me~\ref{theoremeprincipal} est vrai. On peut donc supposer la condition~(iv) v{\'e}rifi{\'e}e.\par Les estimations faites des hauteurs des coefficients des $L_{\mathsf{s}}$ montrent que le th{\'e}or{\`e}me~\ref{theoremeprincipal} d{\'e}coule de l'{\'e}nonc{\'e} suivant (les notations sont celles du th{\'e}or{\`e}me~\ref{theoremeprincipal}).
\begin{theo}
\label{theoremereduit} Supposons que $\{u_{1},\ldots,u_{n}\}$ et $\{\ell_{1},\ldots,\ell_{t}\}$ sont des familles libres sur $\mathbf{Q}$. Alors, pour tout $j\in\{1,\ldots,t\}$, on a $\Lambda_{j}\ne 0$ et \begin{equation*}\log\max_{1\le j\le t}{\vert\Lambda_{j}\vert_{p_{0}}}\ge-(4n)^{90n^{2}}\mathfrak{a}^{1/t}(\log b+\mathfrak{a}\log\mathfrak{e})\prod_{j=1}^{n}{\left(1+\frac{D\log a_{j}}{\log\mathfrak{e}}\right)^{1/t}}.\end{equation*}De plus, si $t=1$ et $\beta_{1,0}\ne 0$, la quantit{\'e} $\log b+\mathfrak{a}\log\mathfrak{e}$ qui est dans le minorant peut {\^{e}}tre remplac{\'e}e par $\log b+\log\mathfrak{e}+D\log\mathfrak{a}$.\end{theo}
Le fait qu'aucun des $\Lambda_{j}$ n'est nul est une cons{\'e}quence du th{\'e}or{\`e}me de Baker qui affirme que la famille $\{1,u_{1},\ldots,u_{n}\}$ est libre sur $\overline{\mathbf{Q}}$ lorsque $\{u_{1},\ldots,u_{n}\}$ l'est sur $\mathbf{Q}$. La constante $(4n)^{91n^{2}}$ qui est dans le th{\'e}or{\`e}me~\ref{theoremeprincipal} est un majorant simple de $(4n)^{90n^{2}}(2n^{3}+1)$, le terme $2n^{3}$ venant de la majoration des hauteurs des $\beta_{\mathsf{j},i}'$ donn{\'e}e plus haut. Le th{\'e}or{\`e}me~\ref{theointro} donn{\'e} dans l'introduction est une cons{\'e}quence imm{\'e}diate de l'{\'e}nonc{\'e}~\ref{theoremereduit} en choisissant $a_{j}=a$ pour tout $j\in\{1,\ldots,n\}$. Il faut toutefois prendre garde que la valeur de $\mathfrak{a}$ dans le th{\'e}or{\`e}me~\ref{theointro} n'est pas tout {\`a} fait celle que l'on obtiendrait en particularisant la valeur de $\mathfrak{a}$ qui est dans le th{\'e}or{\`e}me~\ref{theoremeprincipal} car $\log\prod_{j=1}^{n}{a_{j}}=n\log a$. C'est la raison pour laquelle nous avons la constante $(4n)^{91n^{2}}$ dans le th{\'e}or{\`e}me~\ref{theointro}, constante qui majore $(4n)^{90n^{2}}\max{\{1,\log n\}}^{2}$.

\par La suite de l'article concerne la d{\'e}monstration du th{\'e}or{\`e}me~\ref{theoremereduit}. En particulier, dans toute la suite, nous supposerons que les familles $\{u_{1},\ldots,u_{n}\}$ et $\{\ell_{1},\ldots,\ell_{t}\}$ sont des familles libres sur $\mathbf{Q}$. 
\subsection{Pr{\'e}paratifs}\label{subsecpreparatifs}
Avant de commencer la d{\'e}monstration du th{\'e}or{\`e}me~\ref{theoremereduit}, il nous est utile de modifier les donn{\'e}es brutes du paragraphe~\ref{subsection:donnees}.\par Rappelons que $W_{0}$ d{\'e}signe le sous-espace vectoriel de $k^{n}$ intersection des noyaux des formes lin{\'e}aires $\ell_{i}$, $i\in\{1,\ldots,t\}$. Soit $G_{0}$ le spectre de l'alg{\`e}bre sym{\'e}trique $\mathbf{S}\left((k^{n}/W_{0})^{\mathsf{v}}\right)$. C'est un sch{\'e}ma en groupes affine sur $\spec k$, dont l'espace tangent \`a l'origine s'identifie canoniquement au quotient $k^{n}/W_{0}$. Notons $G$ le $k$-groupe alg{\'e}brique lin{\'e}aire $G_{0}\times\mathbf{G}_{\mathrm{m}}^{n}$. Soit $\lambda:k^{n}\to k^{n}/W_{0}$ la projection canonique et $W$ le sous-espace de l'espace tangent {\`a} l'origine $t_{G}=\left(k^{n}/W_{0}\right)\oplus k^{n}$ d{\'e}fini comme l'ensemble des vecteurs de la forme $(\lambda(y),y)$ avec $y\in k^{n}$. Le fibr{\'e} ad{\'e}lique hermitien $(k^{n},\vert\cdot\vert_{2})$ conf{\`e}re {\`a} $k^{n}/W_{0}$ et {\`a} $t_{G}$ des structures ad{\'e}liques hermitiennes, respectivement par quotient et par somme directe orthogonale. Nous noterons $\overline{k^{n}/W_{0}}$ et $\overline{t_{G}}$ les fibr{\'e}s ad{\'e}liques hermitiens ainsi obtenus. Comme sous-espace vectoriel de $(k^{n},\vert\cdot\vert_{2})$, l'espace $W_{0}$ h{\'e}rite d'une structure de fibr{\'e} ad{\'e}lique hermitien $\overline{W_{0}}$ et l'on note $h(\overline{W_{0}})$ sa hauteur d'Arakelov au sens du \S~\ref{sec:elements}. Par hermitianit{\'e}, cette quantit{\'e} est aussi le degr{\'e} d'Arakelov de $\overline{k^{n}/W_{0}}$. {\`A} une constante qui ne d{\'e}pend que de $n$ pr{\`e}s, elle est born{\'e}e par $\max_{i,j}{\{1,h(\beta_{i,j})\}}$, terme qui appara{\^\i}t dans la d{\'e}finition de $\log b$ du th{\'e}or{\`e}me~\ref{theoremeprincipal}. Plus pr{\'e}cis{\'e}ment on a la
\begin{prop}\label{propositionhauteurwo}La hauteur d'Arakelov de $\overline{W_{0}}$ v{\'e}rifie l'in{\'e}galit{\'e} $$D\max{\{1,h(\overline{W_{0}})\}}\le n^{3}\log b.$$
\end{prop}
\begin{proof}Si $n=t$ alors $W_{0}=\{0\}$ et la proposition est vraie. Nous supposons maintenant $t\le n-1$. L'application $k^{n}\to k^{t}$ qui {\`a} $z=(z_{1},\ldots,z_{n})$ associe $(\ell_{1}(z),\ldots,\ell_{t}(z))$ se factorise en un isomorphisme $\varphi:k^{n}/W_{0}\to k^{t}$ de $k$-espaces vectoriels. Pour tout $i\in\{1,\ldots,t\}$ et pour tout plongement $\sigma:k\hookrightarrow\mathbf{C}_{p}$, l'in{\'e}galit{\'e} de Cauchy-Schwarz fournit la majoration $$\vert\ell_{i}(z_{1},\ldots,z_{n})\vert_{\sigma}\le\Vert\ell_{i}\Vert_{(k^{n},\vert\cdot\vert_{2})^{\mathsf{v}},\sigma}\Vert\overline{z}\Vert_{\overline{k^{n}/W_{0}},\sigma}$$avec \begin{equation*}\Vert\ell_{i}\Vert_{(k^{n},\vert\cdot\vert_{2})^{\mathsf{v}},\sigma}=\moy_{p}(\vert\beta_{i,1}\vert_{\sigma},\ldots,\vert\beta_{i,n}\vert_{\sigma}).\end{equation*}Par cons{\'e}quent, la norme d'op{\'e}rateur $\Vert\varphi\Vert_{\sigma}$ de $\varphi$ entre $\overline{k^{n}/W_{0}}$ et $(k^{t},\vert\cdot\vert_{2})$ est plus petite que $\vert(\beta_{i,j})_{i,j}\vert_{2,\sigma}$. Ainsi on a \begin{equation*}\begin{split}h(\overline{W_{0}})=\widehat{\deg}_{\mathrm{n}}(\overline{k^{n}/W_{0}})&\le\widehat{\deg}_{\mathrm{n}}(k^{t},\vert\cdot\vert_{2})+\frac{t}{D}\sum_{p\in\EuScript{P}}{\sum_{\sigma:k\hookrightarrow\mathbf{C}_{p}}{\log\Vert\varphi\Vert_{\sigma}}}\\ &\le 0+th_{(k^{nt},\vert\cdot\vert_{2})}((\beta_{i,j})_{i,j})\\ &\le\frac{t}{2}\log(nt)+nt^{2}\frac{\log b}{D}\end{split}\end{equation*}(la premi{\`e}re in{\'e}galit{\'e} est une in{\'e}galit{\'e} de pentes classique, voir p.~ex. le lemme~$6.3$ de~\cite{rendiconti}). On conclut au moyen des majorations suivantes~:\begin{align*}\frac{t}{2}\log(nt)+nt^{2}&\le\frac{(n-1)}{2}\log(n(n-1))+n(n-1)^{2}\\ &\le n\log n+n^{3}-n(2n-1)\\ &\le n^{3}-n^{2}+n\le n^{3}\end{align*}car $\log n\le n$.\end{proof}
L'isomorphisme $\varphi$ d{\'e}fini ci-dessus permet de d{\'e}finir $u'_{0}:=\varphi^{-1}(u_{0})\in k^{n}/W_{0}$. \textbf{Dans la suite, nous confondrons $u'_{0}$ avec $u_{0}$.} Nous ne garderons que la notation $u_{0}$, le contexte permettant de distinguer s'il s'agit de $u'_{0}\in k^{n}/W_{0}$ ou bien de $u_{0}\in k^{t}$. \par Par ailleurs, plut{\^{o}}t que minorer $\max{\{\vert\Lambda_{i}\vert_{p_{0}};\ 1\le i\le t\}}$, nous allons minorer la norme de $u_{0}-\lambda(u)$, norme relative {\`a} $k^{n}/W_{0}$ au plongement $\sigma_{0}$. Ceci est rendu possible par le r{\'e}sultat suivant~:
\begin{prop}\label{propcomparaison}On a $\Vert u_{0}-\lambda(u)\Vert_{\overline{k^{n}/W_{0}},\sigma_{0}}\le b^{n^{3}}(\sqrt{t})^{\epsilon_{0}}\max_{1\le i\le t}{\vert\Lambda_{i}\vert_{p_{0}}}$.\end{prop}
\begin{proof}Si $n=t$ alors $W_{0}=\{0\}$ et $\Vert u_{0}-\lambda(u)\Vert_{\overline{k^{n}/W_{0}},\sigma_{0}}=\vert(\Lambda_{1},\ldots,\Lambda_{t})\vert_{2,\sigma_{0}}$. La proposition est vraie dans ce cas et, dans la suite, nous supposons $t\le n-1$. Comme $\varphi(u_{0}-\lambda(u))=-(\Lambda_{1},\ldots,\Lambda_{t})$ on a $$\Vert u_{0}-\lambda(u)\Vert_{\overline{k^{n}/W_{0}},\sigma_{0}}\le\Vert\varphi^{-1}\Vert_{\sigma_{0}}\vert(\Lambda_{1},\ldots,\Lambda_{t})\vert_{2,\sigma_{0}}\le\Vert\varphi^{-1}\Vert_{\sigma_{0}}(\sqrt{t})^{\epsilon_{0}}\max_{1\le i\le t}{\vert\Lambda_{i}\vert_{p_{0}}}$$o{\`u}, si $\sigma$ est un plongement de $k$, $\Vert\varphi^{-1}\Vert_{\sigma}$ d{\'e}signe la norme d'op{\'e}rateur de $\varphi^{-1}_{\sigma}:(\mathbf{C}_{p}^{t},\vert\cdot\vert_{2,\sigma})\to(k^{n}/W_{0}\otimes_{\sigma}\mathbf{C}_{p},\Vert\cdot\Vert_{\overline{k^{n}/W_{0}},\sigma})$. On a $\Vert\varphi^{-1}\Vert_{\sigma}\le\Vert\varphi\Vert_{\sigma}^{t-1}\Vert\det\varphi\Vert_{\sigma}^{-1}$ (voir~\cite[lemme~$7.2$]{rendiconti} p.~ex.) et $\Vert\det\varphi\Vert_{\sigma}\le\Vert\varphi\Vert_{\sigma}^{t}$ par $1\le\Vert\varphi\Vert_{\sigma}\Vert\varphi^{-1}\Vert_{\sigma}$. On en d{\'e}duit \begin{equation*}h(\overline{W_{0}})=\widehat{\deg}_{\mathrm{n}}\overline{k^{n}/W_{0}}=h(\det\varphi)\le\frac{t}{D}\sum_{\sigma}{\log\max{(1,\Vert\varphi\Vert_{\sigma})}}-\frac{\log\Vert\varphi^{-1}\Vert_{\sigma_{0}}}{D}.\end{equation*}De plus $h(\overline{W_{0}})\ge 0$ car si $e_{i_{1}},\ldots,e_{i_{t}}$ sont des vecteurs de la base canonique de $k^{n}$ dont les r{\'e}ductions $\overline{e_{i_{1}}},\ldots,\overline{e_{i_{t}}}$ modulo $W_{0}$ forment une base de $k^{n}/W_{0}$, on a \begin{align*}h(\overline{W_{0}})&=-\frac{1}{D}\sum_{p\in\EuScript{P}}{\sum_{\sigma:k\hookrightarrow\mathbf{C}_{p}}{\log\Vert\overline{e_{i_{1}}}\wedge\cdots\wedge\overline{e_{i_{t}}}\Vert_{\overline{\wedge^{t}(k^{n}/W_{0})},\sigma}}}\\ &\ge -\frac{1}{D}\sum_{p\in\EuScript{P}}{\sum_{\sigma:k\hookrightarrow\mathbf{C}_{p}}{\log(\Vert\overline{e_{i_{1}}}\Vert_{\overline{k^{n}/W_{0}},\sigma}\cdots\Vert\overline{e_{i_{t}}}\Vert_{\overline{k^{n}/W_{0}},\sigma})}}\\ &\ge-\frac{1}{D}\sum_{p\in\EuScript{P}}{\sum_{\sigma:k\hookrightarrow\mathbf{C}_{p}}{\log(\vert e_{i_{1}}\vert_{2,\sigma}\cdots\vert e_{i_{t}}\vert_{2,\sigma})}}=0.\end{align*}Par ailleurs la norme $\Vert\varphi\Vert_{\sigma}$ est inf{\'e}rieure {\`a} $\vert(\beta_{i,j})_{i,j}\vert_{2,\sigma}$. On en d{\'e}duit \begin{align*}\log\Vert\varphi^{-1}\Vert_{\sigma_{0}}^{1/D}&\le \frac{t}{D}\sum_{p\in\EuScript{P}}{\sum_{\sigma:k\hookrightarrow\mathbf{C}_{p}}{\log\max{(1,\vert(\beta_{i,j})_{i,j}\vert_{2,\sigma})}}}\\ & \le t(\log\sqrt{nt}+nt)\frac{\log b}{D}.\end{align*}On conclut en majorant $t(\log\sqrt{nt}+nt)$ par $n^{3}$ comme nous l'avons fait {\`a} la fin de la proposition~\ref{propositionhauteurwo}, en tenant compte de $t\le n-1$. 
\end{proof}

\subsection{Choix des param{\`e}tres}\label{section:choixdesparametres}
 Soit $C_{0}:=(4n)^{10n}$. Posons $y:=0$ si $t=1$ et $\beta_{1,0}\ne 0$ et $y:=1$ sinon. Soit $S_{0}:=C_{0}\mathfrak{a}$ et $S:=C_{0}^{3}\mathfrak{a}$. Soit $U_{0}>0$ un nombre r{\'e}el d{\'e}fini un peu plus loin (voir~\eqref{formuleuo}). Soit $\widetilde{D}_{0},\ldots,\widetilde{D}_{n},\widetilde{T},\widetilde{T}_{0}$ les nombres r{\'e}els donn{\'e}s par les formules suivantes~:
\begin{equation*}\widetilde{T}_{0}:=\frac{C_{0}U_{0}}{S\log\mathfrak{e}},\qquad\widetilde{T}:=C_{0}^{2}\widetilde{T}_{0},\end{equation*}\begin{equation*}\widetilde{D}_{0}:=\frac{U_{0}}{\log b+D\log S+S^{y}\log\mathfrak{e}}\end{equation*}(la pr{\'e}sence du param{\`e}tre $y$ au d{\'e}nominateur explique le raffinement donn{\'e} dans le th{\'e}or{\`e}me~\ref{theoremereduit} lorsqu'il n'y a qu'une seule forme lin{\'e}aire, non homog{\`e}ne) et\begin{equation*}\forall\,i\in\{1,\ldots,n\},\quad\widetilde{D}_{i}:=\frac{U_{0}}{S\log\mathfrak{e}+DS\log a_{i}} \cdotp\end{equation*}Notons $\overline{k}$ une cl{\^{o}}ture alg{\'e}brique de $k$ dans $\mathbf{C}_{p_{0}}$. Un sous-groupe alg{\'e}brique connexe $G'$ de $G$ se d{\'e}compose sur $\overline{k}$ en un produit $G_{0}'\times\mathbf{G}_{\mathrm{m}}'$ avec $G_{0}'$ (\emph{resp}. $\mathbf{G}_{\mathrm{m}}'$) un sous-groupe alg{\'e}brique connexe de $G_{0}$ (\emph{resp}. $\mathbf{G}_{\mathrm{m}}^{n}$). Posons \begin{equation}\label{eq:dimensions}t':=\dim G_{0}',\quad n':=\dim\mathbf{G}_{\mathrm{m}}',\quad r':=\codim_{G}G',\quad \lambda':=\codim_{W}(W\cap t_{G'}).\end{equation}De simples consid{\'e}rations d'alg{\`e}bre lin{\'e}aire montrent que \begin{equation}\label{alglineaire}r'-\lambda'=t-\dim(t_{G_{0}'}+\lambda(t_{\mathbf{G}_{\mathrm{m}}'}))\le\min{\{t-t',n-n'\}}.\end{equation}Par ailleurs, les groupes alg{\'e}briques $G_{0}$ et $G$ admettent des compactifications naturelles $X_{0}:=\mathbf{P}(k\oplus(k^{n}/W_{0})^{\mathsf{v}})$ et $X:=X_{0}\times(\mathbf{P}_{k}^{1})^{n}$. {\`A} une sous-vari{\'e}t{\'e} $V$ de $X_{\overline{k}}$ est associ{\'e} un polyn{\^o}me dit de Hilbert-Samuel $H_{V}$ \`a $n+1$ variables (en prenant l'adh{\'e}rence de Zariski de $V$ dans $X_{\overline{k}}$). Soit $\EuScript{H}(V;X_{0},\ldots,X_{n})$ la partie homog{\`e}ne de plus haut degr{\'e} de $H_{V}$ multipli{\'e}e par $(\dim V)!$. Les coefficients de $\EuScript{H}$ sont des entiers positifs de somme {\'e}gale au degr{\'e} de $V$. Par exemple, on a $\EuScript{H}(G;X_{0},\ldots,X_{n})=\frac{(n+t)!}{t!}X_{0}^{t}X_{1}\cdots X_{n}$ (pour les propri{\'e}t{\'e}s de $\EuScript{H}$ nous renvoyons au texte de Roy~\cite[chapitre~$5$]{miw4}). Posons $$\mathsf{p}:=(u_{0},\alpha_{1},\ldots,\alpha_{n})\quad\text{et}\quad\Sigma_{\mathsf{p}}(S):=\{0_{G},\mathsf{p},2\mathsf{p},\ldots,S\mathsf{p}\}.$$Lorsque $W+t_{G'}\ne t_{G}$ (et en particulier $t'\ne t$), on pose \begin{equation*}x(G'):=\left(\frac{\widetilde{T}^{\lambda'}\card\left(\frac{\Sigma_{\mathsf{p}}(S)+G'(\mathbf{C}_{p_{0}})}{G'(\mathbf{C}_{p_{0}})}\right)\EuScript{H}(G';\widetilde{D}_{0},\ldots,\widetilde{D}_{n})}{C_{0}\EuScript{H}(G;\widetilde{D}_{0},\ldots,\widetilde{D}_{n})}\right)^{\frac{1}{t-t'}} \cdotp\end{equation*}Soit $x$ le minimum des nombres r{\'e}els $x(G')$ lorsque $G'$ parcourt les sous-groupes alg{\'e}briques connexes de $G$ tels que $W+t_{G'}\ne t_{G}$ (l'existence de $x$ d{\'e}coule d'un argument standard, expliqu{\'e} par exemple {\`a} la page~$734$ de~\cite{artepredeux}). Soit $\widetilde{G}=\widetilde{G}_{0}\times\widetilde{\mathbf{G}}_{\mathrm{m}}$ un sous-groupe alg{\'e}brique v{\'e}rifiant ces conditions et tel que $x=x(\widetilde{G})$, auquel on adjoint les entiers $\widetilde{t},\widetilde{n},\widetilde{r},\widetilde{\lambda}$ d{\'e}finis par~\eqref{eq:dimensions}. N{\'e}cessairement on a $\widetilde{n}<n$ car sinon $t_{\widetilde{G}}+W=t_{G}$. Notons $T:=[\widetilde{T}]$, $D_{0}:=[x\widetilde{D}_{0}]$ et, pour tout $i\in\{1,\ldots,n\}$, $D_{i}:=[\widetilde{D}_{i}]$. L'homog{\'e}n{\'e}it{\'e} de la fonction $\EuScript{H}$ permet de voir que la quantit{\'e} $U_{0}^{(r'-\lambda')/(t-t')}x(G')$ ne d{\'e}pend pas de $U_{0}$, lorsque $G'$ varie parmi les sous-groupes autoris{\'e}s. En observant que $\{0\}$ fait partie de ces sous-groupes, nous pouvons choisir $U_{0}$ de sorte que $x(\{0\})\le 1$ (et, en particulier, on a $x\in\,]0,1]$). Concr{\`e}tement, la condition $x(\{0\})\le 1$ est satisfaite si l'on choisit $U_{0}$ {\'e}gal {\`a}\begin{equation}\label{formuleuo}C_{0}^{(3n-1)/t}\left\{\frac{t!\card(\Sigma_{\mathsf{p}}(S))}{(n+t)!}\prod_{i=1}^{n}{\left(1+\frac{D\log a_{i}}{\log\mathfrak{e}}\right)}\right\}^{1/t}\times (\log b+D\log S+S^{y}\log\mathfrak{e}).\end{equation}On notera que $\card(\Sigma_{\mathsf{p}}(S))=S+1$ sauf, {\'e}ventuellement, si toutes les conditions $p_{0}=\infty$, $n=1$, $\beta_{1,0}=0$ et $u_{1}\in\mathbf{Q}i\pi$ sont remplies. Mais ce cas est exclu par la r{\'e}duction~(iv) faite au d{\'e}but du \S~\ref{paragraphereductions}. On a donc $\card(\Sigma_{\mathsf{p}}(S))=S+1$. L'expression~\eqref{formuleuo} de $U_{0}$ permet de justifier les estimations suivantes.
\begin{prop}\label{propparametres}
Les propri{\'e}t{\'e}s suivantes sont satisfaites~:\begin{enumerate}\item[(i)]
 $C_{0}\max{\{1,\widetilde{D}_{0}/S^{1-y},\widetilde{D}_{1},\ldots,\widetilde{D}_{n}\}}\le\widetilde{T}_{0}$,
\item[(ii)] $D_{0}\ne 0$ et $C_{0}\le\max_{1\le j\le n}{\{\widetilde{D}_{j}\}}$ (en particulier l'un au moins des $D_{j}$, $1\le j\le n$, n'est pas nul),
\item[(iii)] $D\log\mathfrak{a}\le2\mathfrak{a}\log\mathfrak{e}$,
\item[(iv)] $T\log(4\widetilde{D}_{0})\le (10n\log C_{0})U_{0}/D$.
\end{enumerate}
\end{prop} 
\begin{proof} Si les in{\'e}galit{\'e}s (i) et (ii) se v{\'e}rifient ais{\'e}ment {\`a} partir des valeurs des param{\`e}tres, il n'est pas enti{\`e}rement {\'e}vident que $D_{0}$ est non nul, en raison du $x$ devant $\widetilde{D}_{0}$. Pour voir cela, on reprend l'argumentation du lemme 5.1, (iii), de~\cite{artepredeux} qui repose sur une propri{\'e}t{\'e} de d{\'e}croissance d'un quotient de fonctions $\EuScript{H}$  et sur l'in{\'e}galit{\'e} $\widetilde{\lambda}\ge n-\widetilde{n}\ge 1$. On a \begin{equation*}\begin{split}\left(x\widetilde{D}_{0}\right)^{t-\widetilde{t}}&\ge\frac{\widetilde{T}^{\widetilde{\lambda}}\EuScript{H}(\widetilde{G};\widetilde{D}_{0},\ldots,\widetilde{D}_{n})\widetilde{D}_{0}^{t-\widetilde{t}}}{C_{0}\EuScript{H}(G;\widetilde{D}_{0},\ldots,\widetilde{D}_{n})}\\ &\ge\frac{\widetilde{T}^{\widetilde{\lambda}}}{C_{0}\binom{n+t}{t}\max_{1\le i\le n}{\{\widetilde{D}_{i}\}}^{n-\widetilde{n}}}\\ &\ge\frac{C_{0}^{3\widetilde{\lambda}-1}}{(n+t)!}\quad(\text{gr{\^{a}}ce {\`a} (i)}).
\end{split}\end{equation*}La valeur de $C_{0}$ entra{\^{\i}}ne alors $(x\widetilde{D}_{0})^{t-\widetilde{t}}\ge 1$ puis $D_{0}\ge 1$. Venons-en maintenant {\`a} (iii) et (iv) qui requi{\`e}rent aussi quelques d{\'e}tails. Tout d'abord, on a \begin{equation*}\log\mathfrak{a}=\log\left(\frac{\mathfrak{a}\log\mathfrak{e}}{D}\frac{D}{\log\mathfrak{e}}\right)\le \frac{2\mathfrak{a}\log\mathfrak{e}}{D}\end{equation*}(car $\log x\le x$ si $x\ge 1$), ce qui donne (iii). Par ailleurs, les d{\'e}finitions de $T$ et $\widetilde{D}_{0}$ combin{\'e}es avec la valeur~\eqref{formuleuo} de $U_{0}$ donn{\'e}e ci-dessus impliquent \begin{equation*}T\log(4\widetilde{D}_{0})\le\frac{U_{0}}{\mathfrak{a}\log\mathfrak{e}}\log\left(C_{0}^{3n+3}\mathfrak{a}\prod_{i=1}^{n}{\left(1+\frac{D\log a_{i}}{\log\mathfrak{e}}\right)}\right)\cdotp
\end{equation*}En majorant le produit qui est dans le logarithme du majorant par \begin{equation*}\left(1+\frac{D}{\log\mathfrak{e}}\log\prod_{i=1}^{n}{a_{i}}\right)^{n}\le\left(1+\exp{\left\{\frac{2\mathfrak{a}\log\mathfrak{e}}{D}\right\}}\right)^{n}\end{equation*}on a \begin{equation*}T\log(4\widetilde{D}_{0})\le\frac{U_{0}}{\mathfrak{a}\log\mathfrak{e}}\left((3n+3)\log C_{0} +\log\mathfrak{a}+3n\frac{\mathfrak{a}\log\mathfrak{e}}{D}\right)\cdotp\end{equation*}On utilise alors (iii) pour majorer $\log\mathfrak{a}$ et conclure.
\end{proof}

D{\'e}signons par $\Omega_{\sigma_{0}}$ l'ensemble $\{0\}\times(2i\pi\mathbf{Z})^{n}$ si $\sigma_{0}$ est archim{\'e}dien et l'ensemble $\{0\}$ sinon.
\begin{defi}\label{defi:casperiodique} Nous dirons que nous sommes dans le \emph{cas p{\'e}riodique} s'il existe un entier $m\in\{1,\ldots,(n+t)S\}$ tel que $m(u_{0},u)\in t_{\widetilde{G}}(\mathbf{C}_{p_{0}})+\Omega_{\sigma_{0}}$, et que nous sommes dans le \emph{cas non p{\'e}riodique} dans le cas contraire. \end{defi} Dans le cas p{\'e}riodique, la lettre $\Upsilon$ d{\'e}signe l'ensemble des couples $(m,\tau)$ avec $m\in\{0,1,\ldots,(n+t)S\}$ et $\tau=(\tau_{1},\ldots,\tau_{n})\in\mathbf{N}^{n}$ qui v{\'e}rifie $\sum_{i=1}^{n}{\tau_{i}}\le2(n+t)T$ et $\tau_{n}\le T_{0}$. Dans le cas non p{\'e}riodique, $\Upsilon$ est l'ensemble des couples $(m,\tau)$ avec $m\in\{0,1,\ldots,S_{0}-1\}$ et $\tau=(\tau_{1},\ldots,\tau_{n})\in\mathbf{N}^{n}$ qui v{\'e}rifie $\sum_{i=1}^{n}{\tau_{i}}\le2(n+t)T$. Dans les deux cas, nous notons $\mu$ le cardinal de $\Upsilon$.
\begin{enonce}{Remarque}\label{remarquecinqsix}Comme la famille $\{u_{1},\ldots,u_{n}\}$ est libre sur $\mathbf{Q}$, {\^{e}}tre dans le cas p{\'e}riodique implique que le plongement $\sigma_{0}$ est archim{\'e}dien et que la famille $\{2i\pi,u_{1},\ldots,u_{n}\}$ est li{\'e}e sur $\mathbf{Q}$ (c'est-{\`a}-dire $\alpha_{1},\ldots,\alpha_{n}$ multiplicativement d{\'e}pendants). Comme il ne peut exister qu'une seule relation, {\`a} multiplication par un scalaire non nul pr{\`e}s, entre $2i\pi$ et les $u_{j}$ (sinon existerait une relation non triviale entre les $u_{j}$), le groupe $\widetilde{\mathbf{G}}_{\mathrm{m}}$ est de dimension $n-1$ et son espace tangent contient $W_{0}$.
\end{enonce}
\subsection{Charni{\`e}re de la d{\'e}monstration et compl{\'e}ments au cas p{\'e}riodique}\label{subseccomplements}Le groupe alg{\'e}brique $\widetilde{\mathbf{G}}_{\mathrm{m}}$ est de dimension $\widetilde{n}\le n-1$. En particulier, comme $\{u_{1},\ldots,u_{n}\}$ est libre sur $\mathbf{Q}$ on a $u\not\in t_{\widetilde{\mathbf{G}}_{\mathrm{m}}}(\mathbf{C}_{p_{0}})$. Par cons{\'e}quent les vecteurs $(u_{0},u)$ et $(\lambda(u),u)$ n'appartiennent pas {\`a} l'espace tangent du groupe alg{\'e}brique $\widetilde{G}(\mathbf{C}_{p_{0}})$. L'intersection $W\cap t_{\widetilde{G}}$ est donc un sous-espace strict de $W$ (car ce dernier contient $(\lambda(u),u)$). Si $\sigma_{0}$ est archim{\'e}dien, fixons alors une base \emph{orthonorm{\'e}e} $w:=(w_{1},\ldots,w_{n})$ de $W\otimes_{\sigma_{0}}\mathbf{C}_{p_{0}}$ qui provient de $W\otimes_{\sigma_{0}}\mathbf{R}$ lorsque $\sigma_{0}$ est un plongement r{\'e}el et qui poss{\`e}de en outre les deux propri{\'e}t{\'e}s suivantes~:\begin{enumerate}\item[(i)] $(w_{1},\ldots,w_{\mathfrak{n}})$ est une base de $(W\cap t_{\widetilde{G}})\otimes_{\sigma_{0}}\mathbf{C}_{p_{0}}$ (le $\mathfrak{n}$ gothique est la dimension de $W\cap t_{\widetilde{G}}$),
\item[(ii)] si $\mathfrak{u}:=(\mathfrak{u}_{1},\ldots,\mathfrak{u}_{n})\in\mathbf{C}_{p_{0}}^{n}$ d{\'e}signe le vecteur des coordonn{\'e}es de $(\lambda(u),u)$ dans la base $(w_{1},\ldots,w_{n})$ alors $\vert\mathfrak{u}_{n}\vert_{p_{0}}=\max{\{\vert\mathfrak{u}_{j}\vert_{p_{0}};\ \mathfrak{n}+1\le j\le n\}}$.\end{enumerate}Dans le cas p{\'e}riodique (et cette hypoth{\`e}se implique maintenant que n{\'e}cessairement $\sigma_{0}$ est archim{\'e}dien puisque $(u_{0},u)\not\in t_{\widetilde{G}}(\mathbf{C}_{p_{0}})$), l'on sait minorer $\vert\mathfrak{u}_{n}\vert_{p_{0}}$ de la mani{\`e}re suivante. Il existe un entier $m\in\{1,\ldots,(n+t)S\}$, $x\in t_{\widetilde{G}}(\mathbf{C})$ et $\omega\in(2i\pi\mathbf{Z})^{n}\setminus\{0\}$ tels que $m(u_{0},u)=x+(0,\omega)$. Soit $\mathrm{d}_{\sigma_{0}}$ la distance sur $t_{G}(\mathbf{C}_{p_{0}})$ induite par la norme $\Vert\cdot\Vert_{\overline{t_{G}},\sigma_{0}}$. L'in{\'e}galit{\'e} de Cauchy-Schwarz montre que $\vert\mathfrak{u}_{n}\vert_{p_{0}}$ est minor{\'e} par \begin{equation*}\begin{split}\frac{\mathrm{d}_{\sigma_{0}}(m(\lambda(u),u),t_{\widetilde{G}}(\mathbf{C}))}{m\sqrt{n}} & =\frac{\mathrm{d}_{\sigma_{0}}((m(\lambda(u)-u_{0}),\omega),t_{\widetilde{G}}(\mathbf{C}))}{m\sqrt{n}}\\ & \ge\frac{\mathrm{d}_{\sigma_{0}}(\omega,t_{\widetilde{\mathbf{G}}_{\mathrm{m}}}(\mathbf{C}))}{(n+t)S\sqrt{n}}\cdotp\end{split}\end{equation*}
{\`A} ce stade, nous aurons besoin d'un r{\'e}sultat que l'on trouve en substance dans~\cite{bertrandphilippon} et de mani{\`e}re plus explicite dans la d{\'e}monstration de la proposition~$6.1$ de~\cite{gaelcompositio}. Pour $x\in\mathbf{R}^{n}$, on note $\vert x\vert_{1}$ la somme des valeurs absolues des composantes de $x$.
\begin{lemm} Soit $H$ un sous-groupe alg{\'e}brique connexe de $\mathbf{G}_{\mathrm{m}}^{n}$, de dimension $h\in\{1,\ldots,n-1\}$. Alors il existe une famille libre $\{\mathsf{h}_{i}:=(h_{i,1},\ldots,h_{i,n})\,;\ 1\le i\le n-h\}$ de $\mathbf{Z}^{n}$ telle que \begin{itemize}
\item[(i)] $H=\{(x_{1},\ldots,x_{n})\in\mathbf{G}_{\mathrm{m}}^{n}\,;\ \forall\,i\in\{1,\ldots,n-h\},\ \prod_{j=1}^{n}{x_{j}^{h_{i,j}}}=1\}$
\item[(ii)] $\prod_{i=1}^{n-h}{\vert\mathsf{h}_{i}}\vert_{1}\le\frac{(n-h)!}{h!}\deg H$ (le degr{\'e} est relatif au plongement usuel $\mathbf{G}_{\mathrm{m}}^{n}\hookrightarrow(\mathbf{P}^{1})^{n}$).
\end{itemize}
\end{lemm}
En utilisant ce lemme, nous allons minorer $\mathrm{d}_{\sigma_{0}}(\omega,t_{\widetilde{\mathbf{G}}_{\mathrm{m}}}(\mathbf{C}))$ de la mani{\`e}re suivante. Ce lemme appliqu{\'e} {\`a} $H=\widetilde{\mathbf{G}}_{\mathrm{m}}$, qui est de dimension $n-1$ d'apr{\`e}s la remarque~\ref{remarquecinqsix}, donne l'existence d'un vecteur $\mathsf{h}$ de $\mathbf{Z}^{n}$ qui forme une base de l'orthogonal de $t_{\widetilde{\mathbf{G}}_{\mathrm{m}}}(\mathbf{C})$ dans $(\mathbf{C}^{n},\vert\cdot\vert_{2,\sigma_{0}})$ avec $\vert\mathsf{h}\vert_{1}\le(\deg\widetilde{\mathbf{G}}_{\mathrm{m}})/(n-1)!$. Ainsi la projection $\omega'$ de $\omega$ sur cet orthogonal est de norme {\'e}gale {\`a} $\mathrm{d}_{\sigma_{0}}(\omega,t_{\widetilde{\mathbf{G}}_{\mathrm{m}}}(\mathbf{C}))$. Par Cauchy-Schwarz, le produit hermitien $\omega.\mathsf{h}$ est de valeur absolue inf{\'e}rieure {\`a} $\vert\mathsf{h}\vert_{2,\sigma_{0}}\vert\omega'\vert_{2,\sigma_{0}}$. De plus ce produit est non nul car $\omega\not\in t_{\widetilde{\mathbf{G}}_{\mathrm{m}}}(\mathbf{C})$ et il appartient {\`a} $2i\pi\mathbf{Z}$. Comme $\vert\mathsf{h}\vert_{2,\sigma_{0}}\le\vert\mathsf{h}\vert_{1}$, on en d{\'e}duit l'in{\'e}galit{\'e} \begin{equation*}\mathrm{d}_{\sigma_{0}}(\omega,t_{\widetilde{\mathbf{G}}_{\mathrm{m}}}(\mathbf{C}))\ge\frac{2\pi(n-1)!}{\deg\widetilde{\mathbf{G}}_{\mathrm{m}}}
\end{equation*}puis $$\vert\mathfrak{u}_{n}\vert\ge\left(\frac{2\pi(n-1)!}{(n+t)\sqrt{n}}\right)\frac{1}{S\deg\widetilde{\mathbf{G}}_{\mathrm{m}}}\ge\left(\frac{\pi(n-1)!}{n\sqrt{n}}\right)\frac{1}{S\deg\widetilde{\mathbf{G}}_{\mathrm{m}}}$$car $t\le n$. On v{\'e}rifie alors que le premier quotient {\`a} droite est sup{\'e}rieur {\`a} $1$ ce qui conduit {\`a} l'{\'e}nonc{\'e} suivant.
\begin{prop}\label{propperiodique}Dans le cas p{\'e}riodique, on a $\vert\mathfrak{u}_{n}\vert^{-1}\le S\deg\widetilde{\mathbf{G}}_{\mathrm{m}}$.
\end{prop}Pour {\^{e}}tre utile, cette proposition devra {\^{e}}tre coupl{\'e}e avec le r{\'e}sultat suivant.
\begin{prop}\label{propannexe}Dans le cas p{\'e}riodique, on a $\deg\widetilde{\mathbf{G}}_{\mathrm{m}}\le C_{0}^{-2}\binom{n+t}{t}\widetilde{D}_{0}$.
\end{prop}
\begin{proof}
Par d{\'e}finition de $\widetilde{G}$, on a $x(\widetilde{G})\le x(\{0\})\le 1$ et la d{\'e}finition de $x(\widetilde{G})$ entra{\^{\i}}ne alors la majoration \begin{equation}\label{ineqinter}\frac{\widetilde{T}^{\widetilde{\lambda}}\EuScript{H}(\widetilde{G};\widetilde{D}_{0},\ldots,\widetilde{D}_{n})}{C_{0}\EuScript{H}(G;\widetilde{D}_{0},\ldots,\widetilde{D}_{n})}\le 1.\end{equation}En d{\'e}composant la fonction $\EuScript{H}$ sur les parties $G_{0}$ et $\mathbf{G}_{\mathrm{m}}^{n}$ et en utilisant la d{\'e}croissance des fonctions partielles $$x_{i}\mapsto\frac{\EuScript{H}(\widetilde{G};\widetilde{D}_{0},\ldots,x_{i},\ldots,\widetilde{D}_{n})}{\EuScript{H}(G;\widetilde{D}_{0},\ldots,x_{i},\ldots,\widetilde{D}_{n})},$$l'in{\'e}galit{\'e}~\eqref{ineqinter} devient \begin{equation}\label{trey}\widetilde{T}^{\widetilde{\lambda}}\deg\widetilde{\mathbf{G}}_{\mathrm{m}}\le C_{0}\binom{n+t}{t}(\widetilde{D}_{0})^{t-\widetilde{t}}\max{\left(\widetilde{D}_{1},\ldots,\widetilde{D}_{n}\right)}^{n-\widetilde{n}}.\end{equation}D'apr{\`e}s la remarque~\ref{remarquecinqsix}, on a $\dim\widetilde{\mathbf{G}}_{\mathrm{m}}=\widetilde{n}=n-1$. La formule~\eqref{alglineaire} donne $\widetilde{\lambda}=\dim(t_{\widetilde{G}_{0}}+\lambda(t_{\widetilde{\mathbf{G}}_{\mathrm{m}}}))-\widetilde{t}+1$ et $t-\dim(t_{\widetilde{G}_{0}}+\lambda(t_{\widetilde{\mathbf{G}}_{\mathrm{m}}}))\le 1$. De~\eqref{trey} on tire alors la majoration $$\left(\frac{\widetilde{T}}{\max{\left(\widetilde{D}_{0},\ldots,\widetilde{D}_{n}\right)}}\right)^{\widetilde{\lambda}}\le\frac{C_{0}\binom{n+t}{t}\widetilde{D}_{0}}{\deg\widetilde{\mathbf{G}}_{\mathrm{m}}}\cdotp$$Si $y=1$ la proposition~\ref{propannexe} d{\'e}coule de $\widetilde{T}\ge C_{0}^{3}\max{\left(\widetilde{D}_{0},\ldots,\widetilde{D}_{n}\right)}$ (proposition~\ref{propparametres}, (i)). Si $y=0$ on a $t=1$ par d{\'e}finition de $y$ et, en utilisant $\widetilde{\lambda}\ge 1$ et $\widetilde{t}\ge 0$, l'in{\'e}galit{\'e}~\eqref{trey} se simplifie en $$\frac{\widetilde{T}}{\max{\left(\widetilde{D}_{1},\ldots,\widetilde{D}_{n}\right)}}\le\frac{C_{0}\binom{n+t}{t}\widetilde{D}_{0}}{\deg\widetilde{\mathbf{G}}_{\mathrm{m}}}\cdotp$$Le minorant est sup{\'e}rieur {\`a} $C_{0}^{3}$ et la proposition s'en d{\'e}duit.\end{proof}

\subsection{Lemme de multiplicit{\'e}s}\label{paragraphemultiplicite} Au \S~\ref{paragraphepuissance}, nous avons introduit la notion de polyn{\^{o}}me multihomog{\`e}ne. Elle va nous permettre de d{\'e}finir un espace vectoriel dans lequel sera la section auxiliaire que nous allons construire au \S~\ref{secconstructiondes}. Soit $E_{0}:=k\oplus(k^{n}/W_{0})$ et, pour $i\in\{1,\ldots,n\}$, $E_{i}:=k\oplus k$. Posons $$E:=S^{D_{0}}(E_{0}^{\mathsf{v}})\otimes_{k}\bigotimes_{i=1}^{n}{S^{D_{i}}(E_{i}^{\mathsf{v}})}.$$Pour un plongement $\sigma:k\hookrightarrow\mathbf{C}_{p}$, soit $s\in E\otimes_{\sigma}\mathbf{C}_{p}$ et consid{\'e}rons $F_{s,\sigma}:(k^{n}/W_{0})\otimes_{\sigma}\mathbf{C}_{p}\times\prod_{i=1}^{n}{\mathcal{T}_{p}}\to\mathbf{C}_{p}$ l'application qui {\`a} $(z_{0},z_{1},\ldots,z_{n})$  associe \begin{equation}\label{defidefs}F_{s,\sigma}(z_{0},z_{1},\ldots,z_{n})=s((1,z_{0}),(1,e^{z_{1}}),\ldots,(1,e^{z_{n}}))\end{equation}(dans cette {\'e}criture nous avons identifi{\'e} $k$ et son dual \textit{via} la base canonique).
\begin{defi} Soit $\ell\in\mathbf{N}$ et $x=(x_{0},x_{1},\ldots,x_{n})\in G(\mathbf{C}_{p})$. On dit qu'un polyn{\^{o}}me $s\in E\otimes_{\sigma}\mathbf{C}_{p}$ s'annule au point $x$ {\`a} l'ordre $\ell$ le long de $W$ si l'application $$(z_{1},\ldots,z_{n})\mapsto s((1,x_{0}+\lambda(z_{1},\ldots,z_{n})),(1,x_{1}e^{z_{1}}),\ldots,(1,x_{n}e^{z_{n}}))$$s'annule {\`a} l'ordre $\ell$ en $(0,\ldots,0)$, \emph{i.\,e.} si la s{\'e}rie formelle en les variables $z_{1},\ldots,z_{n}$ d{\'e}finie par cette application est un {\'e}l{\'e}ment de l'id{\'e}al $(z_{1},\ldots,z_{n})^{\ell}$ de $\mathbf{C}_{p}[[z_{1},\ldots,z_{n}]]$.
\end{defi}Si le point $x$ est l'image du vecteur $(\mathsf{z}_{0},\mathsf{z}_{1},\ldots,\mathsf{z}_{n})\in(k^{n}/W_{0})\otimes_{\sigma}\mathbf{C}_{p}\times\prod_{i=1}^{n}{\mathcal{T}_{p}}$ par l'application exponentielle de $G(\mathbf{C}_{p})$ alors dire que $s$ s'annule au point $x$ {\`a} l'ordre $\ell$ le long de $W$ {\'e}quivaut {\`a} dire que l'application $$(z_{1},\ldots,z_{n})\mapsto F_{s,\sigma}(\mathsf{z}_{0}+\lambda(z_{1},\ldots,z_{n}),\mathsf{z}_{1}+z_{1},\ldots,\mathsf{z}_{n}+z_{n})$$s'annule {\`a} l'ordre $\ell$ au point $(0,\ldots,0)$.\par Passons maintenant au r{\'e}sultat principal de ce paragraphe.
\begin{prop}\label{proplemmedemultiplicites}Aucun {\'e}l{\'e}ment non nul de $E\otimes_{\sigma}\mathbf{C}_{p}$ ne s'annule le long de $W$ {\`a} l'ordre $(n+t)T$ en tous les points de l'ensemble $\{0_{G},\mathsf{p},\ldots,(n+t)S\mathsf{p}\}$.  
\end{prop}
La d{\'e}monstration est presque la m{\^e}me que celle de la proposition~$5.3$ de~\cite{artepredeux}, o{\`u} le r{\^o}le jou{\'e} par $\mathbf{G}_{\mathrm{m}}^{n}$ {\'e}tait d{\'e}volu {\`a} une vari{\'e}t{\'e} ab{\'e}lienne. Ici, il y a quelques simplifications qui augmentent la clart{\'e} de la preuve. C'est la raison pour laquelle nous la reproduisons ci-apr{\`e}s, sous une forme un peu abr{\'e}g{\'e}e.
\begin{proof}
Supposons cet {\'e}nonc{\'e} faux et qu'un tel polyn{\^{o}}me existe. D'apr{\`e}s le lemme de multiplicit{\'e}s~\cite{philippon2}, il existe un sous-groupe alg{\'e}brique connexe $G_{0}^{\star}$ (\emph{resp}. $\mathbf{G}_{\mathrm{m}}^{\star}$) de $G_{0,\mathbf{C}_{p}}$ (\emph{resp}. de $\mathbf{G}_{\mathrm{m},\mathbf{C}_{p}}^{n}$) tel que $G^{\star}:=G_{0}^{\star}\times\mathbf{G}_{\mathrm{m}}^{\star}\ne G_{\mathbf{C}_{p}}$ et \begin{equation}\label{ineqdeppph}
T^{\lambda^{\star}}\card\left(\frac{\Sigma_{\mathsf{p}}(S)+G^{\star}(\mathbf{C}_{p})}{G^{\star}(\mathbf{C}_{p})}\right)\EuScript{H}(G^{\star};D_{0}',\ldots,D_{n}')\le\EuScript{H}(G;D_{0}',\ldots,D_{n}')\end{equation}o{\`u} $D_{i}':=\max{\{1,D_{i}\}}$ pour tout $i\in\{1,\ldots,n\}$. Nous allons montrer que cette in{\'e}galit{\'e} ne peut pas {\^e}tre satisfaite en distinguant deux cas. Quitte {\`a} permuter les facteurs dans $\mathbf{G}_{\mathrm{m}}^{n}$, l'on peut supposer $D_{1}\le\cdots\le D_{n}$ sans perte de g{\'e}n{\'e}ralit{\'e}.\par\textbf{Premier cas}~: $t_{G^{\star}}+W\otimes_{\sigma}\mathbf{C}_{p}=t_{G_{\mathbf{C}_{p}}}$, c'est-{\`a}-dire $\lambda^{\star}=r^{\star}$. De l'in{\'e}galit{\'e}~\eqref{ineqdeppph} l'on d{\'e}duit l'existence la majoration \begin{equation*}\card\left(\frac{\Sigma_{\mathsf{p}}(S)+G^{\star}(\mathbf{C}_{p})}{G^{\star}(\mathbf{C}_{p})}\right)\le\frac{(n+t)!}{t!}\left(\frac{D_{0}}{T}\right)^{t-t^{\star}}\left(\frac{D_{n}}{T}\right)^{n-n^{\star}}\ \cdotp\end{equation*}En observant que le cardinal {\`a} gauche vaut $S+1$ lorsque $t=1$ et $t^{\star}=0$, la premi{\`e}re propri{\'e}t{\'e} de la proposition~\ref{propparametres} contredit cette in{\'e}galit{\'e}.\par\textbf{Second cas}~:  $t_{G^{\star}}+W\otimes_{\sigma}\mathbf{C}_{p}\ne t_{G_{\mathbf{C}_{p}}}$. En particulier on a $t^{\star}\ne t$. Soit $\kappa$ le plus petit entier de $\{1,\ldots,n\}$ tel que $D_{\kappa}\ge 1$. Soit $\pi_{\kappa}:\mathbf{G}_{\mathrm{m}}^{n}\to\mathbf{G}_{\mathrm{m}}^{n-\kappa+1}$ la projection sur les $n-\kappa+1$ derniers facteurs. Posons $G_{\kappa}^{\star}:=G_{0}^{\star}\times\mathbf{G}_{\mathrm{m}}^{\kappa-1}\times\pi_{\kappa}(\mathbf{G}_{\mathrm{m}}^{\star})$. Une propri{\'e}t{\'e} de la fonction $\EuScript{H}$ assure que $$\EuScript{H}(\pi_{\kappa}(\mathbf{G}_{\mathrm{m}}^{\star});D_{\kappa},\ldots,D_{n})\le\EuScript{H}(\mathbf{G}_{\mathrm{m}}^{\star};D_{1},\ldots,D_{n}),$$ce qui entra{\^\i}ne \begin{equation*}\frac{\EuScript{H}(G;D_{0}',\ldots,D_{n}')}{\EuScript{H}(G^{\star};D_{0}',\ldots,D_{n}')}\le n!\frac{\EuScript{H}(G;\widetilde{D}_{0},\ldots,\widetilde{D}_{n})}{\EuScript{H}(G_{\kappa}^{\star};\widetilde{D}_{0},\ldots,\widetilde{D}_{n})}x^{t-t^{\star}}.\end{equation*}De plus, comme $G^{\star}\subseteq G_{\kappa}^{\star}$ on a $\lambda_{\kappa}^{\star}:=\codim_{W}(W\cap t_{G_{\kappa}^{\star}})\le\lambda^{\star}$ et $$\card\left(\frac{\Sigma_{\mathsf{p}}(S)+G_{\kappa}^{\star}(\mathbf{C}_{p})}{G_{\kappa}^{\star}(\mathbf{C}_{p})}\right)\le\card\left(\frac{\Sigma_{\mathsf{p}}(S)+G^{\star}(\mathbf{C}_{p})}{G^{\star}(\mathbf{C}_{p})}\right)\cdotp$$De la sorte l'on obtient une in{\'e}galit{\'e} analogue {\`a}~\eqref{ineqdeppph} avec le sous-groupe alg{\'e}brique $G_{\kappa}^{\star}$~:\begin{equation}\label{ineqdeppphdeux}
T^{\lambda_{\kappa}^{\star}}\card\left(\frac{\Sigma_{\mathsf{p}}(S)+G_{\kappa}^{\star}(\mathbf{C}_{p})}{G_{\kappa}^{\star}(\mathbf{C}_{p})}\right)\EuScript{H}(G_{\kappa}^{\star};\widetilde{D}_{0},\ldots,\widetilde{D}_{n})\le n!\EuScript{H}(G;\widetilde{D}_{0},\ldots,\widetilde{D}_{n})x^{t-t^{\star}}.\end{equation}En reprenant l'{\'e}tude du premier cas (on a bien $G_{\kappa}^{\star}\ne G_{\mathbf{C}_{p}}$ car $t^{\star}\ne t$) et compte tenu du fait que $x\le 1$, on montre que $\lambda_{\kappa}^{\star}\ne r_{\kappa}^{\star}$. Ceci signifie que  $t_{G_{\kappa}^{\star}}+W\otimes_{\sigma}\mathbf{C}_{p}\ne t_{G_{\mathbf{C}_{p}}}$ et la majoration~\eqref{ineqdeppphdeux} entre alors en contradiction avec la d{\'e}finition de $x$.
\end{proof}
\subsection{Fibr{\'e} ad{\'e}lique hermitien des sections auxiliaires}\label{sectiontroissix} Reprenons l'espace vectoriel $E$ d{\'e}fini au d{\'e}but du \S~\ref{paragraphemultiplicite} et notons $\nu:=\dim E$ sa dimension. L'objectif de ce paragraphe est de munir $E$ d'une structure de fibr{\'e} ad{\'e}lique hermitien assez particuli{\`e}re de mani{\`e}re {\`a} permettre d'extrapoler plus tard sur les d{\'e}rivations. Tout d'abord, les espaces $E_{i}$ qui entrent dans la d{\'e}finition de $E$ sont naturellement munis d'une structure ad{\'e}lique hermitienne: pour $E_{0}=k\oplus k^{n}/W_{0}$ on choisit la structure quotient du fibr{\'e} $(k^{n},\vert\cdot\vert_{2})$ sur $k^{n}/W_{0}$ et l'on fait une somme directe hermitienne avec $k$ (muni de sa collection de valeurs absolues); pour $i\in\{1,\ldots,n\}$, on a $E_{i}=k\oplus k$ et l'on prend simplement la norme $\vert\cdot\vert_{2}$ de $k^{2}$.  Comme nous l'avons d{\'e}crit au \S~\ref{paragraphepuissance}, ces structures hermitiennes conf{\`e}rent {\`a} $E$ une premi{\`e}re structure de fibr{\'e} ad{\'e}lique hermitien pur $\overline{E}:=(E,(\Vert\cdot\Vert_{\sigma})_{\sigma})$. Mais ce n'est pas tout {\`a} fait cette structure que l'on va mettre sur $E$. Pour les normes relatives au plongement $\sigma_{0}$ et {\`a} ses conjugu{\'e}s, nous allons ajouter une quantit{\'e} suppl{\'e}mentaire, comme nous l'avons fait au \S~\ref{lemmedepetitesvaleurs} pour {\'e}tablir le lemme de Siegel approch{\'e} absolu.
\subsubsection*{Normes particuli{\`e}res}
Soit $(e_{1},\ldots,e_{n})$ la base canonique de $k^{n}$. Si $\sigma_{0}$ est ultram{\'e}trique, la famille $w:=\{w_{i}:=(\lambda(e_{i}),e_{i})\,;\ 1\le i\le n\}$ forme une base \emph{orthonorm{\'e}e} de $W\otimes_{\sigma_{0}}\mathbf{C}_{p_{0}}$. En effet, pour tout $\mathsf{z}=\sum_{i=1}^{n}{z_{i}e_{i}}\in\mathbf{C}_{p_{0}}^{n}$, on a $z:=(\lambda(\mathsf{z}),\mathsf{z})=\sum_{i=1}^{n}{z_{i}w_{i}}\in W\otimes_{\sigma_{0}}\mathbf{C}_{p_{0}}$ et $$\Vert z\Vert_{\overline{t_{G}},\sigma_{0}}=\max{(\Vert\lambda(\mathsf{z})\Vert_{\overline{k^{n}/W_{0}},\sigma_{0}},\vert\mathsf{z}\vert_{2,\sigma_{0}})}=\vert\mathsf{z}\vert_{2,\sigma_{0}}$$par propri{\'e}t{\'e} de la norme quotient. Dans le cas d'un plongement $\sigma_{0}$ archim{\'e}dien, on note $(w_{1},\ldots,w_{n})$ la base orthonorm{\'e}e de $W\otimes_{\sigma_{0}}\mathbf{C}_{p_{0}}$ introduite au \S~\ref{subseccomplements}. La famille $\{(\lambda(e_{i}),e_{i});\ 1\le i\le n\}$ est une base de $W$, dans laquelle les vecteurs $w_{i}$ peuvent s'{\'e}crire~: il existe des formes lin{\'e}aires $\{\mathsf{l}_{i}(z_{1},\ldots,z_{n});\ 1\le i\le n\}\subseteq(\mathbf{C}_{p_{0}}^{n})^{\mathsf{v}}$ telles que \begin{equation}\label{defiformeslineaires}\forall\,(z_{1},\ldots,z_{n})\in\mathbf{C}_{p_{0}}^{n},\quad\sum_{i=1}^{n}{z_{i}w_{i}}=\sum_{i=1}^{n}{\mathsf{l}_{i}(z_{1},\ldots,z_{n})(\lambda(e_{i}),e_{i})}.\end{equation}Les coefficients des $\mathsf{l}_{i}$ sont dans le compl{\'e}t{\'e} $k_{\sigma_{0}}$ de $\sigma_{0}(k)$ dans $\mathbf{C}_{p_{0}}$~: c'est {\'e}vident si $\sigma_{0}$ est complexe ou si $\sigma_{0}$ est ultram{\'e}trique car $\mathsf{l}_{i}(z_{1},\ldots,z_{n})=z_{i}$, et c'est aussi vrai si $\sigma_{0}$ est un plongement r{\'e}el car dans ce cas $w_{i}$ a {\'e}t{\'e} choisi dans $W\otimes_{k}k_{\sigma_{0}}$. Pour  $s\in E\otimes_{\sigma_{0}}\mathbf{C}_{p_{0}}$, $m\in\mathbf{N}$ et $\tau\in\mathbf{N}^{n}$, notons $\frac{1}{\tau !}\mathrm{D}_{w}^{\tau}F_{s,\sigma_{0}}(m(u_{0},u))$ le $\tau^{\text{{\`e}me}}$ coefficient de Taylor {\`a} l'origine de l'application $$z=(z_{1},\ldots,z_{n})\mapsto F_{s,\sigma_{0}}\left(mu_{0}+\lambda(\mathsf{l}_{1}(z),\ldots,\mathsf{l}_{n}(z)),mu_{1}+\mathsf{l}_{1}(z),\ldots,mu_{n}+\mathsf{l}_{n}(z)\right)$$(voir~\eqref{defidefs} pour la d{\'e}finition de $F_{s,\sigma_{0}}$). Soit $(s_{1},\ldots,s_{\nu})$ une base orthonorm{\'e}e de $(E\otimes_{k}k_{\sigma_{0}},\Vert\cdot\Vert_{\sigma_{0}})$. Pour $z=(z_{0},z_{1},\ldots,z_{n})\in t_{G_{0}}(\mathbf{C}_{p_{0}})\oplus\mathcal{T}_{p_{0}}^{n}$, d{\'e}finissons  \begin{equation}\label{defidealpha}\alpha_{\sigma_{0}}(z):=\prod_{j=1}^{n}{\vert(1,e^{z_{j}})\vert_{2,\sigma_{0}}^{-D_{j}}}\times\moy_{\sigma_{0}}(1,\Vert z_{0}\Vert_{\overline{t_{G}},\sigma_{0}})^{-D_{0}}.\end{equation}Il existe $\alpha_{\sigma_{0}}'(m(u_{0},u))\in k_{\sigma_{0}}$ de valeur absolue $p_{0}$-adique {\'e}gale {\`a} $\alpha_{\sigma_{0}}(m(u_{0},u))$. C'est clair si $\sigma_{0}$ est complexe. Si $\sigma_{0}$ est ultram{\'e}trique, le vecteur $u_{0}\in k^{n}/W_{0}$ a des coordonn{\'e}es dans $k_{\sigma_{0}}$ par rapport {\`a} une $k_{\sigma_{0}}$-base orthonorm{\'e}e de $(k^{n}/W_{0})\otimes_{\sigma_{0}}k_{\sigma_{0}}$. Ainsi $\Vert u_{0}\Vert_{\overline{t_{G}},\sigma_{0}}$, qui n'est rien d'autre que le maximum des valeurs absolues de ces coordonn{\'e}es, est un {\'e}l{\'e}ment de $\vert k_{\sigma_{0}}\vert_{p_{0}}$. C'est aussi le cas pour $\vert(1,\alpha_{j}^{m})\vert_{2,\sigma_{0}}$ qui vaut $1=\vert 1\vert_{p_{0}}$ ou $\vert\alpha_{j}^{m}\vert_{p_{0}}$. On a donc $\alpha_{\sigma_{0}}(m(u_{0},u))\in\vert k_{\sigma_{0}}\vert_{p_{0}}$ et l'existence de $\alpha_{\sigma_{0}}'(m(u_{0},u))$ en d{\'e}coule. Consid{\'e}rons alors la matrice $\mathsf{A}_{\sigma_{0}}$, de taille $\mu\times\nu$, dont les coefficients sont les {\'e}l{\'e}ments de $k_{\sigma_{0}}$ suivants~: pour tout $(m,\tau)\in\Upsilon$, pour tout $i\in\{1,\ldots,\nu\}$, \begin{equation}\label{defidescoefficientsdeA}\mathsf{A}_{\sigma_{0}}[(m,\tau),i]:=\left(\frac{1}{\tau !}\mathrm{D}_{w}^{\tau}F_{s_{i},\sigma_{0}}(m(u_{0},u))\right)\times\alpha_{\sigma_{0}}'(m(u_{0},u)).\end{equation}Nous sommes maintenant en mesure de d{\'e}finir la norme voulue sur $E\otimes_{\sigma_{0}}\mathbf{C}_{p_{0}}$~: Soit $\alpha>0$ (choisi dans la proposition~\ref{propconstructiondes}). Pour tout $s=\sum_{i=1}^{\nu}{x_{i}s_{i}}\in E\otimes_{\sigma_{0}}\mathbf{C}_{p_{0}}$, \begin{equation}\label{defidesnormess}\Vert s\Vert_{\overline{E}_{\alpha},\sigma_{0}}:=\moy_{\sigma_{0}}(\vert x\vert_{2,\sigma_{0}},\alpha\vert\mathsf{A}_{\sigma_{0}}x\vert_{2,\sigma_{0}}).\end{equation}Dans ces expressions, $x$ d{\'e}signe le vecteur colonne de coordonn{\'e}es $x_{1},\ldots,x_{\nu}$. Pour un plongement $\sigma$ conjugu{\'e} {\`a} $\sigma_{0}$, c'est-{\`a}-dire d{\'e}finissant la m{\^{e}}me place $v_{0}$, nous choisissons la norme $\Vert\cdot\Vert_{\overline{E}_{\alpha},\sigma}$ de sorte que la collection de normes $(\Vert\cdot\Vert_{\overline{E}_{\alpha},\sigma})_{\sigma\in v_{0}}$ soit invariante par Galois (au sens de la d{\'e}finition~\ref{defiinvariance}). En {\'e}crivant $\sigma=\iota\circ\sigma_{0}$ avec $\iota$ automorphisme continu de $\mathbf{C}_{p_{0}}$, on observe que $\Vert\cdot\Vert_{\overline{E}_{\alpha},\sigma}$ v{\'e}rifie une {\'e}galit{\'e} du type~\eqref{defidesnormess} en prenant $\mathsf{A}_{\sigma}=\iota(\mathsf{A}_{\sigma_{0}})$ (les bases $(s_{1},\ldots,s_{\nu})$ et $w$ {\'e}tant remplac{\'e}es par leurs images par $\iota$). Par d{\'e}finition, le fibr{\'e} ad{\'e}lique hermitien $\overline{E}_{\alpha}$ est le fibr{\'e} {\'e}gal {\`a} $\overline{E}$ sauf pour les plongement $\sigma\in v_{0}$ o\`u les normes sont donn{\'e}es par $(\Vert\cdot\Vert_{\overline{E}_{\alpha},\sigma})_{\sigma\in v_{0}}$ (on retrouve $\overline{E}$ en prenant $\alpha=0$).
\begin{lemm}\label{lemmecinq}La pente d'Arakelov normalis{\'e}e $\widehat{\mu}(\overline{E})$ du fibr{\'e} ad{\'e}lique hermitien $\overline{E}$ est minor{\'e}e par $-D_{0}h(\overline{W_{0}})/(t+1)$.
\end{lemm}En r{\'e}alit{\'e}, on conna{\^i}t une formule exacte pour cette pente, formule qui entra{\^\i}ne imm{\'e}diatement le lemme. Pour cela, l'on peut se r{\'e}f{\'e}rer au lemme~$7.3.$ de~\cite{rendiconti} en tenant compte du fait que la pente d'Arakelov est additive vis {\`a} vis du produit tensoriel de fibr{\'e}s ad{\'e}liques hermitiens (voir~\cite[proposition~$5.2$]{rendiconti}). 
\subsection{Estimation du rang d'un syst{\`e}me lin{\'e}aire}
\begin{lemm}\label{lemmefacile}Soit $a,x,y$ des entiers naturels. Alors on a $$\binom{x+y+a}{a}-\binom{x+a}{a}\le ay(x+y+1)^{a-1}.$$
\end{lemm}
\begin{proof}
{\`A} $x,a$ fix{\'e}s, la fonction $H:y\mapsto\binom{x+y+a}{a}-\binom{x+a}{a}$ est une fonction polynomiale {\`a} coefficients positifs. Cette fonction est convexe sur $\mathbf{R}^{+}$ car de d{\'e}riv{\'e}e seconde positive. On a donc $H(y)=H(y)-H(0)\le yH'(y)$. Le calcul explicite de $H'(y)$ comme une somme de $a$ produits de $a-1$ termes conduit {\`a} la majoration $H'(y)\le\sum_{i=1}^{a}{(x+y+1)^{a-1}}=a(x+y+1)^{a-1}$.\end{proof}
\begin{prop}\label{proprang} Soit $\rho$ le rang de la matrice $\mathsf{A}_{\sigma_{0}}$ de coefficients~\eqref{defidescoefficientsdeA}. Alors on a $\rho\le2(4n)^{3n}\nu/C_{0}$. 
\end{prop}
\begin{proof}
Consubstantielle {\`a} la distinction entre le cas p{\'e}riodique et le cas non p{\'e}riodique et {\`a} l'introduction du groupe $\widetilde{G}$, l'argumentation n'a pas chang{\'e} depuis l'article fondateur de Philippon \& Waldschmidt~\cite{pph-miw}. En suivant l'approche plus effective de~\cite[\S~$6.3$]{sinnoudeux} et  en posant $$(T_{1},S_{1}):=\begin{cases}(T,S_{0}) & \text{dans le cas non p{\'e}riodique},\\ (T_{0},S) & \text{dans le cas p{\'e}riodique},\end{cases}$$on montre que \begin{equation}\label{ineq:majorationderho}\rho\le X\card\left(\frac{\Sigma_{\mathsf{p}}(S_{1})+\widetilde{G}(\overline{k})}{\widetilde{G}(\overline{k})}\right)(\dim\widetilde{G}+1+\EuScript{H}(\widetilde{G};D_{0}',\ldots,D_{n}'))\end{equation}avec \begin{equation*}X=\begin{cases}\binom{2(n+t)T+\widetilde{\lambda}}{\widetilde{\lambda}} & \text{dans le cas non p{\'e}riodique},\\ \binom{2(n+t)T+\widetilde{\lambda}}{\widetilde{\lambda}} -\binom{2(n+t)T-T_{0}+\widetilde{\lambda}}{\widetilde{\lambda}}& \text{dans le cas p{\'e}riodique}.\end{cases}\end{equation*}Le lemme~\ref{lemmefacile} entra{\^{\i}}ne la majoration $X\le (2(n+t))^{n}T^{\widetilde{\lambda}-1}T_{1}$. Par ailleurs, rappelons que, pour tout $i\in\{0,\ldots,n\}$, la fonction partielle $x_{i}\mapsto\frac{\EuScript{H}(\widetilde{G};x_{0},\ldots,x_{n})}{\EuScript{H}(G;x_{0},\ldots,x_{n})}$ d{\'e}cro{\^{\i}}t sur $]0,+\infty[$. Cette observation et les majorations $\widetilde{D}_{i}/2\le D_{i}'$ pour $i\in\{1,\ldots,n\}$ et $x\widetilde{D}_{0}/2\le D_{0}$ impliquent
  \begin{equation*}\begin{split}\frac{\EuScript{H}(\widetilde{G};D_{0}',\ldots,D_{n}')}{\EuScript{H}(G;D_{0}',\ldots,D_{n}')} &\le 2^{n+t}\frac{\EuScript{H}(\widetilde{G};x\widetilde{D}_{0},\widetilde{D}_{1},\ldots,\widetilde{D}_{n})}{\EuScript{H}(G;x\widetilde{D}_{0},\ldots,\widetilde{D}_{n})}\\ &=\frac{2^{n+t}C_{0}}{\widetilde{T}^{\widetilde{\lambda}}\card\left(\frac{\Sigma_{\mathsf{p}}(S)+\widetilde{G}(\overline{k})}{\widetilde{G}(\overline{k})}\right)}\end{split}\end{equation*}(l'{\'e}galit{\'e} repose sur la d{\'e}finition de $x=x(\widetilde{G})$). En injectant cette estimation dans~\eqref{ineq:majorationderho}, nous avons \begin{equation*}\rho\le 2^{n+t}(2(n+t))^{n+1}C_{0}\max{\left\{\frac{T_{0}}{\widetilde{T}},\frac{S_{0}}{S}\right\}}\EuScript{H}(G;D_{0}',\ldots,D_{n}').\end{equation*}La dimension $\nu$ de $E$ vaut $\binom{D_{0}+t}{t}(D_{1}+1)\cdots(D_{n}+1)$ et, en particulier, on a \begin{equation*}\EuScript{H}(G;D_{0}',\ldots,D_{n}')=\frac{(n+t)!}{t!}D_{0}^{t}D_{1}'\cdots D_{n}'\le(n+t)!\nu.\end{equation*}Par cons{\'e}quent, en majorant $(n+t)!$ par $(n+t)^{n+t-1}$, on obtient \begin{equation*}\rho\le 2(2(n+t))^{2n+t}C_{0}\max{\left\{\frac{T_{0}}{\widetilde{T}},\frac{S_{0}}{S}\right\}}\nu.\end{equation*}Le choix des param{\`e}tres et l'in{\'e}galit{\'e} $t\le n$  permettent alors de conclure.
\end{proof}
\subsection{Estimation d'une d{\'e}riv{\'e}e}Soit $\sigma:k\hookrightarrow\mathbf{C}_{p}$ un plongement de $k$. Soit $s\in E\otimes_{\sigma}\mathbf{C}_{p}$ que l'on {\'e}crit dans une base orthonorm{\'e}e comme dans~\eqref{ecrituredes} avec des coefficients $p_{i}$. On notera $L(s):=\sum_{i}{\vert p_{i}\vert_{p}}$ si $\sigma$ est archim{\'e}dien et $L(s):=\max_{i}{\vert p_{i}\vert_{p}}$ si $\sigma$ est ultram{\'e}trique. Soit $\mathsf{w}=(\mathsf{w}_{1},\ldots,\mathsf{w}_{n})$ une base de $W\otimes_{\sigma}\mathbf{C}_{p}$. En proc{\'e}dant de la m{\^{e}}me mani{\`e}re qu'au paragraphe~\ref{sectiontroissix} (en rempla{\c{c}}ant la base $w$ d'alors par, ici, $\mathsf{w}$), on dispose de la d{\'e}riv{\'e}e divis{\'e}e $\frac{1}{\tau !}\mathrm{D}_{\mathsf{w}}^{\tau}F_{s,\sigma}(z)$ en un point $z=(z_{0},\ldots,z_{n})\in\mathbf{C}_{p}^{t}\times\mathcal{T}_{p}^{n}$ {\`a} l'ordre $\tau=(\tau_{1},\ldots,\tau_{n})\in\mathbf{N}^{n}$. Voici une premi{\`e}re estimation, {\'e}l{\'e}mentaire malgr{\'e} son aspect technique, d'un usage fr{\'e}quent dans la suite.
\begin{prop}\label{propmajorationarchimediennesimple} Dans les conditions ci-dessus ($\sigma,s,\mathsf{w},z,\tau$ quelconques), on a \begin{equation*}\left\vert\frac{1}{\tau !}\mathrm{D}_{\mathsf{w}}^{\tau}F_{s,\sigma}(z)\right\vert_{p}\cdot\alpha_{\sigma}(z)\le\left\{\prod_{i=1}^{n}{\Vert\mathsf{w}_{i}\Vert_{\overline{t_{G}},\sigma}^{\tau_{i}}}\right\}L(s)\times\begin{cases} e^{2n\sqrt{n}\max_{0\le i\le n}{\{D_{i}\}}} & \text{si $p=\infty$},\\ e^{\vert\tau\vert} & \text{si $p\ne\infty$}\end{cases}
\end{equation*}o{\`u}\begin{equation*}\alpha_{\sigma}(z):=\prod_{j=1}^{n}{\vert(1,e^{z_{j}})\vert_{2,\sigma}^{-D_{j}}}\times\moy_{\sigma}(1,\Vert z_{0}\Vert_{\overline{t_{G}},\sigma})^{-D_{0}}.\end{equation*}
\end{prop}
\begin{proof}
Tout d'abord, observons que si l'on multiplie le vecteur $\mathsf{w}_{i}$ par un nombre $\theta_{i}\in\mathbf{C}_{p}$ alors $\frac{1}{\tau !}\mathrm{D}_{\mathsf{w}}^{\tau}F_{s,\sigma}(z)$ est multipli{\'e} par $\prod_{i=1}^{n}{\theta_{i}^{\tau_{i}}}$. On peut donc supposer que, pour tout $i\in\{1,\ldots,n\}$, $\Vert\mathsf{w}_{i}\Vert_{\overline{t_{G}},\sigma}=1$. Consid{\'e}rons maintenant des formes lin{\'e}aires $\mathsf{l}_{1},\ldots,\mathsf{l}_{n}\in(\mathbf{C}_{p}^{n})^{\mathsf{v}}$ telles que, pour tout $x=(x_{1},\ldots,x_{n})\in\mathbf{C}_{p}^{n}$, on a $\sum_{i=1}^{n}{x_{i}\mathsf{w}_{i}}=\sum_{i=1}^{n}{\mathsf{l}_{i}(x)(\lambda(e_{i}),e_{i})}$ (rappelons que $e_{i}$ d{\'e}signe le $i^{\text{{\`e}me}}$ vecteur de la base canonique de $\mathbf{C}_{p}^{n}$). La somme $\sum_{i=1}^{n}{\mathsf{l}_{i}(e_{j})\lambda(e_{i})}$ est la composante de $\mathsf{w}_{j}$ sur $k^{n}/W_{0}\otimes_{\sigma}\mathbf{C}_{p}$. Elle est donc de norme plus petite que $1$.\par Soit $(e_{0,1},\ldots,e_{0,t})$ la base orthonorm{\'e}e de $k^{n}/W_{0}$ choisie pour {\'e}crire les coefficients de $s$ (cf.~\eqref{ecrituredes}). Le terme $\frac{1}{\tau !}\mathrm{D}_{\mathsf{w}}^{\tau}F_{s,\sigma}(z)$ que l'on cherche {\`a} {\'e}valuer est le coefficient devant $x^{\tau}$ de l'application \begin{equation*}x=(x_{1},\ldots,x_{n})\mapsto s((1,z_{0}+\lambda(\mathsf{l}_{1}(x),\ldots,\mathsf{l}_{n}(x))),(1,e^{z_{1}+\mathsf{l}_{1}(x)}),\ldots,(1,e^{z_{n}+\mathsf{l}_{n}(x)})).\end{equation*}Celle-ci est une somme de termes de la forme \begin{equation}\label{eq:coeffsdl}p_{\mathsf{h}}\left\{\prod_{a=1}^{t}{\left(e_{0,a}^{\mathsf{v}}(z_{0})+\sum_{i=1}^{n}{\mathsf{l}_{i}(x)e_{0,a}^{\mathsf{v}}(\lambda(e_{i}))}\right)^{\mathsf{h}_{0,a}}}\right\}\prod_{i=1}^{n}{e^{\mathsf{h}_{i}z_{i}+\mathsf{h}_{i}\mathsf{l}_{i}(x)}}\end{equation}o\`u $\mathsf{h}_{0}=(\mathsf{h}_{0,1},\ldots,\mathsf{h}_{0,t})\in\mathbf{N}^{t}$, $\mathsf{h}=(\mathsf{h}_{0},\mathsf{h}_{1},\ldots,\mathsf{h}_{n})\in\mathbf{N}^{t}\times\mathbf{N}^{n}$ avec, pour tout $i\in\{0,\ldots,n\}$, $\vert\mathsf{h}_{i}\vert\le D_{i}$, et $p_{\mathsf{h}}\in\mathbf{C}_{p}$ est l'un des coefficients de $s$. Comme seul le coefficient devant $x^{\tau}$ nous int{\'e}resse, on peut remplacer l'exponentielle $\exp{\{\sum_{i=1}^{n}{\mathsf{h}_{i}\mathsf{l}_{i}(x)}\}}$ qui est dans~\eqref{eq:coeffsdl} par son d{\'e}veloppement de Taylor {\`a} l'ordre $\vert\tau\vert$~: \begin{equation}\label{definitiondeqhtau}Q_{\mathsf{h},\tau}(x):=\sum_{j=0}^{\vert\tau\vert}{\frac{1}{j!}(\mathsf{h}_{1}\mathsf{l}_{1}(x)+\cdots+\mathsf{h}_{n}\mathsf{l}_{n}(x))^{j}}.\end{equation}De la sorte la quantit{\'e}~\eqref{eq:coeffsdl} peut {\^{e}}tre {\'e}chang{\'e}e par l'expression polynomiale \begin{equation}\label{expressionvingt}p_{\mathsf{h}}\left\{\prod_{a=1}^{t}{\left(e_{0,a}^{\mathsf{v}}(z_{0})+\sum_{j=1}^{n}{x_{j}e_{0,a}^{\mathsf{v}}\left(\sum_{i=1}^{n}{\mathsf{l}_{i}(e_{j})\lambda(e_{i})}\right)}\right)^{\mathsf{h}_{0,a}}}\right\}\left(\prod_{i=1}^{n}{e^{\mathsf{h}_{i}z_{i}}}\right)Q_{\mathsf{h},\tau}(x).\end{equation}La fonction longueur $\mathsf{P}\mapsto L(\mathsf{P})$ est sous-multiplicative~: $L(\mathsf{P}\mathsf{Q})\le L(\mathsf{P})L(\mathsf{Q})$ (voir~\cite[p.~$76$]{miw4}). Cette propri{\'e}t{\'e} appliqu{\'e}e {\`a} la somme sur $\mathsf{h}$ des expressions~\eqref{expressionvingt} permet de majorer $\left\vert\frac{1}{\tau !}\mathrm{D}_{\mathsf{w}}^{\tau}F_{s,\sigma}(z)\right\vert_{p}$ par \begin{equation*}L(s)(\Vert z_{0}\Vert_{\overline{t_{G}},\sigma}+n)^{D_{0}}\left\{\prod_{i=1}^{n}{\max{(1,\vert e^{z_{i}}\vert_{p})^{D_{i}}}}\right\}\max_{\mathsf{h}}{\{L(Q_{\mathsf{h},\tau})\}}\end{equation*} (on notera que le vecteur $e_{0,a}^{\mathsf{v}}$ est de norme $1$ et donc $\vert e_{0,a}^{\mathsf{v}}(\mathsf{z})\vert_{p}\le\Vert\mathsf{z}\Vert_{\overline{t_{G}},\sigma}$ pour tout $\mathsf{z}\in k^{n}/W_{0}\otimes_{\sigma}\mathbf{C}_{p}$). Pour conclure, il suffit d'estimer le maximum qui appara{\^{\i}}t {\`a} la fin de cette expression. Lorsque le plongement $\sigma$ est archim{\'e}dien, la longueur de $Q_{\mathsf{h},\tau}$ est major{\'e}e par
$$\sum_{j=0}^{\vert\tau\vert}{\frac{1}{j!}\left(\sum_{i=1}^{n}{D_{i}\vert\mathsf{l}_{i}(1,\ldots,1)\vert}\right)^{j}}.$$Or $\sum_{i=1}^{n}{\mathsf{l}_{i}(1,\ldots,1)e_{i}}$ est la composante sur $\mathbf{C}_{p}^{n}$ de $\mathsf{w}_{1}+\cdots+\mathsf{w}_{n}$ et donc $\sum_{i=1}^{n}{\vert\mathsf{l}_{i}(1,\ldots,1)\vert^{2}}\le n^{2}$. Ainsi on a $\sum_{i=1}^{n}{D_{i}\vert\mathsf{l}_{i}(1,\ldots,1)\vert}\le\sqrt{D_{1}^{2}+\cdots+D_{n}^{2}}\times n\le n\sqrt{n}\max{\{D_{1},\ldots,D_{n}\}}$. On en d{\'e}duit dans ce cas \begin{equation}\label{eqmajorantdeq}L(Q_{\mathsf{h},\tau})\le e^{n\sqrt{n}\max{\{D_{1},\ldots,D_{n}\}}}.\end{equation}Lorsque $\sigma$ est ultram{\'e}trique alors, pour tout $j\in\mathbf{N}$, on a $\vert j!\vert_{p}^{-1}\le p^{j/(p-1)}\le e^{j}$. Dans ce cas la longueur de $Q_{\mathsf{h},\tau}$ est plus petite que $e^{\vert\tau\vert}$ et la proposition~\ref{propmajorationarchimediennesimple} est d{\'e}montr{\'e}e.
\end{proof}
\begin{coro}\label{coro:evaluationcoeffsa}Pour tout $\sigma\in v_{0}$, la norme d'op{\'e}rateur de la matrice $\mathsf{A}_{\sigma}$, matrice d{\'e}finie au \S~\ref{sectiontroissix}, est inf{\'e}rieure {\`a} $e^{4nU_{0}}$.
\end{coro}
\begin{proof}Il suffit de le montrer pour $\sigma=\sigma_{0}$. On applique la proposition~\ref{propmajorationarchimediennesimple} {\`a} chacune des fonctions $F_{s_{i},\sigma_{0}}$ qui intervient dans la d{\'e}finition~\eqref{defidescoefficientsdeA} des coefficients de $\mathsf{A}_{\sigma_{0}}$ avec $\sigma=\sigma_{0}$, $\mathsf{w}=w$ et $z=m(u_{0},u)$. Supposons dans un premier temps que $\sigma_{0}$ est archim{\'e}dien. D'apr{\`e}s la fin du paragraphe~\ref{paragraphepuissance}, la longueur de $s_{i}$ est plus petite que \begin{equation}\label{longueurs}\sqrt{t}^{D_{0}}\sqrt{2}^{D_{1}+\cdots+D_{n}}\le(n2^{n})^{\max{\{D_{0},\ldots,D_{n}\}}/2}\le e^{n\max{\{D_{0},\ldots,D_{n}\}}}.\end{equation}Ainsi les coefficients de $\mathsf{A}_{\sigma_{0}}$ sont inf{\'e}rieurs {\`a} $e^{3n\sqrt{n}\max_{0\le i\le n}{\{D_{i}\}}}$, quantit{\'e} elle-m{\^{e}}me plus petite que $e^{nU_{0}}$ d'apr{\`e}s le choix des param{\`e}tres (d{\'e}but du \S~\ref{section:choixdesparametres}) et car $\log C_{0}\ge 3n$. En la comparant {\`a} la norme de Hilbert-Schmidt, on constate que la norme d'op{\'e}rateur de $\mathsf{A}_{\sigma_{0}}$ est plus petite que $(\nu\mu)^{1/2}e^{nU_{0}}$. La dimension $\nu$ de $E$ vaut $\binom{D_{0}+t}{t}(D_{1}+1)\cdots(D_{n}+1)$ et donc \begin{equation}\label{ineqsurnu}\nu\le(D_{0}+1)^{t}(D_{1}+1)\cdots(D_{n}+1)\le e^{tD_{0}+D_{1}+\cdots+D_{n}}\le e^{2n\max{\{D_{0},\ldots,D_{n}\}}}\le e^{U_{0}/D}\le e^{U_{0}}\end{equation}car le maximum des $\widetilde{D}_{i}$ est plus petit que $U_{0}/(D\log C_{0})$. De m{\^{e}}me, $\mu$ est le cardinal de $\Upsilon$ (introduit {\`a} la suite de la d{\'e}finition~\ref{defi:casperiodique}). Il est inf{\'e}rieur {\`a} \begin{align*}\binom{2(n+t)T+n}{n}S(n+t)&\le\left(2(n+t)T+1\right)^{n}S(n+t)\\ &\le(4n+1)^{n}(2n)T^{n}S\\ &\le\frac{(4n+1)^{n}(2n)(n+1)!}{C_{0}}\frac{U_{0}^{n+1}}{(n+1)!}\quad\text{car $S\le U_{0}/C_{0}$}\\ &\le1\times e^{U_{0}}=e^{U_{0}}.\end{align*}Ainsi on trouve $\Vert\mathsf{A}_{\sigma_{0}}\Vert_{\sigma_{0}}\le e^{(n+1)U_{0}}\le e^{4nU_{0}}$ lorsque $\sigma_{0}$ est archim{\'e}dien. Dans le cas ultram{\'e}trique, la norme d'op{\'e}rateur de $\mathsf{A}_{\sigma_{0}}$ est le maximum des coefficients de $\mathsf{A}_{\sigma_{0}}$ (voir proposition~\ref{propnop}). D'apr{\`e}s la proposition~\ref{propmajorationarchimediennesimple}, on a $\Vert\mathsf{A}_{\sigma_{0}}\Vert_{\sigma_{0}}\le e^{2(n+t)T}\le e^{4nU_{0}}$ car chacun des $s_{i}$ est de longueur $1$.
\end{proof}
\subsection{Construction d'une section auxiliaire}\label{secconstructiondes}
\begin{prop}\label{propconstructiondes}Soit $\alpha>0$ tel que $\log\alpha=C_{0}^{3/2}U_{0}$. Il existe une section $s\in E\otimes_{\sigma_{0}}\overline{k}$, non nulle, telle que \begin{equation}\label{ineq:sbati}h_{\overline{E}_{\alpha}}(s)\le\frac{[k_{\sigma_{0}}:\mathbf{Q}_{p_{0}}]\times(4n)^{4n}C_{0}^{1/2}U_{0}}{D}.\end{equation}
\end{prop}
\begin{proof}
D'apr{\`e}s le lemme de Siegel approch{\'e} absolu du~\S~\ref{lemmedepetitesvaleurs}, il existe un vecteur $s\in E\otimes_{\sigma_{0}}\overline{k}$, non nul, tel que 
\begin{equation}\label{majorationhauteurdes}h_{\overline{E}_{\alpha}}(s)\le\frac{\rho}{\nu D}\left(\sum_{\sigma\in v_{0}}{\log\moy_{\sigma}(1,\alpha)+\log\max{\{1,\Vert\mathsf{A}_{\sigma}\Vert_{\sigma}\}}}\right)+\frac{1}{2}\log\nu-\widehat{\mu}(\overline{E})\end{equation}(on notera que le rang $\rho$ de $\mathsf{A}_{\sigma_{0}}$ est aussi celui de $\mathsf{A}_{\sigma}$). Dans le membre de droite, on majore $\rho/\nu$ avec la proposition~\ref{proprang} et $\Vert\mathsf{A}_{\sigma}\Vert_{\sigma}$ est {\'e}valu{\'e} au moyen du corollaire~\ref{coro:evaluationcoeffsa} par $e^{4nU_{0}}$. De plus, gr{\^{a}}ce au lemme~\ref{lemmecinq}, {\`a} la proposition~\ref{propositionhauteurwo}, au choix des param{\`e}tres $D_{0}$ et $x\le 1$, on a $$-\widehat{\mu}(\overline{E})\le \frac{D_{0}h(\overline{W_{0}})}{t+1}\le\frac{n^{3}U_{0}}{D}.$$Enfin, on a $\frac{1}{2}\log\nu\le U_{0}/(2D)$ d'apr{\`e}s~\eqref{ineqsurnu}. En rempla{\c c}ant ces estimations dans~\eqref{majorationhauteurdes}, on trouve \begin{equation*}h_{\overline{E}_{\alpha}}(s)\le\frac{2(4n)^{3n}[k_{\sigma_{0}}:\mathbf{Q}_{p_{0}}]}{DC_{0}}\times\left(\log\sqrt{2}+C_{0}^{3/2}U_{0}+4nU_{0}\right)+\frac{U_{0}}{2D}+\frac{n^{3}U_{0}}{D}\end{equation*}et l'on conclut en utilisant la valeur de $C_{0}=(4n)^{10n}$.\end{proof}
La d{\'e}monstration du th{\'e}or{\`e}me~\ref{theoremereduit} s'effectue avec la section $s$ que l'on vient de construire dans cette proposition, qui v{\'e}rifie~\eqref{ineq:sbati}. \textit{A priori} cette section n'est pas d{\'e}finie sur $k$ mais sur une extension finie $K$ de $k$. Cette complication technique n'a pas de cons{\'e}quence. En effet, le degr{\'e} relatif $[K:k]$ n'intervient pas car les estimations des jets de $s$ en les plongements $\sigma'$ de $K$ qui prolonge $\sigma:k\hookrightarrow\mathbf{C}_{p}$ sont de la forme $\mathfrak{c}_{\sigma}\Vert s\Vert_{\overline{E}_{\alpha},\sigma'}$, o{\`u} $\mathfrak{c}_{\sigma}$ ne d{\'e}pend que de $\sigma$ (et des autres donn{\'e}es) et pas de $\sigma'$. C'est la raison pour laquelle nous \emph{supposerons} --- et ceci sans perte de g{\'e}n{\'e}ralit{\'e} --- que $s$ est d{\'e}finie sur $k$. 
\subsection{Estimations g{\'e}n{\'e}rales}
Soit $(m,\ell)\in\mathbf{N}^{2}$. Consid{\'e}rons la section $s$ construite dans la proposition~\ref{propconstructiondes}. Soit $\mathsf{w}=(\mathsf{w}_{1},\ldots,\mathsf{w}_{n})$ une $k$-base de $W$. Il existe des formes lin{\'e}aires $\mathsf{l}_{1},\ldots,\mathsf{l}_{n}$ sur $k^{n}$ telles que, si $x=(x_{1},\ldots,x_{n})$ est un $n$-uplet de variables, on ait $$\sum_{i=1}^{n}{x_{i}\mathsf{w}_{i}}=(\lambda(\mathsf{l}_{1}(x),\ldots,\mathsf{l}_{n}(x)),\mathsf{l}_{1}(x),\ldots,\mathsf{l}_{n}(x)).$$Tous les coefficients de Taylor {\`a} l'origine --- $a(s,\mathsf{w},m,\tau)$, $\tau\in\mathbf{N}^{n}$ --- de la s{\'e}rie formelle \begin{equation}\label{fonctionfs}x=(x_{1},\ldots,x_{n})\mapsto s\left((1,mu_{0}+\lambda(\mathsf{l}_{1}(x),\ldots,\mathsf{l}_{n}(x))),(1,\alpha_{1}^{m}e^{\mathsf{l}_{1}(x)}),\ldots,(1,\alpha_{n}^{m}e^{\mathsf{l}_{n}(x)})\right)\end{equation}sont des {\'e}l{\'e}ments de $k$. {\'E}tant donn{\'e} un plongement $\sigma:k\hookrightarrow\mathbf{C}_{p}$ pour laquelle $m\mathsf{p}$ poss{\`e}de un logarithme $z\in t_{G}(\mathbf{C}_{p})$, l'image du coefficient $a(s,\mathsf{w},m,\tau)$ dans $\sigma(k)$ est {\'e}gale {\`a} $\frac{1}{\tau !}\mathrm{D}_{\mathsf{w}}^{\tau}F_{s,\sigma}(z)$ d{\'e}finie au \S~\ref{sectiontroissix}.
\begin{defi}\label{definitionjet}Soit $(\mathsf{w}_{1}^{\mathsf{v}},\ldots,\mathsf{w}_{n}^{\mathsf{v}})$ la base duale de $\mathsf{w}$. Le \emph{jet de $s$ d'ordre $\ell$ le long de $W$ au point $m\mathsf{p}$}, not{\'e} $\jet_{W}^{\ell}s(m\mathsf{p})$, est le vecteur \begin{equation}\label{defidujet}\jet_{W}^{\ell}s(m\mathsf{p}):=\sum_{\genfrac{}{}{0pt}{}{\tau=(\tau_{1},\ldots,\tau_{n})\in\mathbf{N}^{n}}{\vert\tau\vert=\ell}}{a(s,\mathsf{w},m,\tau)\cdot\prod_{i=1}^{n}{(\mathsf{w}_{i}^{\mathsf{v}})^{\tau_{i}}}}\end{equation}de $S^{\ell}(W^{\mathsf{v}})$.
\end{defi}
Tel que nous venons de le d{\'e}finir, le terme $\jet_{W}^{\ell}s(m\mathsf{p})$ d{\'e}pend du choix de la base $\mathsf{w}$. Toutefois, lorsque $\jet_{W}^{h}s(m\mathsf{p})=0$ pour tout entier $h\in\{0,\ldots,\ell-1\}$, ce n'est plus le cas. Consid{\'e}rons le couple $(m,\ell)\in\mathbf{N}^{2}$ tel que $m\in\{0,\ldots,(n+t)S\}$, $\ell\in\{0,\ldots,(n+t)T\}$ et $(m,\ell)$ \emph{minimal} pour la propri{\'e}t{\'e} $\jet_{W}^{\ell}s(m\mathsf{p})\ne 0$. L'adjectif minimal s'entend par rapport {\`a} l'ordre lexicographique sur $\mathbf{N}^{2}$. 
\subsubsection*{Seconde torsion de norme} Il serait naturel (et possible) de travailler dans la suite avec $S^{\ell}(W^{\mathsf{v}})$ muni de la structure de puissance sym{\'e}trique d{\'e}finie au \S~\ref{paragraphepuissance}. Toutefois il s'av{\`e}re plus pratique de modifier l{\'e}g{\`e}rement les normes sur $\overline{S^{\ell}(W^{\mathsf{v}})}$ de la mani{\`e}re suivante~: pour un plongement $\sigma$ quelconque de $k$, consid{\'e}rons le coefficient $\alpha_{\sigma}(m(u_{0},u))$ introduit lors de la d{\'e}finition de $\mathsf{A}_{\sigma_{0}}$  (voir~\eqref{defidealpha}, p.~\pageref{defidealpha}). L'{\'e}galit{\'e} $\alpha_{\iota\circ\sigma}(m(u_{0},u))=\alpha_{\sigma}(m(u_{0},u))$ pour tout automorphisme continu $\iota$ de $\mathbf{C}_{p}$ l{\'e}gitime la
\begin{defi}
Le fibr{\'e} ad{\'e}lique hermitien $\overline{\mathfrak{Jet}}$ est le fibr{\'e} sur $k$ d'espace sous-jacent $S^{\ell}(W^{\mathsf{v}})$ et sa norme en $\sigma:k\hookrightarrow\mathbf{C}_{p}$ est donn{\'e}e par~: pour tout $x\in S^{\ell}(W^{\mathsf{v}})\otimes_{\sigma}\mathbf{C}_{p}$, on a $$\Vert x\Vert_{\overline{\mathfrak{Jet}},\sigma}:=\alpha_{\sigma}(m(u_{0},u))\Vert x\Vert_{\overline{S^{\ell}(W^{\mathsf{v}})},\sigma}.$$
\end{defi}
\begin{prop}\label{propositionpentejet} La pente maximale de $\overline{\mathfrak{Jet}}$ est plus petite que $7n^{2}U_{0}/D$.
\end{prop}
Pour {\'e}tablir cette majoration, nous nous appuierons sur deux r{\'e}sultats auxiliaires.
\begin{lemm}\label{lemmepentemaxw}La pente maximale de $\overline{S^{\ell}(W^{\mathsf{v}})}$ est inf{\'e}rieure {\`a} $2n^{2}U_{0}/D$.\end{lemm}
\begin{proof}Nous allons montrer que $\widehat{\mu}_{\mathrm{max}}\left(\overline{W^{\mathsf{v}}}\right)$ est major{\'e} par $(\log 2)/2$. La proposition~\ref{propestipentesym} et le fait que $\ell\le(n+t)U_{0}/D$ permettront alors de conclure. L'application $k$-lin{\'e}aire $\iota:k^{n}\to W$, $\iota(y):=(\lambda(y),y)$, est un isomorphisme d'espaces vectoriels. Les normes $\sigma$-adiques de l'application duale $\iota^{\mathsf{v}}:W^{\mathsf{v}}\to (k^{n})^{\mathsf{v}}$ sont plus petites que $(\sqrt{2})^{\epsilon_{\sigma}}$ avec $\epsilon_{\sigma}:=1$ si $\sigma$ est archim{\'e}dien et $0$ sinon. L'in{\'e}galit{\'e} de pentes $$\widehat{\mu}_{\mathrm{max}}(\overline{W^{\mathsf{v}}})\le\widehat{\mu}_{\mathrm{max}}((k^{n},\vert\cdot\vert_{2})^{\mathsf{v}})+h(\iota^{\mathsf{v}})$$(voir~\cite[lemme~$6.4$]{rendiconti}) entra{\^{\i}}ne le r{\'e}sultat voulu car la premi{\`e}re quantit{\'e} du majorant est nulle et l'autre inf{\'e}rieure {\`a} $(\log 2)/2$ (on notera que $2n((\log 2)/2+2\log n)\le 2n^{2}$).
\end{proof}
\begin{lemm}\label{lemmehauteuruo}Pour tout entier $m$, la hauteur $h_{\overline{E_{0}}}(1,mu_{0})$ de $(1,mu_{0})$ relative au fibr{\'e} hermitien $\overline{E_{0}}=\overline{k\oplus(k^{n}/W_{0})}$ est inf{\'e}rieure {\`a}  $\log m+(2n^{2}\log b)/D$.
\end{lemm}
\begin{proof}
Dans cet {\'e}nonc{\'e}, $u_{0}$ d{\'e}signe la pr{\'e}image de $-(\beta_{1,0},\ldots,\beta_{t,0})$ par l'application $\varphi:k^{n}/W_{0}\to k^{t}$, qui, {\`a} la classe d'{\'e}quivalence de $z=(z_{1},\ldots,z_{n})$ modulo $W_{0}$, associe $(\ell_{1}(z),\ldots,\ell_{t}(z))$ (voir \S~\ref{subsecpreparatifs}, apr{\`e}s la d{\'e}monstration de la proposition~\ref{propositionhauteurwo}). Ainsi, pour tout $z\in k^{n}$ tel que $\ell_{i}(z)=-\beta_{i,0}$ pour tout $i\in\{1,\ldots,t\}$, la norme de $u_{0}$ relative {\`a} $\overline{k^{n}/W_{0}}$ est inf{\'e}rieure {\`a} la norme de $z$ dans $(k^{n},\vert\cdot\vert_{2})$. En particulier, on a $h_{\overline{E_{0}}}(1,mu_{0})\le h_{(k^{n+1},\vert\cdot\vert_{2})}(1,mz)$. Au moyen d'une matrice extraite, il existe un ensemble $\mathsf{I}\subseteq\{1,\ldots,n\}$ {\`a} $t$ {\'e}l{\'e}ments, des d{\'e}terminants $(\Delta_{\mathsf{i}})_{\mathsf{i}\in\mathsf{I}},\Delta$ de matrices {\`a} coefficients dans $\{\pm\beta_{\mathsf{u},\mathsf{v}};\ 0\le\mathsf{u}\le t,\ 1\le\mathsf{v}\le n\}$ et une solution $z$ du syst{\`e}me $\ell_{i}(z)=-\beta_{i,0}$, $1\le i\le t$, telle que $z_{\mathsf{i}}=\Delta_{\mathsf{i}}/\Delta$ pour $\mathsf{i}\in\mathsf{I}$ et $z_{\mathsf{i}}=0$ sinon. On a alors  \begin{equation*}h_{(k^{n+1},\vert\cdot\vert_{2})}(1,mz)=h_{(k^{t+1},\vert\cdot\vert_{2})}(\Delta,m(\Delta_{\mathsf{i}})_{\mathsf{i}\in\mathsf{I}}).\end{equation*}Des majorations grossi{\`e}res des d{\'e}terminants qui apparaissent ici conduisent {\`a} \begin{equation*}\vert(\Delta,m(\Delta_{\mathsf{i}})_{\mathsf{i}\in\mathsf{I}})\vert_{2,\sigma}\le\prod_{\mathsf{u},\mathsf{v}}{\max{\{1,\vert\beta_{\mathsf{u},\mathsf{v}}\vert_{\sigma}\}}} \times\begin{cases}(1+m^{2}t)^{1/2}t! & \text{si $p=\infty$},\\ 1 & \text{si $p\ne\infty$}.\end{cases}\end{equation*}Par cons{\'e}quent, si $t\le n-1$, on a $$h_{(k^{t+1},\vert\cdot\vert_{2})}(\Delta,m(\Delta_{\mathsf{i}})_{\mathsf{i}\in\mathsf{I}})\le(t+1)n\max_{\mathsf{u},\mathsf{v}}{\{h(\beta_{\mathsf{u},\mathsf{v}})\}}+\log m+\log n!$$et la d{\'e}finition de $\log b$ permet de conclure. Si $t=n$ on a $W_{0}=\{0\}$ et n{\'e}cessairement $z=u_{0}=-(\beta_{1,0},\ldots,\beta_{n,0})$. Dans ce cas, on a $$\vert(1,mz)\vert_{2,\sigma}\le(1+m^{2}n)^{\epsilon_{\sigma}/2}\prod_{i=1}^{n}\max{(1,\vert\beta_{i,0}\vert_{\sigma})}$$pour tout plongement $\sigma:k\hookrightarrow\mathbf{C}_{p}$ (rappelons que $\epsilon_{\sigma}\in\{0,1\}$ est nul si et seulement si $p\ne\infty$). On obtient alors $$h_{(k^{n+1},\vert\cdot\vert_{2})}(1,mz)\le\log m+\frac{\log(n+1)}{2}+\frac{n\log b}{D}.$$On conclut en utilisant $(\log b)/D\ge 1$ et $n+(\log(n+1))/2\le 2n^{2}$.
\end{proof}
\begin{proof}[D{\'e}monstration de la proposition~\ref{propositionpentejet}]
Comme la pente maximale est un maximum, il existe un sous-espace $J\subseteq S^{\ell}(W^{\mathsf{v}})$ tel que $\widehat{\mu}_{\mathrm{max}}(\overline{\mathfrak{Jet}})=\widehat{\mu}(J,(\Vert\cdot\Vert_{\overline{\mathfrak{Jet}},\sigma})_{\sigma})$, pente qui est elle-m{\^{e}}me {\'e}gale {\`a} $\widehat{\mu}(J,(\Vert\cdot\Vert_{\overline{S^{\ell}(W^{\mathsf{v}})},\sigma})_{\sigma})-\frac{1}{D}\sum_{\sigma}{\log\alpha_{\sigma}(m(u_{0},u))}$. Dans cette diff{\'e}rence, la premi{\`e}re quantit{\'e} est inf{\'e}rieure {\`a} $\widehat{\mu}_{\mathrm{max}}(\overline{S^{\ell}(W^{\mathsf{v}})})$, qui a {\'e}t{\'e} estim{\'e}e dans le lemme~\ref{lemmepentemaxw} ci-dessus. L'autre, avec le signe moins, est exactement {\'e}gale {\`a} $$D_{0}h_{\overline{k\oplus(k^{n}/W_{0})}}(1,mu_{0})+\sum_{j=1}^{n}{D_{j}h_{(k^{2},\vert\cdot\vert_{2})}(1,\alpha_{j}^{m})}.$$En vertu du lemme~\ref{lemmehauteuruo}, le premier terme dans cette expression est major{\'e} par $D_{0}(\log((n+t)S)+(2n^{2}\log b)/D)\le 3n^{2}U_{0}/D$. Par ailleurs, comme $(1+x^{2})^{1/2}\le\sqrt{2}\max{\{1,\vert x\vert\}}$, on a \begin{equation*}\begin{split}h_{(k^{2},\vert\cdot\vert_{2})}(1,\alpha_{j}^{m})&\le\log\sqrt{2}+h(\alpha_{j}^{m})=\log\sqrt{2}+mh(\alpha_{j})\\ &\le\log\sqrt{2}+(n+t)S\log a_{j}\end{split}\end{equation*}puis \begin{equation*}\begin{split}\sum_{j=1}^{n}{D_{j}h_{(k^{2},\vert\cdot\vert_{2})}(1,\alpha_{j}^{m})}&\le\left(\sum_{j=1}^{n}{D_{j}}\right)\log\sqrt{2}+(n+t)\sum_{j=1}^{n}{D_{j}S\log a_{j}}\\ &\le\frac{U_{0}n}{D}\log\sqrt{2}+n(2n-1)\frac{U_{0}}{D}\\ &\le \frac{2n^{2}U_{0}}{D}\cdotp\end{split}\end{equation*}En regroupant ces diff{\'e}rentes majorations, on trouve \begin{equation*}\widehat{\mu}_{\mathrm{max}}(\overline{\mathfrak{Jet}})\le\frac{2n^{2}U_{0}}{D}+\frac{3n^{2}U_{0}}{D}+\frac{2n^{2}U_{0}}{D}=\frac{7n^{2}U_{0}}{D}\cdotp\end{equation*}
\end{proof}
Dor{\'e}navant, la majeure partie de la d{\'e}monstration du th{\'e}or{\`e}me~\ref{theoremereduit} va consister {\`a} {\'e}valuer chacune des normes $\sigma$-adiques de $\jet_{W}^{\ell}s(m\mathsf{p})$ dans le but d'estimer sa hauteur, relative au fibr{\'e} ad{\'e}lique hermitien $\overline{\mathfrak{Jet}}$.\subsubsection{Estimations archim{\'e}diennes}
\begin{prop}\label{propestimationarchi} Pour tout plongement complexe $\sigma:k\hookrightarrow\mathbf{C}$,  on a \begin{equation*}\Vert\jet_{W}^{\ell}s(m\mathsf{p})\Vert_{\overline{\mathfrak{Jet}},\sigma}\le \exp{\{nU_{0}/D\}}\Vert s\Vert_{\overline{E}_{\alpha},\sigma}.\end{equation*}
\begin{proof}Soit $\sigma:k\hookrightarrow\mathbf{C}$. Dans l'expression~\eqref{defidujet} du jet, choisissons pour $\mathsf{w}$ une base \emph{orthonorm{\'e}e} de $W\otimes_{\sigma}\mathbf{C}$ et fixons un logarithme $mu$ de $m\mathsf{p}\in G(k_{\sigma})$. La norme $\Vert\jet_{W}^{\ell}s(m\mathsf{p})\Vert_{\overline{\mathfrak{Jet}},\sigma}$ est major{\'e}e par \begin{equation}\label{majorationdelanormedujet}\binom{n+\ell-1}{\ell}^{1/2}\max_{\vert\tau\vert=\ell}{\left\{\left\vert\frac{1}{\tau !}\mathrm{D}_{\mathsf{w}}^{\tau}F_{s,\sigma}(mu)\right\vert_{\sigma}\alpha_{\sigma}(m(u_{0},u))\right\}}\cdotp\end{equation}La racine du facteur binomial est inf{\'e}rieure {\`a} $$(\ell+1)^{(n-1)/2}\le((n+t)T+1)^{(n-1)/2}$$ et donc {\`a} $$(2n+1)^{(n-1)/2}\left(\frac{U_{0}}{D}\right)^{(n-1)/2}\le (2n+1)^{(n-1)/2}(n-1)!\,e^{\sqrt{U_{0}/D}}\le e^{U_{0}/(2D)}$$ car $T\le U_{0}/D$ et $x^{n-1}\le(n-1)!e^{x}$ pour $x\ge 0$. La proposition~\ref{propmajorationarchimediennesimple} majore la d{\'e}riv{\'e}e de $F_{s,\sigma}$ et l'on obtient alors $$\Vert\jet_{W}^{\ell}s(m\mathsf{p})\Vert_{\overline{\mathfrak{Jet}},\sigma}\le e^{U_{0}/(2D)}e^{2n\sqrt{n}\max_{0\le i\le n}{\{D_{i}\}}}L(s).$$La longueur de $s$ qui appara{\^{\i}}t ici peut {\^{e}}tre estim{\'e}e comme {\`a} la fin du \S~\ref{paragraphepuissance}, en utilisant~\eqref{longueurs}, et l'on a $L(s)\le\Vert s\Vert_{\overline{E}_{\alpha},\sigma}\exp{\{n\max_{0\le i\le n}{D_{i}}\}}$. Le choix des param{\`e}tres montre que le maximum des $D_{i}$ est plus petit que $U_{0}/(D\log C_{0})$, ce qui permet de conclure.
\end{proof}
\end{prop}
\subsubsection{Estimations ultram{\'e}triques}
\'Etant donn{\'e} des entiers $\ell$ et $h$ strictement positifs, on d{\'e}finit l'entier $\delta_{\ell}(h)$ comme le ppcm des produits $i_{1}\cdots i_{h'}$ o{\`u} $h'\in\{1,\ldots,h\}$, $i_{j}\in\mathbf{N}\setminus\{0\}$ et $i_{1}+\cdots+i_{h'}\le\ell$. Le th{\'e}or{\`e}me des nombres premiers assure l'existence d'une constante absolue $c>0$ telle que $\log\delta_{\ell}(h)\le\ell\log(ch)$. Au \S~\ref{sectionconclusion}, nous utiliserons la valeur $c=4$ obtenue par Bruiltet dans~\cite{bruiltet}.
\begin{prop}\label{propestimationultra}
Pour tout plongement ultram{\'e}trique $\sigma:k\hookrightarrow\mathbf{C}_{p}$, pour tous entiers $m,\ell\ge 0$, pour toute section globale $s$ de $M$, le jet de $s$, d'ordre $\ell$ le long de $W$ au point $m\mathsf{p}$, est de norme $\sigma$-adique inf{\'e}rieure {\`a} $\Vert s\Vert_{\overline{E},\sigma}/\vert\delta_{\ell}(D_{0})\vert_{p}$. 
\end{prop}
La d{\'e}monstration est une variante plus simple de l'{\'e}nonc{\'e} {\'e}quivalent montr{\'e} dans le cadre des vari{\'e}t{\'e}s ab{\'e}liennes au \S~$5.8$ de~\cite{artepredeux}, au moyen du proc{\'e}d{\'e} de changement de variables de Chudnovsky. Nous renvoyons le lecteur {\`a} cet article pour plus de pr{\'e}cisions.
\subsection{Extrapolation sur les d{\'e}rivations}\label{section:extrapolation}
L'objectif de ce paragraphe est d'{\'e}valuer, pour tout plongement $\sigma$ conjugu{\'e} {\`a} $\sigma_{0}$, la $\sigma$-norme de $\jet_{W}^{\ell}s(m\mathsf{p})$ en tenant compte de la construction de $s$. Comme les $\sigma\in v_{0}$ jouent tous le m{\^{e}}me r{\^{o}}le, il suffit d'{\'e}tablir les estimations pour $\sigma=\sigma_{0}$ (ce que nous supposons dor{\'e}navant). Comme de coutume, nous allons devoir distinguer les cas p{\'e}riodique et non p{\'e}riodique, distinction {\`a} laquelle se superposera la distinction entre $\sigma_{0}$ plongement complexe et $\sigma_{0}$ ultram{\'e}trique.\par Soit $w$ la base particuli{\`e}re de $W\otimes_{\sigma_{0}}\mathbf{C}_{p_{0}}$ introduite au \S~\ref{sectiontroissix}, utilis{\'e}e dans la d{\'e}finition de la matrice $\mathsf{A}_{\sigma_{0}}$. D'apr{\`e}s la majoration~\eqref{majorationdelanormedujet} appliqu{\'e}e avec $\sigma=\sigma_{0}$, $\mathsf{w}=w$ et $z=m(u_{0},u)$, le probl{\`e}me est d'estimer $\frac{1}{\tau!}\mathrm{D}_{w}^{\tau}F_{s,\sigma_{0}}(m(u_{0},u))$. Si $(m,\tau)\in\Upsilon$ alors la d{\'e}finition~\eqref{defidescoefficientsdeA} des coefficients de la matrice $\mathsf{A}_{\sigma_{0}}$ implique \begin{equation}\label{ineqmajconstruction}\left\vert\frac{1}{\tau!}\mathrm{D}_{w}^{\tau}F_{s,\sigma_{0}}(m(u_{0},u))\right\vert_{p_{0}}\alpha_{\sigma_{0}}(m(u_{0},u))\le\alpha^{-1}\Vert s\Vert_{\overline{E}_{\alpha},\sigma_{0}}\end{equation}et l'on obtient directement une estimation de $\frac{1}{\tau!}\mathrm{D}_{w}^{\tau}F_{s,\sigma_{0}}(m(u_{0},u))$, qui est meilleure que celle que l'on aura pour $(m,\tau)\not\in\Upsilon$. C'est la raison pour laquelle nous supposerons maintenant que $(m,\tau)\not\in\Upsilon$. On ram{\`e}ne alors le probl{\`e}me {\`a} $\frac{1}{\tau!}\mathrm{D}_{w}^{\tau}F_{s,\sigma_{0}}(m(\lambda(u),u))$ gr{\^a}ce au lemme de comparaison suivant~:
\begin{lemm}\label{lemmecomparaison} Soit $\mathsf{w}=(\mathsf{w}_{1},\ldots,\mathsf{w}_{n})$ une base de $W\otimes_{\sigma_{0}}\mathbf{C}_{p_{0}}$. Pour tout $\mathsf{m}\in\mathbf{N}$ tel que $\Vert\mathsf{m}(u_{0}-\lambda(u))\Vert_{\overline{t_{G}},\sigma_{0}}\le 1$ et pour tout $\tau=(\tau_{1},\ldots,\tau_{n})\in\mathbf{N}^{n}$, de longueur $\le 2(n+t)T$, la valeur absolue $p_{0}$-adique de la diff{\'e}rence $\frac{1}{\tau!}\mathrm{D}_{\mathsf{w}}^{\tau}F_{s,\sigma_{0}}(\mathsf{m}(u_{0},u))-\frac{1}{\tau!}\mathrm{D}_{\mathsf{w}}^{\tau}F_{s,\sigma_{0}}(\mathsf{m}(\lambda(u),u))$ est major{\'e}e par \begin{equation*}e^{4nU_{0}}\left\{\prod_{i=1}^{n}{\Vert\mathsf{w}_{i}\Vert_{\overline{t_{G}},\sigma_{0}}^{\tau_{i}}}\right\}\Vert\mathsf{m}(u_{0}-\lambda(u))\Vert_{\overline{t_{G}},\sigma_{0}}\alpha_{\sigma_{0}}(\mathsf{m}(u_{0},u))^{-1}\Vert s\Vert_{\overline{E}_{\alpha},\sigma_{0}}.\end{equation*}\end{lemm}
\begin{proof}Reprenons le d{\'e}but de la d{\'e}monstration de la proposition~\ref{propmajorationarchimediennesimple} (avec $\sigma=\sigma_{0}$) et les formes lin{\'e}aires $\mathsf{l}_{1},\ldots,\mathsf{l}_{n}$ associ{\'e}es {\`a} la base $\mathsf{w}$~:
$$\forall\,x=(x_{1},\ldots,x_{n})\in\mathbf{C}_{p_{0}}^{n},\quad\sum_{i=1}^{n}{x_{i}\mathsf{w}_{i}}=(\lambda(\mathsf{l}_{1}(x),\ldots,\mathsf{l}_{n}(x)),\mathsf{l}_{1}(x),\ldots,\mathsf{l}_{n}(x)).$$On peut supposer $\Vert\mathsf{w}_{i}\Vert_{\overline{t_{G}},\sigma_{0}}=1$ pour tout $i\in\{1,\ldots,n\}$. D{\'e}signons par $(p_{\mathsf{h}})_{\mathsf{h}\in\mathbf{N}^{t}\times\mathbf{N}^{n}}$ les coefficients de $s$ dans une {\'e}criture du type~\eqref{ecrituredes}. La diff{\'e}rence $\frac{1}{\tau!}\mathrm{D}_{\mathsf{w}}^{\tau}F_{s,\sigma_{0}}(\mathsf{m}(u_{0},u))-\frac{1}{\tau!}\mathrm{D}_{\mathsf{w}}^{\tau}F_{s,\sigma_{0}}(\mathsf{m}(\lambda(u),u))$ est la somme sur $\mathsf{h}$ des coefficients devant $x^{\tau}$ des fonctions \begin{equation}\begin{split}\label{termeimportant}&p_{\mathsf{h}}\times\left(\prod_{i=1}^{n}{\alpha_{i}^{\mathsf{m}\mathsf{h}_{i}}}\right)\times e^{\mathsf{h}_{1}\mathsf{l}_{1}(x)+\cdots+\mathsf{h}_{n}\mathsf{l}_{n}(x)}\\ &\times\left\{\prod_{a=1}^{t}{\left(\mathsf{m}e_{0,a}^{\mathsf{v}}(u_{0})+\sum_{i=1}^{n}{\mathsf{l}_{i}(x)e_{0,a}^{\mathsf{v}}(\lambda(e_{i}))}\right)^{\mathsf{h}_{0,a}}}\right.\\ &\qquad\quad\left.-\prod_{a=1}^{t}{\left(\mathsf{m}e_{0,a}^{\mathsf{v}}(\lambda(u))+\sum_{i=1}^{n}{\mathsf{l}_{i}(x)e_{0,a}^{\mathsf{v}}(\lambda(e_{i}))}\right)^{\mathsf{h}_{0,a}}}\right\}.\end{split}\end{equation}Dans cette expression, $\mathsf{h}_{0}=(\mathsf{h}_{0,1},\ldots,\mathsf{h}_{0,t})\in\mathbf{N}^{t}$, $\mathsf{h}_{1},\ldots,\mathsf{h}_{n}\in\mathbf{N}$ sont les coordonn{\'e}es de $\mathsf{h}$. Comme $\mathsf{h}_{0}$ est de longueur $D_{0}$, la diff{\'e}rence des quantit{\'e}s entre accolades est une expression du type \begin{equation}\label{expressiondeladifference}\prod_{i=1}^{D_{0}}{(a_{i}+b_{i})}-\prod_{i=1}^{D_{0}}{(a_{i}'+b_{i})}\end{equation}o{\`u} $a_{i},a_{i}'$ sont les composantes de $\mathsf{m}u_{0},\mathsf{m}\lambda(u)$ respectivement dans la base orthonorm{\'e}e $e_{0}:=(e_{0,1},\ldots,e_{0,t})$ et o{\`u} $b_{i}$ est l'une des composantes de $\mathsf{l}(x):=\lambda(\mathsf{l}_{1}(x),\ldots,\mathsf{l}_{n}(x))$ dans la base $e_{0}$ (r{\'e}p{\'e}t{\'e}e le bon nombre de fois pour avoir la puissance $\mathsf{h}_{0,a}$ dans~\eqref{termeimportant}). Or la diff{\'e}rence~\eqref{expressiondeladifference} est {\'e}gale {\`a} $$\sum_{j=1}^{D_{0}}{(a_{j}-a_{j}')\left(\prod_{i=1}^{j-1}{(a_{i}'+b_{i})}\right)\left(\prod_{i=j+1}^{D_{0}}{(a_{i}+b_{i})}\right)}$$(l'on {\'e}crit $a_{j}-a_{j}'=(a_{j}+b_{j})-(a_{j}'+b_{j})$ et l'on d{\'e}veloppe; les produits se simplifient deux {\`a} deux). Ainsi la diff{\'e}rence des produits entre accolades dans~\eqref{termeimportant}, que l'on peut {\'e}crire symboliquement sous la forme $$\left(\mathsf{m}e_{0}^{\mathsf{v}}(u_{0})+e_{0}^{\mathsf{v}}(\mathsf{l}(x))\right)^{\mathsf{h}_{0}}-\left(\mathsf{m}e_{0}^{\mathsf{v}}(\lambda(u))+e_{0}^{\mathsf{v}}(\mathsf{l}(x))\right)^{\mathsf{h}_{0}},$$est une somme de $D_{0}$ termes de la forme $$\delta(a,\mathsf{h}_{0})(x):=e_{0,a}^{\mathsf{v}}(\mathsf{m}(u_{0}-\lambda(u)))\times\left(\mathsf{m}e_{0}^{\mathsf{v}}(u_{0})+e_{0}^{\mathsf{v}}(\mathsf{l}(x))\right)^{\mathsf{h}_{0}'}\left(\mathsf{m}e_{0}^{\mathsf{v}}(\lambda(u))+e_{0}^{\mathsf{v}}(\mathsf{l}(x))\right)^{\mathsf{h}_{0}''}$$o{\`u} $\mathsf{h}_{0}',\mathsf{h}_{0}''$ sont des vecteurs entiers (qui d{\'e}pendent de $a$) dont la somme des longueurs vaut $D_{0}-1$. On notera que ce coefficient $\delta(a,\mathsf{h}_{0})$ est une fonction polynomiale de la variable $x$. La suite de la d{\'e}monstration repose alors sur les m{\^{e}}mes consid{\'e}rations que celle de la proposition~\ref{propmajorationarchimediennesimple}. Plus pr{\'e}cis{\'e}ment, soit $Q_{\mathsf{h},\tau}$ le polyn{\^{o}}me de Taylor d{\'e}fini par~\eqref{definitiondeqhtau}. Le coefficient de $x^{\tau}$ de~\eqref{termeimportant} est {\'e}gal {\`a} celui du polyn{\^{o}}me \begin{equation}\sum_{\text{$D_{0}$ termes}}{p_{\mathsf{h}}\delta(a,\mathsf{h}_{0})(x)\left(\prod_{i=1}^{n}{\alpha_{i}^{\mathsf{m}\mathsf{h}_{i}}}\right)Q_{\mathsf{h},\tau}(x)}.\end{equation}Ce coefficient est en valeur absolue plus petit que la longueur du polyn{\^{o}}me. Ainsi la quantit{\'e} $\left\vert\frac{1}{\tau!}\mathrm{D}_{\mathsf{w}}^{\tau}F_{s,\sigma_{0}}(\mathsf{m}(u_{0},u))-\frac{1}{\tau!}\mathrm{D}_{\mathsf{w}}^{\tau}F_{s,\sigma_{0}}(\mathsf{m}(\lambda(u),u))\right\vert_{p_{0}}$ est inf{\'e}rieure {\`a} \begin{equation}\begin{split}&\Vert\mathsf{m}u_{0}-\mathsf{m}\lambda(u)\Vert_{\overline{t_{G}},\sigma_{0}}\times L(s)\max_{\mathsf{h}}{\{L(Q_{\mathsf{h},\tau})\}}\times\prod_{i=1}^{n}{\max{\{1,\vert\alpha_{i}\vert_{p_{0}}\}}^{\mathsf{m}D_{i}}}\\ &\qquad\times\begin{cases}D_{0}(\mathsf{m}\Vert u_{0}\Vert_{\overline{t_{G}},\sigma_{0}}+1+n)^{D_{0}-1} & \text{si $p_{0}=\infty$}\\ \max{\{1,\Vert\mathsf{m}u_{0}\Vert_{\overline{t_{G}},\sigma_{0}}\}}^{D_{0}-1} & \text{si $p_{0}\ne\infty$}\end{cases}\end{split}
\end{equation}(l'hypoth{\`e}se $\Vert\mathsf{m}(u_{0}-\lambda(u))\Vert_{\overline{t_{G}},\sigma_{0}}\le 1$ a permis ici de majorer la norme de $\mathsf{m}\lambda(u)$ et l'estimation de la longueur de $\mathsf{l}$ repose sur les estimations des sommes $\sum_{j=1}^{n}{\mathsf{l}_{i}(e_{j})\lambda(e_{i})}$ introduites lors de la d{\'e}monstration de la proposition~\ref{propmajorationarchimediennesimple}). Pour conclure et si $\sigma_{0}$ est complexe, on majore $L(s)$ par $e^{U_{0}}\Vert s\Vert_{\overline{E}_{\alpha},\sigma_{0}}$ (fin du \S~\ref{paragraphepuissance}), on majore $\max_{\mathsf{h}}{\{L(Q_{\mathsf{h},\tau})\}}$ par $e^{U_{0}}$ gr{\^{a}}ce {\`a}~\eqref{eqmajorantdeq} et l'on majore $D_{0}(\mathsf{m}\Vert u_{0}\Vert_{\overline{t_{G}},\sigma_{0}}+1+n)^{D_{0}-1}$ par \begin{equation*}D_{0}(n+2)^{D_{0}-1}\max{\{1,\Vert\mathsf{m}u_{0}\Vert_{\overline{t_{G}},\sigma_{0}}\}}^{D_{0}}\le e^{2nU_{0}}\moy_{\sigma_{0}}(1,\Vert\mathsf{m}u_{0}\Vert_{\overline{t_{G}},\sigma_{0}})^{D_{0}}.\end{equation*}Si $\sigma_{0}$ est ultram{\'e}trique, on a $L(s)=\Vert s\Vert_{\overline{E}_{\alpha},\sigma_{0}}$ et $\max_{\mathsf{h}}{\{L(Q_{\mathsf{h},\tau})\}}\le e^{\vert\tau\vert}\le e^{2(n+t)T}\le e^{4nU_{0}}$.
\end{proof}
L'estimation de $\frac{1}{\tau!}\mathrm{D}_{w}^{\tau}F_{s,\sigma_{0}}(m(\lambda(u),u))$ repose sur une version l{\'e}g{\`e}rement affaiblie des lemmes d'interpolation de Waldschmidt~\cite{miw1980} (cas archim{\'e}dien) et de Roy~\cite{Roy2001} (cas ultram{\'e}trique, voir aussi~\cite[lemme~$4.21$]{casrationnel} pour les simplifications faites ici). Si $x$ est un nombre r{\'e}el positif et $f$ une fonction d{\'e}finie sur le disque ferm{\'e} $\overline{D}(0,x)=\{z\in\mathbf{C}_{p_{0}}\,;\ \vert z\vert_{p_{0}}\le x\}$, on note $\vert f\vert_{x}$ la borne sup{\'e}rieure des $\vert f(z)\vert$, $z\in\overline{D}(0,x)$. Rappelons que $\epsilon_{0}=0$ si $\sigma_{0}$ est ultram{\'e}trique et $1$ sinon.
\begin{lemm}\label{lemmedinterpolation}Soit $S_{1},T_{1}$ des entiers naturels $\ge 2$ et $\mathsf{R},\mathsf{r}$ des nombres r{\'e}els v{\'e}rifiant $\mathsf{R}\ge\mathsf{r}\ge(2S_{1})^{\epsilon_{0}}$. Soit $f$ une fonction analytique dans le disque $\overline{D}(0,\mathsf{R})$. On a alors \begin{equation*}\begin{split}\vert f\vert_{\mathsf{r}}\le S_{1}^{T_{1}}\mathrm{max}&\left\{\left(\frac{(2r_{p_{0}})^{\epsilon_{0}}\mathsf{r}}{r_{p_{0}}\mathsf{R}}\right)^{T_{1}S_{1}}\vert f\vert_{\mathsf{R}},\right.\\ &\quad\left.\left(\left(\frac{10r_{p_{0}}}{S_{1}}\right)^{\epsilon_{0}}\frac{\mathsf{r}}{r_{p_{0}}}\right)^{T_{1}S_{1}}\max_{\genfrac{}{}{0pt}{}{0\le h<T_{1}}{0\le\mathsf{m}<S_{1}}}{\left\{\left\vert\frac{1}{h!}f^{(h)}(\mathsf{m})\right\vert_{\sigma_{0}}\right\}}\right\}\cdotp\end{split}\end{equation*} 
\end{lemm}
La suite de la d{\'e}monstration repose sur un raisonnement par l'absurde, au moyen de l'hypoth{\`e}se suivante.
\begin{enonce}{Hypoth{\`e}se}\label{hypocru}
$\Vert u_{0}-\lambda(u)\Vert_{\overline{t_{G}},\sigma_{0}}\le \exp{\left(-(C_{0}^{3/2}+4n)U_{0}\right)}$. 
\end{enonce}
Notons $\mathbf{D}$ le disque ouvert $\{z\in\mathbf{C}_{p_{0}}\,;\ \vert z\vert_{p_{0}}<\min{\{r_{p_{0}}/\vert u_{i}\vert_{p_{0}}\,;\ 1\le i\le n\}}\}$ si $\sigma_{0}$ est ultram{\'e}trique et l'ensemble $\mathbf{C}$ des nombres complexes sinon. On notera que les disques ferm{\'e}s, centr{\'e}s en $0$ et de rayons $1$ et $\mathfrak{e}/r_{p_{0}}$ respectivement, sont tous les deux inclus dans $\mathbf{D}$.
\subsubsection{Cas non p{\'e}riodique}
Soit $\tau\in\mathbf{N}^{n}$ de longueur $\ell$ (on rappelle que cet entier $\ell$, fix{\'e} apr{\`e}s la d{\'e}finition~\ref{definitionjet}, est inf{\'e}rieur {\`a} $(n+t)T$). Soit $\mathrm{f}:\mathbf{D}\to\mathbf{C}_{p_{0}}$ l'application analytique d{\'e}finie par~: $$\forall\,z\in\mathbf{D},\quad\mathrm{f}(z)=\frac{1}{\tau !}\mathrm{D}_{w}^{\tau}F_{s,\sigma_{0}}(z(\lambda(u),u)).$$ Observons que, pour tout $h\in\mathbf{N}$ et tout $z\in\mathbf{D}$, on a \begin{equation}\begin{split}\label{deriveedef}\frac{\mathrm{f}^{(h)}(z)}{h!} & =\frac{1}{\tau!h!}\mathrm{D}_{w}^{\tau}\mathrm{D}_{(\lambda(u),u)}^{h}F_{s,\sigma_{0}}(z(\lambda(u),u))\\ &=\sum_{\genfrac{}{}{0pt}{}{j\in\mathbf{N}^{n}}{\vert j\vert=h}}{\binom{\tau+j}{\tau}\mathfrak{u}^{j}\frac{\mathrm{D}_{w}^{\tau+j}}{(\tau+j)!}F_{s,\sigma_{0}}(z(\lambda(u),u))}\end{split}\end{equation}($\mathfrak{u}=(\mathfrak{u}_{1},\ldots,\mathfrak{u}_{n})\in\mathbf{C}_{p_{0}}^{n}$ est le vecteur des coordonn{\'e}es de $(\lambda(u),u)$ dans la base $(w_{1},\ldots,w_{n})$). Pour estimer $\frac{1}{\tau!}\mathrm{D}_{w}^{\tau}F_{s,\sigma_{0}}(m(u_{0},u))$, nous allons utiliser le lemme d'interpolation~\ref{lemmedinterpolation} afin de majorer $\vert\mathrm{f}(m)\vert_{p_{0}}$, avec les param{\`e}tres suivants~:\begin{enumerate}\item[(i)] Si $\sigma_{0}$ est un plongement \emph{complexe} alors $S_{1}:=S_{0}$, $\mathsf{r}:=2m$, $\mathsf{R}:=2\mathsf{r}\mathfrak{e}$, $T_{1}:=(n+t)T$ et $f:=\mathrm{f}$.\item[(ii)] Si $\sigma_{0}$ est un plongement \emph{ultram{\'e}trique} alors $S_{1}:=S_{0}$, $\mathsf{r}:=1$, $\mathsf{R}:=\mathfrak{e}/r_{p_{0}}$, $T_{1}:=(n+t)T$ et $f:=\mathrm{f}$.\end{enumerate}Avec ces choix, la proposition~\ref{propmajorationarchimediennesimple} fournit l'estimation suivante~:
\begin{prop}\label{propositiondefr} On a $\alpha_{\sigma_{0}}(m(u_{0},u))\vert\mathrm{f}\vert_{\mathsf{R}}\le \exp{\{17n^{2}U_{0}\}}\Vert s\Vert_{\overline{E}_{\alpha},\sigma_{0}}$.
\end{prop}
\begin{proof}En remarque liminaire, mentionnons que les majorations qui suivent sont volontairement tr{\`e}s larges, afin de rester valide dans le cas p{\'e}riodique, qui sera {\'e}tudi{\'e} au prochain paragraphe.\par 
D'apr{\`e}s la proposition~\ref{propmajorationarchimediennesimple}, on a \begin{equation*}\vert\mathrm{f}(z)\vert_{p_{0}}\cdot\alpha_{\sigma_{0}}(z(\lambda(u),u))\le L(s)\times\begin{cases}e^{2n\sqrt{n}\max_{0\le j\le n}{D_{j}}} & \text{si $p_{0}=\infty$}\\ e^{4nT} & \text{si $p_{0}\ne\infty$}.\end{cases}\end{equation*}Comme nous l'avons d{\'e}j{\`a} vu (voir fin de la d{\'e}monstration de la proposition~\ref{propestimationarchi}), la longueur de $s$ est plus petite que $e^{n\max_{0\le j\le n}{D_{j}}}\Vert s\Vert_{\overline{E}_{\alpha},\sigma_{0}}$ si $p_{0}=\infty$ et $L(s)\le\Vert s\Vert_{\overline{E}_{\alpha},\sigma_{0}}$ si $p_{0}\ne\infty$. Le choix des param{\`e}tres donne alors \begin{equation*}\vert\mathrm{f}(z)\vert_{p_{0}}\cdot\alpha_{\sigma_{0}}(z(\lambda(u),u))\le e^{4nU_{0}}\Vert s\Vert_{\overline{E}_{\alpha},\sigma_{0}}.\end{equation*}En utilisant la d{\'e}finition~\eqref{defidealpha} de $\alpha_{\sigma_{0}}$, la quantit{\'e} $\alpha_{\sigma_{0}}(z(\lambda(u),u))^{-1}$ est, pour $\vert z\vert_{p_{0}}\le\mathsf{R}$, major{\'e}e par  \begin{equation*}\begin{cases}(1+64(nS\mathfrak{e}\Vert\lambda(u)\Vert_{\overline{t_{G}},\sigma_{0}})^{2})^{D_{0}/2}\prod_{j=1}^{n}{(1+e^{16nS\mathfrak{e}\vert u_{j}\vert_{p_{0}}})^{D_{j}/2}} & \text{si $p_{0}=\infty$},\\ \max{\{1,2\mathfrak{e}\Vert\lambda(u)\Vert_{\overline{t_{G}},\sigma_{0}}\}}^{D_{0}} & \text{si $p_{0}\ne\infty$}\end{cases}\end{equation*}(dans le deuxi{\`e}me cas, nous avons major{\'e} $r_{p_{0}}^{-1}=p_{0}^{\frac{1}{p_{0}-1}}$ par $2$). Dans le cas ultram{\'e}trique l'hypoth{\`e}se~\ref{hypocru} implique $\Vert\lambda(u)\Vert_{\overline{t_{G}},\sigma_{0}}=\Vert u_{0}\Vert_{\overline{t_{G}},\sigma_{0}}$ et $$\max{\{1,2\mathfrak{e}\Vert u_{0}\Vert_{\overline{t_{G}},\sigma_{0}}\}}^{D_{0}}\le(2\mathfrak{e})^{D_{0}}\max{\{1,\Vert u_{0}\Vert_{\overline{t_{G}},\sigma_{0}}\}}^{D_{0}}.$$On trouve ainsi $$\frac{\alpha_{\sigma_{0}}(z(\lambda(u),u))^{-1}}{\alpha_{\sigma_{0}}(m(u_{0},u))^{-1}}\le(2\mathfrak{e})^{D_{0}}\le e^{2U_{0}}.$$Le r{\'e}sultat est alors acquis car $2+4n\le 17n^{2}$. Dans le cas archim{\'e}dien, on majore $\Vert\lambda(u)\Vert_{\overline{t_{G}},\sigma_{0}}^{2}$ par $2(\Vert\lambda(u)-u_{0}\Vert_{\overline{t_{G}},\sigma_{0}}^{2}+\Vert u_{0}\Vert_{\overline{t_{G}},\sigma_{0}}^{2})$. Le quotient $$\frac{\alpha_{\sigma_{0}}(z(\lambda(u),u))^{-1}}{\alpha_{\sigma_{0}}(m(u_{0},u))^{-1}}$$est major{\'e} par 
\begin{equation*}\sqrt{2}^{\sum_{j=0}^{n}{D_{j}}}(8n\mathfrak{e}S)^{D_{0}}(1+\Vert u_{0}\Vert_{\overline{t_{G}},\sigma_{0}}^{2})^{D_{0}/2}\exp{\left\{8nS\mathfrak{e}\left(\sum_{j=1}^{n}{D_{j}\vert u_{j}\vert_{p_{0}}}\right)\right\}}.\end{equation*}Le choix des param{\`e}tres entra{\^{\i}}ne \begin{enumerate}\item[(i)] $8nS\mathfrak{e}\left(\sum_{j=1}^{n}{D_{j}\vert u_{j}\vert_{p_{0}}}\right)\le 8n^{2}U_{0 }$,\item[(ii)] $\sqrt{2}^{\sum_{j=0}^{n}{D_{j}}}(8n\mathfrak{e}S)^{D_{0}}\le e^{3U_{0}}$,\item[(iii)] $\moy_{\sigma_{0}}(1,\Vert u_{0}\Vert_{\overline{t_{G}},\sigma_{0}})^{D_{0}}\le\exp{(DD_{0}h_{\overline{E_{0}}}(1,u_{0}))}\le e^{2n^{2}U_{0}}$ (lemme~\ref{lemmehauteuruo}).\end{enumerate}On trouve ainsi \begin{equation}\label{majdealpha}\frac{\alpha_{\sigma_{0}}(z(\lambda(u),u))^{-1}}{\alpha_{\sigma_{0}}(m(u_{0},u))^{-1}}\le e^{(10n^{2}+3)U_{0}}
\quad (\text{pour $\vert z\vert_{p_{0}}\le\mathsf{R}$})\end{equation}et l'on conclut avec $10n^{2}+3+4n\le 17n^{2}$.
\end{proof}
Par ailleurs, au moyen de la formule~\eqref{deriveedef}, du lemme~\ref{lemmecomparaison} et de la construction de la section $s$, on a \begin{equation*}\max_{\genfrac{}{}{0pt}{}{0\le h<T_{1}}{0\le\mathsf{m}<S_{1}}}{\left\{\left\vert\frac{1}{h!}\mathrm{f}^{(h)}(\mathsf{m})\right\vert_{p_{0}}\right\}}\le \theta_{0}\left(e^{4nU_{0}}\Vert u_{0}-\lambda(u)\Vert_{\overline{t_{G}},\sigma_{0}}+\alpha^{-1}\right)\Vert s\Vert_{\overline{E}_{\alpha},\sigma_{0}}\end{equation*}o{\`u} \begin{equation*}\theta_{0}:=\max_{0\le\mathsf{m}<S_{0}}{\{\alpha_{\sigma_{0}}(\mathsf{m}(u_{0},u))^{-1}\}}\times\begin{cases} (n+\sum_{j=1}^{n}{\vert\mathfrak{u}_{j}\vert_{p_{0}}})^{2(n+t)T} & \text{si $p_{0}=\infty$},\\ 1 & \text{si $p_{0}\ne\infty$}.\end{cases}\end{equation*}La somme des $\vert\mathfrak{u}_{j}\vert_{p_{0}}$ est major{\'e}e par $\sqrt{n}\Vert(\lambda(u),u)\Vert_{\overline{t_{G}},\sigma_{0}}$, lui-m{\^{e}}me plus petit que $\sqrt{2n}\vert u\vert_{2,\sigma_{0}}$. On a $\vert u\vert_{2,\sigma_{0}}\le\sum_{j=1}^{n}{\vert u_{j}\vert_{p_{0}}}$ et avec la d{\'e}finition de $\log a_{j}$ on trouve \begin{equation*}\begin{split}n+\sum_{j=1}^{n}{\vert\mathfrak{u}_{j}\vert_{p_{0}}}&\le n+\sqrt{2n}\left(\frac{D}{\mathfrak{e}}\right)\log\prod_{j=1}^{n}{a_{j}}\\ & \le n+\sqrt{2n}\exp{\left\{2\mathfrak{a}\log\mathfrak{e}\right\}}.
\end{split}\end{equation*}En posant $\mathfrak{x}:=\mathfrak{a}\log\mathfrak{e}\ge 1$, on a $$(n+\sqrt{2n}e^{2\mathfrak{x}})^{4n/\mathfrak{x}}\le e^{8n}(ne^{-2}+\sqrt{2n})^{4n}\le e^{10n^{2}}.$$La majoration $T\le U_{0}/(\mathfrak{a}\log\mathfrak{e})$ entra{\^{\i}}ne alors \begin{equation*}\theta_{0}\le\max_{0\le\mathsf{m}<S_{1}}{\{\alpha_{\sigma_{0}}(\mathsf{m}(u_{0},u))^{-1}\}}\exp{\left\{10n^{2}U_{0}\right\}}\cdot\end{equation*}En proc{\'e}dant de la m{\^{e}}me mani{\`e}re que pour la majoration~\eqref{majdealpha}, on a, pour tout entier $\mathsf{m}$ entre $1$ et $S_{1}$, \begin{equation*}\begin{split}\alpha_{\sigma_{0}}(\mathsf{m}(u_{0},u))^{-1}&\le (2nS)^{D_{0}}e^{2n^{2}U_{0}}\prod_{j=1}^{n}{\left(\sqrt{2}e^{D\mathsf{m}\log a_{j}}\right)^{D_{j}}}\\ &\le \exp{\left(\left(1+\frac{\log(2n)}{\log C_{0}}+2n^{2}+\frac{n\log\sqrt{2}}{C_{0}}+n(n+t)\right)U_{0}\right)}\\ &\le e^{5n^{2}U_{0}}.\end{split}\end{equation*}Comme $\alpha_{\sigma_{0}}(m(u_{0},u))\le 1$, on en d{\'e}duit \begin{equation*}\alpha_{\sigma_{0}}(m(u_{0},u))\cdot\max_{0\le\mathsf{m}<S_{1}}{\{\alpha_{\sigma_{0}}(\mathsf{m}(u_{0},u))^{-1}\}} \le e^{5n^{2}U_{0}}.\end{equation*}De plus $\log\alpha$ a {\'e}t{\'e} choisi sup{\'e}rieur {\`a} $C_{0}^{3/2}U_{0}$ (proposition~\ref{propconstructiondes}) et $\Vert u_{0}-\lambda(u)\Vert_{\overline{t_{G}},\sigma_{0}}$ est plus petit que $e^{-(C_{0}^{3/2}+4n)U_{0}}$ d'apr{\`e}s l'hypoth{\`e}se~\ref{hypocru}. Les calculs ci-dessus montrent alors que \begin{align}\label{petitessealpha}\alpha_{\sigma_{0}}(m(u_{0},u))\max_{\genfrac{}{}{0pt}{}{0\le h<T_{1}}{0\le\mathsf{m}<S_{1}}}{\left\{\left\vert\frac{1}{h!}\mathrm{f}^{(h)}(\mathsf{m})\right\vert_{p_{0}}\right\}}&\le e^{15n^{2}U_{0}}\left(e^{4nU_{0}}\Vert u_{0}-\lambda(u)\Vert_{\overline{t_{G}},\sigma_{0}}+\alpha^{-1}\right)\Vert s\Vert_{\overline{E}_{\alpha},\sigma_{0}}\\ &\le e^{(16n^{2}-C_{0}^{3/2})U_{0}}\Vert s\Vert_{\overline{E}_{\alpha},\sigma_{0}}.\end{align}Nous sommes maintenant en mesure d'appliquer le lemme d'interpolation~\ref{lemmedinterpolation}, qui, \textit{via} la majoration $\vert\mathrm{f}(m)\vert_{p_{0}}\le\vert\mathrm{f}\vert_{\mathsf{r}}$, permet de majorer $\alpha_{\sigma_{0}}(m(u_{0},u))\vert\mathrm{f}(m)\vert_{p_{0}}$ par \begin{equation*} e^{2nU_{0}(\log C_{0}+(\log\mathfrak{e})^{-1})}\max{\left\{e^{(17n^{2}-C_{0})U_{0}},\,(40nC_{0}^{2})^{2nC_{0}U_{0}}e^{(16n^{2}-C_{0}^{3/2})U_{0}}\right\}}\Vert s\Vert_{\overline{E}_{\alpha},\sigma_{0}},\end{equation*}quantit{\'e} plus petite que $e^{-(3/4)C_{0}U_{0}}\Vert s\Vert_{\overline{E}_{\alpha},\sigma_{0}}$. En effet, dans le maximum, la premi{\`e}re quantit{\'e} $e^{17n^{2}-C_{0}}$ l'emporte sur la seconde car $$C_{0}\left(2n\log(40 n)+4n\log C_{0}+1\right)\le C_{0}^{3/2}\le C_{0}^{3/2}+n^{2}$$ gr{\^{a}}ce au choix de $C_{0}=(4n)^{10n}$. Gr{\^{a}}ce {\`a} $\mathfrak{e}\ge e$, on a alors $$2n(\log C_{0}+(\log\mathfrak{e})^{-1})+17n^{2}-C_{0}\le -3C_{0}/4.$$En utilisant une nouvelle fois le lemme~\ref{lemmecomparaison} pour passer de $\mathrm{f}(m)$ {\`a} $\frac{1}{\tau !}\mathrm{D}_{w}^{\tau}F_{s,\sigma_{0}}(m(u_{0},u))$, on trouve \begin{equation*}\alpha_{\sigma_{0}}(m(u_{0},u))\left\vert\frac{1}{\tau !}\mathrm{D}_{w}^{\tau}F_{s,\sigma_{0}}(m(u_{0},u))\right\vert_{p_{0}}\le 2e^{-(3/4)C_{0}U_{0}}\Vert s\Vert_{\overline{E}_{\alpha},\sigma_{0}}.\end{equation*}Cette estimation est valide pour tout $\tau\in\mathbf{N}^{n}$ de longueur $\ell$ et, \textit{a priori}, pour $(m,\tau)\not\in\Upsilon$. Mais l'in{\'e}galit{\'e}~\eqref{ineqmajconstruction} montre qu'elle reste vraie pour $(m,\tau)\in\Upsilon$ puisque $\alpha^{-1}\le\exp{\{-C_{0}^{3/2}U_{0}\}}$. La d{\'e}finition~\ref{definitionjet} du jet de $s$ (voir aussi~\eqref{majorationdelanormedujet}) conduit alors au r{\'e}sultat suivant.
\begin{prop}\label{estimationjetvo}On a $\Vert\jet_{W}^{\ell}s(m\mathsf{p})\Vert_{\overline{\mathfrak{Jet}},\sigma}\le\exp{\{-C_{0}U_{0}/2\}}\Vert s\Vert_{\overline{E}_{\alpha},\sigma}$ pour tout plongement $\sigma$ conjugu{\'e} {\`a} $\sigma_{0}$.
\end{prop}
\subsubsection{Cas p{\'e}riodique}\label{sectioncasperiodique} Rappelons tout d'abord qu'il a {\'e}t{\'e} vu au \S~\ref{subseccomplements} que le cas p{\'e}riodique n'est possible que lorsque $\sigma_{0}$ est complexe. Soit $\tau=(\tau_{1},\ldots,\tau_{n})\in\mathbf{N}^{n}$ de longueur $\ell$ et $\tau':=(\tau_{1},\ldots,\tau_{n-1},0)$. Soit $\mathrm{f}:\mathbf{C}\to\mathbf{C}$ l'application analytique d{\'e}finie par~: $$\forall\,z\in\mathbf{C},\quad\mathrm{f}(z)=\frac{1}{\tau'!}\mathrm{D}_{w}^{\tau'}F_{s,\sigma_{0}}(z(\lambda(u),u)).$$Afin d'{\'e}valuer $\jet_{W}^{\ell}s(m\mathsf{p})$, nous allons donner une majoration fine de la d{\'e}riv{\'e}e $\mathrm{f}^{(\tau_{n})}(m)/\tau_{n}!$ au moyen du lemme d'interpolation~\ref{lemmedinterpolation}. Observons tout d'abord que \begin{equation}\label{ineqcauchy}\left\vert\frac{1}{\tau_{n}!}\mathrm{f}^{(\tau_{n})}(m)\right\vert\le\vert\mathrm{f}\vert_{1+m}\end{equation}gr{\^{a}}ce {\`a} l'in{\'e}galit{\'e} de Cauchy usuelle pour les fonctions holomorphes complexes. Le terme $\vert\mathrm{f}\vert_{1+m}$ peut {\^{e}}tre major{\'e} au moyen du lemme d'interpolation~\ref{lemmedinterpolation} avec les param{\`e}tres suivants~: $S_{1}:=(n+t)S$, $\mathsf{r}:=2(n+t)S$, $\mathsf{R}:=2\mathsf{r}\mathfrak{e}$, $T_{1}:=T_{0}$ et $f:=\mathrm{f}$. Dans ce lemme d'interpolation intervient $\vert\mathrm{f}\vert_{\mathsf{R}}$ que l'on peut estimer en suivant la d{\'e}monstration de la proposition~\ref{propositiondefr}, dont les calculs ont {\'e}t{\'e} faits de telle sorte qu'ils conviennent encore ici. Par ailleurs, gr{\^a}ce {\`a} la formule~\eqref{deriveedef}, au lemme de comparaison~\ref{lemmecomparaison}, et par construction de $s$, les d{\'e}riv{\'e}es de $\mathrm{f}$ jusqu'{\`a} l'ordre $T_{0}$ en les nombres $z\in\{0,\ldots,(n+t)S\}$ sont petites comme dans l'estimation~\eqref{petitessealpha}. En r{\'e}alit{\'e} une fois que l'on a observ{\'e} qu'avec les choix de param{\`e}tres ci-dessus le produit $T_{1}S_{1}$ {\'e}tait presque le m{\^e}me que dans le cas pr{\'e}c{\'e}dent, tout se passe {\`a} l'identique, si ce n'est que l'on doit se servir de~\eqref{ineqcauchy}. Un dernier passage par le lemme de comparaison~\ref{lemmecomparaison} permet alors d'obtenir \textit{in fine}~: \begin{equation}\label{ineqcasp}\alpha_{\sigma_{0}}(m(u_{0},u))\left\vert\frac{1}{\tau!}\mathrm{D}_{w}^{\tau'}\mathrm{D}_{(\lambda(u),u)}^{\tau_{n}}F_{s,\sigma_{0}}(m(u_{0},u))\right\vert\le 2\exp{\{-(3/4)C_{0}U_{0}\}}\Vert s\Vert_{\overline{E}_{\alpha},\sigma_{0}}\end{equation}(majoration valide pour tout $\tau=(\tau',\tau_{n})$ de longueur $\ell$). M{\^e}me si nous n'obtenons pas une estimation de $\frac{1}{\tau!}\mathrm{D}_{w}^{\tau}F_{s,\sigma_{0}}(m(u_{0},u))$, ceci n'a pas d'importance car ce qui nous int{\'e}resse est le jet d'ordre $\ell$ de $s$, qui lui ne d{\'e}pend pas du choix de la base consid{\'e}r{\'e}e pour le d{\'e}finir (voir le commentaire qui suit la d{\'e}finition~\ref{definitionjet}). Choisissons la base $\mathsf{w}:=(w_{1},\ldots,w_{n-1},(\lambda(u),u))$ dans la d{\'e}finition~\ref{definitionjet}. Si $x\in W\otimes_{\sigma_{0}}\mathbf{C}$ a pour coordonn{\'e}es $(x_{1},\ldots,x_{n})$ dans la base $w$ alors $$\forall\,i\in\{1,\ldots,n-1\},\quad w_{i}^{\mathsf{v}}(x)=x_{i}-x_{n}\frac{\mathfrak{u}_{i}}{\mathfrak{u}_{n}}\quad\text{et}\quad(\lambda(u),u)^{\mathsf{v}}(x)=\frac{x_{n}}{\mathfrak{u}_{n}}\cdotp$$Par cons{\'e}quent on a$$\forall\,i\in\{1,\ldots,n-1\},\quad\Vert w_{i}^{\mathsf{v}}\Vert_{\overline{t_{G}^{\mathsf{v}}},\sigma_{0}}=\left(1+\left(\frac{\vert\mathfrak{u}_{i}\vert}{\vert\mathfrak{u}_{n}\vert}\right)^{2}\right)^{1/2}\quad\text{et}\quad\Vert(\lambda(u),u)^{\mathsf{v}}\Vert_{\overline{t_{G}^{\mathsf{v}}},\sigma_{0}}=\frac{1}{\vert\mathfrak{u}_{n}\vert}\cdotp$$De ce fait et de la majoration~\eqref{ineqcasp}, l'on d{\'e}duit \begin{equation}\label{estimationperiodique}\Vert\jet_{W}^{\ell}s(m\mathsf{p})\Vert_{\overline{\mathfrak{Jet}},\sigma_{0}}\le 2\left(\frac{\max{\{1,\Vert(\lambda(u),u)\Vert_{\overline{t_{G}},\sigma_{0}}\}}}{\vert\mathfrak{u}_{n}\vert}\right)^{\ell}\exp{\{-(3/4)C_{0}U_{0}\}}\Vert s\Vert_{\overline{E}_{\alpha},\sigma_{0}}.\end{equation}L'hypoth{\`e}se~\ref{hypocru} et l'in{\'e}galit{\'e} triangulaire donnent la majoration $$\Vert(\lambda(u),u)\Vert_{\overline{t_{G}},\sigma_{0}}\le1+\Vert(u_{0},u)\Vert_{\overline{t_{G}},\sigma_{0}}\le 1+\Vert u_{0}\Vert_{\overline{k^{n}/W_{0}},\sigma_{0}}+\sum_{j=1}^{n}{\vert u_{j}\vert}.$$Avec la r{\'e}duction (iii) du \S~\ref{paragraphereductions}, on a $\Vert u_{0}\Vert_{\overline{k^{n}/W_{0}},\sigma_{0}}\le 1$ et, comme $$\vert u_{j}\vert\le (D\log a_{j})/\mathfrak{e}\le(\log a_{j})\exp{\{(\mathfrak{a}\log\mathfrak{e})/D\}},$$ on en d{\'e}duit $$\Vert(\lambda(u),u)\Vert_{\overline{t_{G}},\sigma_{0}}\le\left(1+\log\prod_{j=1}^{n}{a_{j}}\right)\exp{\{(\mathfrak{a}\log\mathfrak{e})/D\}}\le\exp{\{(2\mathfrak{a}\log\mathfrak{e})/D\}}.$$De plus, en vertu des propositions~\ref{propperiodique} et~\ref{propannexe}, on a $\vert\mathfrak{u}_{n}\vert^{-1}\le C_{0}\mathfrak{a}\binom{n+t}{n}\widetilde{D}_{0}$. La proposition~\ref{propparametres}, (iv), permet de majorer $T\log\widetilde{D}_{0}$ par $(10n\log C_{0})U_{0}$. On trouve alors $$2\left(\frac{\max{\{1,\Vert(\lambda(u),u)\Vert_{\overline{t_{G}},\sigma_{0}}\}}}{\vert\mathfrak{u}_{n}\vert}\right)^{\ell}\le\exp{\{(C_{0}/4)U_{0}\}}.$$Conjugu{\'e} {\`a} l'in{\'e}galit{\'e}~\eqref{estimationperiodique}, ce r{\'e}sultat montre que la proposition~\ref{estimationjetvo} reste vraie dans le cas p{\'e}riodique.
\subsection{Conclusion}\label{sectionconclusion}
 La fin de la d{\'e}monstration du th{\'e}or{\`e}me~\ref{theoremereduit} repose sur le lemme~\ref{liouville} qui, compte tenu de la proposition~\ref{proplemmedemultiplicites}, s'{\'e}crit ici sous la forme \begin{equation}\label{ineqcrucialeliouville} h(\jet_{W}^{\ell}s(m\mathsf{p}))\ge-\widehat{\mu}_{\mathrm{max}}\left(\overline{\mathfrak{Jet}}\right).\end{equation}Nous allons montrer que cette minoration n'est pas compatible avec la majoration de la hauteur du jet que l'on peut trouver par le biais de l'hypoth{\`e}se~\ref{hypocru}. En effet les propositions~\ref{propestimationarchi}, \ref{propestimationultra} et~\ref{estimationjetvo} entra{\^\i}nent \begin{equation*} h(\jet_{W}^{\ell}s(m\mathsf{p}))\le \frac{nU_{0}}{D}+T\log (4D_{0})-\frac{C_{0}U_{0}[k_{\sigma_{0}}:\mathbf{Q}_{p_{0}}]}{2D}+h_{\overline{E}_{\alpha}}(s). 
\end{equation*}La proposition~\ref{propparametres}, (iv), qui permet de majorer $T\log(4D_{0})$ (en se rappelant que $D_{0}\le\widetilde{D}_{0}$ car $x\le 1$) et la construction de $s$ (proposition~\ref{propconstructiondes}) conduisent {\`a} \begin{equation}\label{ineq:majfinaledujet}h(\jet_{W}^{\ell}s(m\mathsf{p}))\le\left(n+10n\log C_{0}-\frac{C_{0}}{2}+(4n)^{4n}C_{0}^{1/2}\right)\frac{[k_{\sigma_{0}}:\mathbf{Q}_{p_{0}}]U_{0}}{D}\cdotp\end{equation}Par ailleurs la pente maximale de $\overline{\mathfrak{Jet}}$ a {\'e}t{\'e} major{\'e}e par  $7n^{2}U_{0}/D$ (proposition~\ref{propositionpentejet}). L'on s'aper{\c{c}}oit alors que le choix de $C_{0}$ met en contradiction~\eqref{ineqcrucialeliouville} et~\eqref{ineq:majfinaledujet}. L'hypoth{\`e}se~\ref{hypocru} est donc fausse et $\log\Vert u_{0}-\lambda(u)\Vert_{\overline{t_{G}},\sigma_{0}}\ge-(C_{0}^{3/2}+4n)U_{0}$. La proposition~\ref{propcomparaison} entra{\^{i}}ne $$\log\max_{1\le i\le t}{\vert\Lambda_{i}\vert_{p_{0}}}\ge-n^{3}\log b-\log\sqrt{n}-(C_{0}^{3/2}+4n)U_{0}.$$Pour conclure et obtenir le th{\'e}or{\`e}me~\ref{theoremereduit}, on observe que \begin{equation}\label{majorationdeuo}U_{0}\le C_{0}^{3n+2}(C_{0}^{3}+2)\mathfrak{a}^{1/t}(\log b+\mathfrak{a}\log\mathfrak{e})\prod_{j=1}^{n}{\left(1+\frac{D\log a_{j}}{\log\mathfrak{e}}\right)^{1/t}}.\end{equation}En effet, dans la formule~\eqref{formuleuo} d{\'e}finissant $U_{0}$, on majore $$\left(\frac{C_{0}^{3n-1}t!\card\Sigma_{\mathsf{p}}(S)}{(n+t)!}\right)^{1/t}\le C_{0}^{3n+2}\mathfrak{a}^{1/t}.$$De plus, la quantit{\'e} $\mathfrak{y}:=\log b+D\log S+S^{y}\log\mathfrak{e}$ v{\'e}rifie \begin{equation*}\begin{split}\mathfrak{y}&\le\log b+3D\log C_{0}+D\log\mathfrak{a}+C_{0}^{3}\mathfrak{a}^{y}\log\mathfrak{e}\\ &\le (1+3\log C_{0})\log b+(C_{0}^{3}+2)\mathfrak{a}\log\mathfrak{e}\quad\text{(proposition~\ref{propparametres}, (iii))}\\ & \le (C_{0}^{3}+2)(\log b+\mathfrak{a}\log\mathfrak{e}),\end{split}\end{equation*}ce qui donne~\eqref{majorationdeuo}. Le th{\'e}or{\`e}me~\ref{theoremereduit} d{\'e}coule alors de $$C_{0}^{3n+2}(C_{0}^{3/2}+4n)(C_{0}^{3}+2)+n^{3}\le C_{0}^{10n}=(4n)^{90n^{2}}$$(on a major{\'e} $n^{3}\log b+\log\sqrt{n}$ par $n^{3}\log(b\mathfrak{e})$). Si $t=1$ et $\beta_{1,0}\ne 0$, on a $y=0$. La majoration de $\mathfrak{y}$ ci-dessus est remplac{\'e}e par $$\mathfrak{y}\le(1+3\log C_{0})(\log b+\log\mathfrak{e}+D\log\mathfrak{a}),$$qui se r{\'e}percute dans~\eqref{majorationdeuo}. En majorant $1+3\log C_{0}$ par $C_{0}^{3}+2$ on obtient alors le th{\'e}or{\`e}me~\ref{theoremereduit} dans ce cas.

\backmatter

\bibliographystyle{plain}

\end{document}